\documentclass[11pt]{amsart}
\usepackage{amsmath,amssymb, amsfonts}
\usepackage[all]{xy}
\usepackage{verbatim}
\usepackage{enumerate}

\newtheorem{thm}{Theorem}[section]
\newtheorem{cor}[thm]{Corollary}
\newtheorem{example}[thm]{Example}
\newtheorem{lem}[thm]{Lemma}
\newtheorem{obs}[thm]{Observation}

\newtheorem{prop}[thm]{Proposition}
\theoremstyle{definition}
\newtheorem{defn}[thm]{Definition}

\newtheorem{rem}[thm]{Remark}
\numberwithin{equation}{section}

\DeclareFontFamily{U}{rsf}{} \DeclareFontShape{U}{rsf}{m}{n}{
  <5> <6> rsfs5 <7> <8> <9> rsfs7 <10->  rsfs10}{}
\DeclareMathAlphabet{\mathscr}{U}{rsf}{m}{n}

\renewcommand{\imath}{\sqrt{-1}}

\DeclareMathOperator{\Hom}{Hom}
\DeclareMathOperator{\colim}{colim}
\DeclareMathOperator{\spec}{spec}
\DeclareMathOperator{\im}{im}
\DeclareMathOperator{\coim}{coim}
\DeclareMathOperator{\coker}{coker}

\DeclareMathOperator{\eq}{eq}
\DeclareMathOperator{\Sym}{Sym}
\DeclareMathOperator{\Frac}{Frac}
\DeclareMathOperator{\straightS}{S}
\DeclareMathOperator{\sF}{F}
\DeclareMathOperator{\sG}{G}
\DeclareMathOperator{\sR}{R}
\DeclareMathOperator{\sL}{L}
\DeclareMathOperator{\straightF}{F}
\DeclareMathOperator{\sK}{K}

\DeclareMathOperator{\sD}{D}
\DeclareMathOperator{\sP}{P}
\DeclareMathOperator{\sI}{I}
\DeclareMathOperator{\su}{u}
\DeclareMathOperator{\straightK}{K}
\DeclareMathOperator{\straightN}{N}
\DeclareMathOperator{\straightD}{D}
\DeclareMathOperator{\straightC}{C}
\DeclareMathOperator{\straighth}{h}
\DeclareMathOperator{\sS}{S}

\DeclareMathOperator{\LH}{LH}

\newcommand{\uHom}{\underline{\Hom}}
\newcommand{\ootimes}{\overline{\otimes}}
\newcommand{\wotimes}{\widehat{\otimes}}

\newcommand{\mA}{\mathcal{A}}
\newcommand{\mB}{\mathcal{B}}
\newcommand{\mC}{\mathcal{C}}

\newcommand{\tC}{\text{\bfseries\sf{C}}}
\newcommand{\tD}{\text{\bfseries\sf{D}}}

\newcommand{\tL}{\text{\bfseries\sf{L}}}
\newcommand{\tT}{\text{\bfseries\sf{T}}}
\newcommand{\tS}{\text{\bfseries\sf{S}}}

\newcommand{\tP}{\text{\bfseries\sf{P}}}
\newcommand{\tI}{\text{\bfseries\sf{I}}}

\newcommand{\mM}{\mathcal{M}}
\newcommand{\mN}{\mathcal{N}}
\newcommand{\mE}{\mathcal{E}}
\newcommand{\mF}{\mathcal{F}}

\newcommand{\mG}{\mathcal{G}}

\newcommand{\mV}{\mathcal{V}}

\newcommand{\ttMod}{{\text{\bfseries\sf{Mod}}}}

\newcommand{\ttSNrm}{{\text{\bfseries\sf{SNrm}}}}

\newcommand{\ttBan}{{\text{\bfseries\sf{Ban}}}}

\newcommand{\Mod}{{\text{\bfseries\sf{Mod}}}}
\newcommand{\ttVect}{{\text{\bfseries\sf{Vect}}}}
\newcommand{\ttComm}{{\text{\bfseries\sf{Comm}}}}
\newcommand{\ttGpd}{{\text{\bfseries\sf{Gpd}}}}
\newcommand{\ttCat}{{\text{\bfseries\sf{Cat}}}}
\newcommand{\ttL}{{\tL}}
\newcommand{\ttC}{{\tC}}
\newcommand{\ttD}{{\tD}}
\newcommand{\ttP}{{\tP}}
\newcommand{\ttI}{{\tI}}

\newcommand{\ttT}{{\tT}}
\newcommand{\ttS}{{\tS}}
\newcommand{\ttsC}{{\text{\bfseries\sf{sC}}}}
\newcommand{\ttA}{{\text{\bfseries\sf{A}}}}
\newcommand{\ttE}{{\text{\bfseries\sf{E}}}}
\newcommand{\ttSch}{{\text{\bfseries\sf{Sch}}}}

\newcommand{\ttAn}{{\text{\bfseries\sf{An}}}}
\newcommand{\ttPr}{{\text{\bfseries\sf{Pr}}}}
\newcommand{\ttSPr}{{\text{\bfseries\sf{SPr}}}}
\newcommand{\ttAff}{{\text{\bfseries\sf{Aff}}}}
\newcommand{\ttAfnd}{{\text{\bfseries\sf{Afnd}}}}
\newcommand{\ttSh}{{\text{\bfseries\sf{Sh}}}}
\newcommand{\ttSt}{{\text{\bfseries\sf{St}}}}
\newcommand{\ttSet}{{\text{\bfseries\sf{Set}}}}
\newcommand{\ttTop}{{\text{\bfseries\sf{Top}}}}

\newcommand{\ttInd}{{\text{\bfseries\sf{Ind}}}}
\newcommand{\ttsInd}{{\text{\bfseries\sf{sInd}}}}

\newcommand{\hZar}{{\text{\bfseries\sf{hZar}}}}
\newcommand{\hZarAff}{{\text{\bfseries\sf{hZarAff}}}}
\numberwithin{equation}{section}

\begin{document}

\title[]{Non-Archimedean analytic geometry as relative algebraic geometry}
\author{Oren Ben-Bassat, Kobi Kremnizer}\thanks{}%
\address{Mathematical Institute,
University of Oxford,
Andrew Wiles Building,
Radcliffe Observatory Quarter,
Woodstock Road,
Oxford,
OX2 6GG, England}
\email{Oren.Ben-Bassat@maths.ox.ac.uk, Yakov.Kremnitzer@maths.ox.ac.uk}%
\dedicatory{}
\thanks{We would like to thank  F. Baldassarri, F. Bambozzi, V. Berkovich, J. Block, D. Joyce, T. Pantev, F. Paugam, P. Schapira, B. To\"{e}n, M. Vaqui\'{e} and G. Vezzosi for interesting conversations.  Thank you to Ivan Fesenko for his input and encouragement. We would especially like to thank Konstantin Ardakov, Antoine Ducros and Michael Temkin for lots of help and discussions relating to Berkovich analytic geometry. The first author acknowledges the support of the European Commission under the Marie Curie Programme for the IEF grant which enabled this research to take place. The contents of this article reflect the views of the two authors and not the views of the European Commission.}

\begin{abstract}We show that non-Archimedean analytic geometry can be viewed as relative algebraic geometry in the sense of To\"{e}n--Vaqui\'{e}--Vezzosi over the category of non-Archimedean Banach spaces. For any closed symmetric monoidal quasi-abelian category we define a topology on certain subcategories of the  category of (relative) affine schemes. In the case that the monoidal category is the category of abelian groups, the topology reduces to the ordinary Zariski topology. By examining this topology in the case that the monoidal category is the category of Banach spaces we recover the G-topology or the topology of admissible subsets on affinoids which is used in analytic geometry. This gives a functor of points approach to non-Archimedean analytic geometry. We demonstrate that the category of Berkovich analytic spaces (and also rigid analytic spaces) embeds fully faithfully into the category of (relative) schemes in our version of relative algebraic geometry. We define a notion of quasi-coherent sheaf on analytic spaces which we use to characterize surjectivity of covers. Along the way, we use heavily the homological algebra in quasi-abelian categories developed by Schneiders.

\noindent 
\end{abstract}
\maketitle 
\tableofcontents
\section{Introduction} 
Berkovich analytic spaces and rigid analytic spaces \cite{Ber1990}, \cite{Ber2009} have the advantage that both schemes of finite type over a field and formal completions of such schemes along their closed subschemes can be thought of as living in the same category of analytic spaces. Punctured tubular neighborhoods of algebraic subvarieties inside ambient varieties over any field can be defined \cite{BeTe}, \cite{Th} as Berkovich analytic spaces.  In this article, we consider non-Archimedean analytic spaces from the perspective of algebraic geometry relative to the closed symmetric monoidal categories of Banach spaces. This language is very universal and provides a place to compare different geometries (rigid analytic spaces, Berkovich spaces and others). To\"{e}n and Vaqui\'{e} introduced in \cite{TV} algebraic geometry relative to any closed symmetric monoidal category. This idea had also been pursued by Hakim \cite{Ha} and Deligne \cite{Del}. We examine one of To\"{e}n and Vezzosi's \cite{TVe5} topologies (which we call the homotopy Zariski topology) restricted to the opposite category of affinoid algebras over a non-Archimedean field and show that Berkovich and rigid analytic geometry embed fully and faithfully into the resulting category of schemes. 
We use the framework of To\"{e}n, Vaqui\'{e} and Vezzosi to suggest a new approach to analytic geometry which will extend in a natural way to the setting of higher and derived analytic stacks over Banach rings. This of course produces the difficult task of comparing the abstract definitions to the existing ones.  In order to make this article accessible to a broad range of mathematicians, we have included as many details as was possible. In forthcoming work \cite{BaBe}, \cite{BaBeKr}, we incorporate complex analytic geometry and give a more general approach that covers dagger analytic geometry and Stein geometry over a general valuation field (Archimedean or non-Archimedean). Some of the properties of the localization maps which we look at were proven from a complex analytic or differential geometric point of view in \cite{Block, Pir, Tay}. Some similar work to ours was done \cite{Mac} by A. Macpherson who uses an abstract, categorical notion of localizations instead of the condition on derived categories which we use. See also \cite{Yu} and the references therein for the relationship with mirror symmetry. 

Using the homological algebra from the work of Schneiders \cite{SchneidersQA}, we give a new interpretation of the homotopy Zariski topology from \cite{TVe5} which is suited to deal with quasi-abelian categories. Two of the theorems which we prove (Theorems \ref{thm:TensLocs} and \ref{thm:localForm}) show that that the finitely presented morphisms of affinoid algebras $\mA\to \mB$ which are \emph{homotopy epimorphisms} (meaning that $\mB\wotimes^{\mathbb{L}}_{\mA}\mB \cong \mB$) correspond geometrically to the affinoid subdomains (affinoids which are the unions of rational domains). In terms of modules, this says that a morphism of affinoids $\mM(\mB) \to \mM(\mA)$ is a domain embedding if and only if the corresponding morphism $D^{-}(\mB)\to D^{-}(\mA)$ is fully faithful. Our categorical characterization of affinoid domains answers an open question poised by Soibelman in section 1.4 of \cite{So}. He identifies this question as being a key step in the development of non-commutative analytic geometry. The other main feature of the G-topology is that the covers in this topology have finite subcovers $\pi:\coprod_{i} \mM(\mA_{V_i}) \to \mM(\mA)$ of an affinoid by affinoid subdomains which are surjective. In order to understand this surjectivity via modules, we introduce the notion of a RR-quasicoherent module for a commutative monoid $A$ relative to a symmetric monoidal quasi-abelian category in Definition \ref{defn:RRqc}. We call this category of modules $\ttMod^{RR}(A)$.   We prove in Lemmas \ref{lem:ConsImpliesSur} and \ref{lem:CoversRconservative} that $\pi:\coprod_{i} \mM(\mA_{V_i}) \to \mM(\mA)$ is surjective if and only if a morphism $f:M \to N$ of RR-quasicoherent modules is an isomorphism if and only $\pi^{*}f: \pi^{*}M \to \pi^{*}N$ is an isomorphism. These modules satisfy a version of Tate acyclicity (see Remark \ref{rem:TateAcylicity}) for the version we use of the Tate complex which uses only completed tensor products. More work on the category of quasi-coherent sheaves in analytic geometry will appear in \cite{BaBeKr, BeKr2}.

The condition we use on morphisms of algebras: homotopy epimorphism, is the key idea behind the holomorphic functional calculus studied by \cite{Tay} by Taylor in complex analytic geometry, who called these morphisms absolute localizations. Pirkovskii in \cite{Pir} shows that open embeddings of complex Stein manifolds satisfy the same abstract condition. 


 In \cite{BeKr2} we construct monoidal model structures on categories of simplicial objects and on categories of complexes (negatively graded or unbounded) in a quasi-abelian category $\ttC$ satisfying certain extra conditions. We prove in \cite{BeKr2} that these structures form HAG contexts in the sense of \cite{TVe3}.  A HAG context is essentially a monoidal model category $M$ satisfying a long list of axioms that make it easy to work with. The model structure is induced from the model structure on simplicial sets. For instance, a morphism $X_{\bullet} \to Y_{\bullet}$ in the category $M=\ttsC$ of simplicial objects in $\ttC$  is a fibration (resp. weak equivalence) if for all projectives $P$ in $\ttC$ (defined in \ref{defn:Projective}) that $\Hom(P, X_{\bullet}) \to \Hom(P, Y_{\bullet})$ is a fibration (resp. weak equivalence) of simplicial sets. In particular, our main application is when the quasi-abelian category $\ttC=\ttInd(\ttBan_{R})$ is the category of Ind-Banach spaces over a Banach ring. In \cite{TVe3} derived and homotopy algebraic geometry is developed relative to a HAG context.  Derived algebraic geometry (the case where $\ttC$ is the abelian category of abelian groups) is only one special case of their work. Our goal is to show that derived analytic geometry is another. A key feature in their work is definitions of (Grothendieck) topologies on the categories of affine schemes relative to $M$. Affine schemes are defined to be opposite category to $\ttComm(M)$, the category of commutative monoids relative to $M$. The main way to describe these topologies is to explain which morphisms should appear in the cover, and describe when such a collection of morphisms is a cover. Both of these can be explained in terms of geometrically induced functors on the homotopy categories of modules. One topology that can be defined from their point of view says 
\begin{defn}\label{defn:TVtopologyIntro}
 A cover of $\spec(A) \in \ttComm(M)^{op}$ is a collection of morphisms $\spec(B_i) \to \spec(A)$ in $\ttComm(M)^{op}$ for $i \in I$ such that for some finite set $J \subset I$ we have 
\begin{enumerate}
\item the push forward functor from the derived category of modules on $B_{i}$ to that on $A$ is fully faithful for all $i \in J$
\item a morphism in the derived category of modules on $A$ is an isomorphism if and only if it becomes an isomorphism when pulled back in the derived sense to a morphism in the derived category of modules on $B_i$ for all $i \in J$.
\end{enumerate} 
For more details see Definition 1.2.6.1 (3) of \cite{TVe3} for the first item and Definition 1.2.5.1 of \cite{TVe3} for the second.
\end{defn}
In our type of derived analytic geometry, the category of affine schemes over a valuation field is opposite to the category of commutative monoids relative to $M=\ttsInd(\ttBan_{k})$. Theorem \ref{thm:BigRestriction} shows that on the subcategory opposite to affinoid algebras, the topology in Definition \ref{defn:TVtopologyIntro} restricts to the weak G-topology on affinoids. This theorem relies on \cite{BeKr2}.

While condition (1) of Definition \ref{defn:TVtopologyIntro} restricts by our results to something completely classical and well understood, Condition (2) of Definition \ref{defn:TVtopologyIntro} seems difficult to check on modules. Therefore, one might try to use a similar condition on the (underived) quasi-abelian category of modules using the underived pullback. However, this does not correspond to the surjectivity condition on covers, forcing us to reconsider which modules we care about and introduce the RR-quasicoherent modules.  For this reason, we defined the formal homotopy Zariski topology on a subcategory $\ttA$ of the category of affine schemes ($\ttComm(\ttC)^{op}$) of a closed symmetric monoidal quasi-abelian category $\ttC$ where $\ttA$ and $\ttC$ satisfy a few simple conditions. The covers are a collection of morphisms $A\to B_{i}$ of in $\ttA$ for $i \in I$ such that for some finite set $J \subset I$
\begin{enumerate}
\item $A \to B_{i}$  are homotopy epimorphisms in $\ttComm(\ttC)$ for all $i \in J$
\item a morphism $f:M \to N$ in $\ttMod^{RR}(A)$ is an isomorphism if and only if the induced morphisms $M\ootimes_{A}B_i \to N\ootimes_{A}B_i$ are for all $i \in J$.
\end{enumerate}
In the special case that $\ttC=\ttBan_{k}$ and $\ttA$ is the category of affinoids our Theorem \ref{thm:covers} says that this topology restricts to the weak G-topology on affinoids. Our work, put together with the work of To\"{e}n and Vezzosi establishes foundations for derived (and homotopy) analytic geometry and the study of higher and derived stacks over general Banach rings. There is another approach due to F. Paugam. In the complex analytic case, J. Lurie and M. Porta have a different approach to derived analytic geometry which should be compared to ours.

\section{Algebras and modules in closed symmetric monoidal categories} We assume that the reader is comfortable with the idea of a symmetric monoidal category. In particular the bifunctor $\ootimes$ admits an adjoint  $\uHom :\ttC^{op} \times \ttC \to \ttC$ in the sense that there are natural isomorphisms
\[\Hom(U \ootimes V,W) \cong \Hom(U, \uHom(V,W))\] for any $U,V,W \in \ttC.$
Consider a category $\ttC$ which has all finite limits and colimits equipped with extra data $(\underline{\Hom}, \overline{\otimes}, \text{id})$ making it into a closed, symmetric monodial category. 
In such a category one always has the following natural isomorphisms for any $U,V,W \in \ttC$

\begin{enumerate}
\item \[\Hom(U,V) \cong \Hom(\text{id}, \uHom(U,V))\]
\item \[\uHom(U \ootimes V,W) \cong \uHom(U, \uHom(V,W))\]
which lifts by taking the $\Hom$ from $\text{id}$ the corresponding isomorphisms for $\Hom$ instead of $\uHom$
\item \[ U \cong \uHom(\text{id},U).\]
\end{enumerate}

Also, there are natural morphisms for any $T,U,V,W \in \ttC$

\begin{enumerate}
\item \[\uHom(U,U) \ootimes U \to U\]
\item \[\uHom(V,W) \ootimes \uHom(U,V) \to \uHom(U,W) \] satisfying associativity 
\item \[\uHom(T,U) \ootimes \uHom(V,W) \to \uHom(T \ootimes V, U \ootimes W)  \]
\item \[\text{id} \to \uHom(U,U)\]
\end{enumerate}

which are compatible with the corresponding morphisms for $\Hom$ instead of $\uHom$. 
\begin{rem}
The morphisms in the second item are adjoint to morphisms 
\[\uHom(V,W) \to \uHom( \uHom(U,V) ,\uHom(U,W))
\]
and
\[\uHom(U,V) \to \uHom( \uHom(V,W) ,\uHom(U,W))
\]
and hence maps
\[\Hom(V,W) \to \Hom( \uHom(U,V) ,\uHom(U,W))
\]
and
\[\Hom(U,V) \to \Hom( \uHom(V,W) ,\uHom(U,W)).
\]
\end{rem}
Such categories admit clear definitions of categories of commutative unital algebras objects in them. We denote these categories by $\ttComm(\ttC)$. The opposite category to $\ttComm(\ttC)$ we will be denoted $\ttAff(\ttC)= \ttComm(\ttC)^{op}$ and if $A \in \ttComm(\ttC)$ we use $\spec(A)$ to indicate this object in the opposite category.

For any fixed commutative unital algebra object $A \in \ttComm(\ttC)$ there is a category of unital left modules over such an algebra which we denote $\ttMod(A)$. See Chapter 1.3 of \cite{M} for details.  
Because we will be considering the unital notions only, we drop the word unital from our descriptions. This category is also a  closed symmetric monoidal category and has all finite limits and colimits. We write $\ootimes_{A}$ for the monoidal structure and $\Hom_{A}$ for the morphisms in this category.  Given a morphism $p:\spec(B) \to \spec(A)$, we sometimes write $p^{*}$ for the functor 
\[\ttMod(A) \to \ttMod(B)
\]
given by
\[
M \mapsto B \ootimes_{A} M. 
\]
\begin{defn}\label{defn:TensA}Let us consider $A \in \ttComm(\ttC)$ and $E,F \in \ttMod(A)$. Let $a_{E}:A \ootimes E \to E$ and $a_{F}:A \ootimes F \to F$ denote the action morphisms. 
The set of morphisms $\Hom_{A}(E,F)$ is defined as the limit of the diagram

\[
\xymatrix{\Hom(E,F) \ar@/^2pc/[rr]^{f \mapsto f \circ a_{E}} \ar[dr]_{h \mapsto \text{id}_{A} \otimes h} &   & \Hom(A \ootimes E,  F) \\
& \Hom(A \ootimes E, A \ootimes F)\ar[ur]_{g \mapsto a_{F} \circ g} & 
}.
\] 
In order to describe the tensor product, let us use $\sigma:A \ootimes F \to F \ootimes A$ to denote the symmetric structure. Let $E \ootimes_{A} F\in \ttMod(A)$ be the element of $\ttC$ given as the colimit of the diagram 
\begin{equation}\label{equation:TensDef}\xymatrix{E \ootimes A \ootimes F \ar@/^/[rr]^{\text{id}_{E} \otimes a_{F}} \ar@/_/[rr]_{(a_{E} \otimes \text{id}_{F})\circ (\sigma \otimes \text{id}_{F})} & & E \ootimes F}
\end{equation}
endowed with the obvious action of $A$.

\end{defn}
Notice that for any $E\in \Mod(A)$ and $F \in \ttC$ can define 
\[L_{E,F}:A\ootimes \uHom(E,F) \to  \uHom(E,F)
\]
as the composition 
\[A\ootimes \uHom(E,F) \to  \uHom(E,E) \ootimes \uHom(E,F)\to  \uHom(E,F)
\]
where the morphism $A \to \uHom(E,E)$ is adjoint to the action morphisms $A \ootimes E \to E.$ Similarly for any $E\in \ttC$ and $F \in \Mod(A)$, one can define a morphism 
\[R_{E,F}:A \ootimes \uHom(E,F)  \to  \uHom(E,F)
\]
as the composition 
\[ A \ootimes \uHom(E,F)  \to  \uHom(F,F) \ootimes \uHom(E,F) \to   \uHom(E,F) \ootimes \uHom(F,F)  \to \uHom(E,F).
\]
Both $L_{E,F}$ and $R_{E,F}$ endow  $\uHom(E,F)$ with the structure of an element of $\Mod(A)$ satisfying various natural properties and we use these structures without further comment.
\begin{lem}\label{lem:ALinHom}Suppose now that $\ttC$ is a closed, symmetric monoidal additive category with all finite limits and colimits and $A \in \ttComm(\ttC)$. Then $\Mod(A)$ is a closed symmetric monoidal category with all finite limits and colimits as well. These limits and colimits can be computed in $\ttC$.
\end{lem}
{\bf Proof.}
By tensoring the morphisms $\text{id} \to \uHom(A,A)$ with the identity of $\uHom(E,F)$ one can form 
\[\uHom(E,F) \to \uHom(A,A) \ootimes \uHom(E,F) \to \uHom(A\ootimes E,A\ootimes F). 
\] Now by using the functor $\uHom(A\ootimes E, \;)$ and the contravariant functor $\uHom(\;,F)$ applied to $A\ootimes F \to F$ and $A\ootimes E \to E$ respectively one can define $\uHom_{A}(E,F)$ with the same type of limit as in \ref{defn:TensA} replacing $\Hom$ with $\uHom$. Both $L_{E,F}$ and $R_{E,F}$ induce well defined (and equal) morphisms on  $\uHom_{A}(E,F)$ and give  $\uHom_{A}(E,F)$ the structure of an element of $\ttMod(A)$.
\ \hfill $\Box$

We have natural isomorphisms 
\begin{equation}\label{equation:ClosedPlus}\uHom(M \ootimes_{A} N, L) \cong \uHom(M, \uHom_{A}(N,L))
\end{equation}
which satisfy the relevant axioms of a closed, symmetric monoidal category.
It is easy to see that (as in) \cite{MeyerDerived} for any $E\in \Mod(A)$
\begin{lem}\label{lem:MeyerProperties}\[\uHom_{A}(A,E) \to E \]
\[f \mapsto f(1)
\]
is an isomorphism and so for any modules $M,N \in \Mod(A)$ and any $V \in \ttC$ we have 
\begin{enumerate}
\item $A \ootimes_{A} M \cong M$
\item $\uHom_{A}(A \ootimes V,E) \cong \uHom(V,E)$
\item $\uHom_{A}(E, \uHom(A,V)) \cong \uHom(E,V)$
\item $\Hom_{A}(A, \uHom_{A}(M,N)) \cong \Hom_{A}(M,N)$
\end{enumerate} 
\end{lem}

We will discuss in this article various closed symmetric monoidal categories whose objects are vector spaces over Archimedean or non-Archimedean fields equipped with extra structures. The unit object in these categories is the field itself and it is easy to check the associativity, commutativity, and unit constraints. It is also clear that all of our closed symmetric monoidal categories are $k$-linear additive categories, as are the categories $\ttMod(A)$ for $A \in \ttComm(\ttC)$. Because $\ttC$ has all finite limits and colimits, $\ttMod(A)$ does as well. 

\begin{defn}A functor $\straightF:\ttT\to \ttS$ is 
said to be conservative under that 
the condition that a morphism $f$ in 
$\ttT$ is an isomorphism if and only 
if $\straightF(f)$ is an isomorphism.
\end{defn}
Notice that for every $p:\spec (B) \to \spec (A)$ in $\ttAff(\ttC)$ the functor $p^{*}:\ttMod(A) \to \ttMod(B)$ has a conservative right adjoint given by considering a $B$ module as an $A$ module 
\[p_{*}:\ttMod(B) \to \ttMod(A).\]

\begin{defn}\label{Symmetric}

If $\ttC$ has countable coproducts, the forgetful functor $\straightF:\ttComm(\ttC)\to \ttC$ has a left adjoint $\straightS:\ttC\to \ttComm(\ttC)$, given on objects 
by 
\begin{equation*}
\straightS(V)=\coprod_{n\geq 0} \straightS^n(V)
\end{equation*}

where $\straightS^n(V)$ are the co-invariants for the symmetric group action on $V^{\ootimes n}$. If $A$ is a commutative algebra in $\ttC$ we get in a similar way a functor 
$\straightS_A:\ttMod(A) \to \ttComm(\ttMod(A))$ which is left adjoint to the forgetful functor. 
\end{defn}

\begin{defn}\label{defn:morphismSymAlg}
Let $\ttC$ be a closed symmetric monoidal category with countable coproducts. Let $E$ and $E'$ be two objects, and assume given a map of algebras $\alpha:\straightS(E')\to \straightS(E)$. Let us denote by $\straightS(E,E',\alpha)$ the quotient of $\straightS(E)$ by the ideal generated by the image of $\alpha$ restricted to the augmentation ideal. In the case that $A\in \ttComm(\ttC)$ and $E$ and $E'$ are in $\ttMod(A)$ we use the notation $\straightS_{A}(E,E',\alpha)$ to denote the same construction taken in the closed symmetric monoidal category $\ttMod(A)$ in place of $\ttC$.
An equivalent way of describing it is as the quotient of $\straightS(E)$ by the ideal generated by the image of the object $E'$ under the map $\alpha$.
\end{defn}
 
 \hfill $\Box$
\section{Algebraic geometry relative to closed symmetric monoidal categories}\label{TV}

In this section, we review several notions from the article \cite{TV} by To\"{e}n and Vaqui\'{e}. We review their definitions of a flat morphism, a Zariski open immersion, and the Zariski and fpqc topologies on the opposite of the category of algebra objects in a closed, symmetric monoidal category admitting all finite limits and colimits. We also include definitions of a formal Zariski open immersion and the formal Zariski topology. We also review the definitions of schemes and (higher) stacks in this setting.
\begin{defn}
The category of presheaves of sets on $\ttAff(\ttC)$ is the category whose objects are contravariant functors $\ttAff(\ttC) \to \ttSet$ and whose morphisms are natural transformations. It will be denoted $\ttPr(\ttAff(\ttC))$. Given a Grothendieck topology $T$ on the category $\ttAff(\ttC)$ we will often consider the full subcategories of sheaves of sets $\ttPr(\ttAff(\ttC)^{T}) \subset \ttPr(\ttAff(\ttC))$, where $\ttAff(\ttC)^{T}$ is the site consisting of the underlying category $\ttAff(\ttC)$ along with the Grothendieck topology $T$.
\end{defn}
\begin{defn}
A family of functors $\{\straightF_{i}:\ttC \to \tD_{i}\}_{i \in I}$ 
is said to be conservative if a morphism $f$ in $\ttC$ 
is an isomorphism if and only if the 
morphisms $\straightF_{i}(f)$ in $\tD_i$ are 
isomorphisms for all $i \in I$.
\end{defn}
Since we do not assume the existence of arbitrary limits and colimits, we need to define finite presentation in a non-standard way.
\begin{defn}\label{defn:AbstractFinitePres}Let $\ttC$ be a closed symmetric monoidal category with all finite limits and colimits. An object $W\in \ttC$ is called of finite presentation if for every filtered diagram of objects $V_{i} \in \ttC$ such that $\colim_{i}V_{i}$ exists, the natural morphism
\begin{equation}\label{eqn:AbstractFinitePres}\colim_{i}\Hom_{\ttC}(W,V_{i}) \to \Hom_{\ttC}(W,\colim_{i}V_{i})
\end{equation}
is an isomorphism.
\end{defn}
\begin{rem}Notice that any finite colimit of finitely 
presented objects is finitely presented.
\end{rem}
\begin{defn}\label{defn:FinitePres}
Let $\ttC$ be a closed symmetric monoidal category with all finite limits and colimits. A morphism $A\to B$ in $\ttComm(\ttC)$ is of finite presentation if for every filtered diagram of objects $A'_{i}\in A/\ttComm(\ttC)$ such that $\colim_{i}A_{i}'$ exists in $A/\ttComm(\ttC)$, the natural morphism
\begin{equation}\label{eqn:FinitePres}\colim_{i}\Hom_{A/\ttComm(\ttC)}(B,A_{i}') \to \Hom_{A/\ttComm(\ttC)}(B,\colim_{i}A_{i}')
\end{equation}
is an isomorphism.
\end{defn}
\begin{rem}Notice that a morphism $A\to B$ in $\ttComm(\ttC)$ is of finite presentation if and only if $B$ is of finite presentation with respect to the category $A/\ttComm(\ttC)$.
\end{rem}
\begin{rem}\label{rem:epi}Let $\ttC$ be a closed symmetric monoidal category with all finite limits and colimits. We will need to pay special attention to the {\it epimorphisms} in $\ttComm(\ttC)$: those morphisms $p:A \to B$ in $\ttComm(\ttC)$ such that for all $C \in \ttComm(\ttC)$, the induced map $\Hom_{\ttComm(\ttC)}(B,C) \to \Hom_{\ttComm(\ttC)}(A,C)$ is injective. These correspond to monomorphisms in $\ttAff(\ttC)$. Notice that $p:A \to B$ is an epimorphism if and only if the multiplication $B\ootimes_{A}B \to B$ is an isomorphism.
\end{rem}
\begin{defn}\label{defn:AbstractFlat}Let $\ttC$ be a closed symmetric monoidal category with all finite limits and colimits. An object $V\in \ttC$ is called flat when the functor $\ttC\to\ttC$ given by $W \mapsto V\ootimes W$ is exact (commutes with finite limits).
\end{defn}
\begin{rem}\label{rem:SummandOfFlat} Suppose that $\ttC$ be a closed symmetric monoidal category with all finite limits and colimits and that finite products and coproducts agree in $\ttC$. Then a finite coproduct of elements of $\ttC$ is flat if and only if each of them is individually flat.
\end{rem}
\begin{defn}\label{defn:TVf}Let $\ttC$ be a closed symmetric monoidal category with all finite limits and colimits. A morphism $p:A \to B$ in $\ttComm(\ttC)$ is flat when the morphism \[p^{*}:\ttMod(A) \to \ttMod(B)\] is exact (commutes with finite limits). This precisely says that $B$ is flat in $\ttComm(\ttMod(A))$.
\end{defn}
Say that $q:\spec(C) \to \spec(A)$ is arbitrary. Consider the Cartesian diagram 
\begin{equation}\label{basechange}\xymatrix{  \spec(C\otimes_{A}B) \ar[r]^{q'} \ar[d]_{p'} & \spec(B) \ar[d]^{p} \\
\spec(C) \ar[r]_{q}& \spec(A). 
}
\end{equation}

Using the notation of diagram (\ref{basechange}) There is a natural equivalence \cite{TV}
\begin{equation}\label{eqn:basechangeequation}p^{*}q_{*} \Longrightarrow q'_{*}p'^{*}
\end{equation}
called base change. 

 \begin{defn}A base change of a morphism $p:\spec(B) \to \spec(A)$ is the morphism $p'$ appearing in diagram (\ref{basechange}) for some $q$.
\end{defn}

\begin{lem}\label{InvBaseChange}
Let $\ttC$ be a closed symmetric monoidal category with all finite limits and colimits. Suppose that $p:A \to B$ is a flat morphism in $\ttComm(\ttC)$ then any base change $p'$ of $p$ is also a flat morphism in $\ttComm(\ttC)$. 
\end{lem}

{\bf Proof.}
By assumption, $p^{*}: \ttMod(A) \to \ttMod(B)$ is exact. Let $\ttL$ be a finite category and consider the categories $\ttMod(D)^{\ttL},$ the category of functors from $\ttL$ to $\ttMod(D)$. Consider the functor 
\[\lim_{\ttL,D}: \ttMod(D)^{\ttL} \to \ttMod(D)
\]
which takes a diagram in $\ttMod(D)$ to its limit.
Consider the following commutative diagram of natural transformations of functors $\ttMod(C)^{\ttL} \to \ttMod(B)$: 
\[
\xymatrix{q'_{*}\lim_{\ttL,C\otimes_{A}B} p'^{*} 
\ar@{=>}@/_2pc/[rrrr] 
& \ar@{=>}[l] \lim_{\ttL,B} q'_{*} p'^{*} & \ar@{=>}[l] \lim_{\ttL,B} p^{*}q_{*} \ar@{=>}[r] & p^{*}q_{*} \lim_{\ttL,C} \ar@{=>}[r] & q'_{*} p'^{*}\lim_{\ttL,C}.}
\]
\\
All the natural transformations except the curved one on the bottom are invertible. This follows from the base change natural equivalences and the fact that the pushforward functors are right adjoints and so commute with limits. Therefore, the natural transformation on the bottom is also invertible.
Since $q'_{*}$ is conservative, the natural transformation $\lim_{\ttL,C\otimes_{A}B} p'^{*} \Longrightarrow p'^{*}\lim_{\ttL,C}$ of functors $\ttMod(C)^{\ttL} \to \ttMod(C \otimes_{A} B)$ is also invertible and so $p'^{*}$ commutes with finite limits. Therefore, $p'$ is flat.
\ \hfill $\Box$
\begin{lem}\label{lem:left}Let $\ttC$ be a closed symmetric monoidal category with all finite limits and colimits. Let $p$ be a morphism in $\ttAff(\ttC)$ and let $p'$ be a base change of $p$.  If $p$ is a monomorphism then so is $p'$.  If $p$ is of finite presentation, so is $p'$.
\end{lem}
{\bf Proof.} Left to the reader.
\ \hfill $\Box$

Consider the following very slight modification of Proposition 2.4 and its proof from \cite{TV} (we only require finite limits and colimits and consider only a special case of their proposition).
\begin{lem}\label{lem:BaseChangeConserv}Let $\ttC$ be a closed symmetric monoidal category with all finite limits and colimits. Suppose that a family $\{p_{i}:X_{i}\to X\}$ in $\ttAff(\ttC)$ is such that the family $\{p^{*}_{i}:\ttMod(X) \to \ttMod(X_i)\}$ has a finite conservative subfamily. Then any pull-back family $\{p_{i}:X_{i}\times_{X} Y\to Y\}$ coming from a base change $Y\to X$ has the same property.
\end{lem}
{\bf Proof.} In order to show the base change property, consider $q:Y \to X$. Choose a finite set $J \subset I$ such that $\prod_{i \in J} p^{*}_{i}$ is conservative. Consider the functor $\prod_{i \in J} p^{'*}_{i}$ where $q'_{i}$, $p'_i$ and $p_i$ play the role of $q'$, $p'$ and $p$ in diagram (\ref{basechange}). In order to show it is conservative, its enough to show that $\prod_{i \in J} q'_{i*}p^{'*}_{i}$ is conservative but using equation (\ref{eqn:basechangeequation}) this is isomorphic to $(\prod_{i \in J} p^{*}_{i}) q_{*}$ which is conservative since $q_{*}$ is conservative.
\ \hfill $\Box$
\begin{prop}\label{TV1}
Let $\ttC$ be a closed symmetric monoidal category with all finite limits and colimits. Consider the families $\{p_{i}:X_{i}\to X\}_{i \in I}$ in $\ttAff(\ttC)$ such that the family $\{p^{*}_{i}:\ttMod(X) \to \ttMod(X_i)\}_{i \in I}$ has a finite conservative subfamily and that each $p_{i}^{*}$ is left exact. These families define a pretopology on $\ttAff(\ttC)$. 
\end{prop}
{\bf Proof.} In order to show the base change property, consider $q:Y \to X$ in $\ttAff(\ttC)$ and let  $q'_{i}$, $p'_i$ and $p_i$ play the role of $q'$, $p'$ and $p$ in diagram (\ref{basechange}). Lemma \ref{lem:BaseChangeConserv} implies that the family $\{p^{*}_{i}\}$ has a finite conservative subfamily. The fact that the $p_{i}^{*}$ are exact follows from Lemma \ref{InvBaseChange}.

\ \hfill $\Box$
\begin{defn}\label{defn:fpqc}
For any closed symmetric monoidal category $\ttC$ which has all finite limits and colimits, the topology coming from Proposition \ref{TV1} is called the fpqc topology on $\ttAff(\ttC)= \ttComm(\ttC)^{op}$. When equipped with this topology, we denote this category by $\ttAff(\ttC)^{fpqc}.$ The category of sheaves of sets is denoted $\ttSh(\ttAff(\ttC)^{fpqc}).$
\end{defn}
\begin{defn}\label{defn:fTVZ}
 The morphism $\spec(B) \to \spec(A)$ is called a formal Zariski open immersion if the corresponding morphism $A \to B$ in $\ttComm(\ttC)$ is a flat epimorphism (defined in Remark \ref{rem:epi} and Definition \ref{defn:TVf}).
\end{defn}
\begin{defn}\label{defn:TVZ}
 The morphism $\spec(B) \to \spec(A)$ is called a Zariski open immersion if the corresponding morphism $A \to B$ in $\ttComm(\ttC)$ is a flat epimorphism of finite presentation (defined in Remark \ref{rem:epi} and Definitions \ref{defn:FinitePres} and \ref{defn:TVf}).
\end{defn}
\begin{lem}\label{lem:ConsAgain}If a family $\{A\to B_{i}\}_{i\in I}$ is conservative and $A'$ is any $A$-algebra then the family $\{A'\to B_{i}\otimes_{A} A'\}_{i\in I}$ is conservative.
\end{lem}
{\bf Proof.}
This has already been shown in Proposition \ref{TV1}.
\ \hfill $\Box$

\begin{prop}\label{prop:formalZarPretop}
There is a pretopology whose covering families $\{A\to B_{i}\}_{i\in I}$ are those families where each $A\to B_{i}$ is a formal Zariski open immersion and the family $\{A\to B_{i}\}_{i\in I}$ has a finite conservative subfamily.
\end{prop}
{\bf Proof.}
This follows immediately from Lemmas \ref{InvBaseChange}, \ref{lem:ConsAgain} and \ref{lem:left}.
\ \hfill $\Box$
\begin{defn}\label{defn:fTV2} The formal Zariski topology on $\ttAff(\ttC)$ is the topology associated to the pretopology from Proposition \ref{prop:formalZarPretop}. When equipped with this topology, we denote the category by $\ttAff(\ttC)^{fZar}$. The category of sheaves of sets is denoted $\ttSh(\ttAff(\ttC)^{fZar}).$
\end{defn}
\begin{prop}\label{prop:ZarPretop}
There is a pretopology whose covering families $\{A\to B_{i}\}_{i\in I}$ are those families where each $A\to B_{i}$ is a Zariski open immersion and the family $\{A\to B_{i}\}_{i\in I}$ has a finite conservative subfamily.
\end{prop}
{\bf Proof.}
This follows immediately from Lemmas \ref{InvBaseChange}, \ref{lem:ConsAgain} and \ref{lem:left}.
\ \hfill $\Box$
\begin{defn}\label{defn:TV2} The Zariski topology on $\ttAff(\ttC)$ is the topology associated to the pretopology from Proposition \ref{prop:ZarPretop}. When equipped with this topology, we denote the category by $\ttAff(\ttC)^{Zar}$. The category of sheaves of sets is denoted $\ttSh(\ttAff(\ttC)^{Zar}).$
\end{defn}
\begin{defn}\label{defn:hRep} For any affine object $\spec(A),$ $A \in \ttComm(\tC),$ the presheaf of sets $\straighth_{A}$ is given by
\[\ttAff(\ttC) \to \ttSet
\]
\[\spec(B) \mapsto \Hom_{\ttComm(\tC)}(A,B).
\]
\end{defn}
Cor. 2.11 of \cite{TV} implies that 
\begin{prop}\label{TV3}
For any $A \in \ttComm(\tC)$, the preseheaf $\straighth_{A}$ is a sheaf for fpqc, the formal Zariski and the Zariski topologies. 
\end{prop}

\begin{defn}\label{defn: ZarOpenSch}\cite{TV} Let $X \in \ttAff(\tC)$ and $F\in \ttSh(\ttAff(\tC)^{Zar})$ be a subsheaf of $X$. Then $F$ is a called a Zariski open of $X$ if there is a family of Zariski opens $\{X_i \to X\}_{i \in I}$ such that $F$ is the image of the morphism of sheaves $\coprod_{i \in I} X_i \to X$. A morphism $F \to G$ in $\ttSh(\ttAff(\tC)^{Zar})$ is called is a Zariski open immersion if for every $X \in \ttAff(\tC)$ and every $X \to G$ the induced morphism $F \times_{G} X \to X$ is a monomorphism whose image is a Zariski open in $X$.
The category $\ttSch(\ttAff(\tC)^{Zar})$ of schemes is defined to be the full subcategory of $\ttSh(\ttAff(\tC)^{Zar})$ of sheaves $F$ such that there exists a family of $X_{i} \in \ttAff(\tC)$ for $i \in I$ and a morphism $p: \coprod_{i \in I} X_i \to F$ such that $p$ is an epimorphism of sheaves and for each $i$ the morphism $X_i \to F$ is a Zariski open.
\end{defn}

\begin{defn}\label{defn:GeneralSch} Suppose that $\ttA$ is a full subcategory of $\ttAff(\tC)$ and suppose that $\tau$ is a subcategory of $\ttA$ with all objects and such that all morphisms in $\tau$ are monomorphisms and so that the base change of a morphism in $\tau$ by an arbitrary morphism of $\ttA$ is in $\tau$. Say that $T$ is a pre-topology whose covers consist of families of morphisms where each morphism in the cover belongs to $\tau$.  Let $X \in \ttA$ and $F\in \ttSh(\ttA^{T})$ be a subsheaf of $X$. Then $F$ is a called a $\tau$-open if there is a family of morphisms in $\tau$ written  $\{X_i \to X\}_{i \in I}$ such that $F$ is the image of the morphism of sheaves $\coprod_{i \in I} X_i \to X$. A morphism $F \to G$ in $\ttSh(\ttA^{T})$ is called is a $\tau$-open immersion if for every $X \in \ttA$ and every $X \to G$ the induced morphism $F \times_{G} X \to X$ is a monomorphism whose image is a $\tau$-open in $X$.
The category of schemes $\ttSch(\ttA, T, \tau)$ is defined to be the full subcategory of $\ttSh(\ttA^{T})$ of sheaves $F$ such that there exists a family of $X_{i} \in \ttA$ for $i \in I$ and a morphism $p: \coprod_{i \in I} X_i \to F$ such that $p$ is an epimorphism of sheaves and for each $i$ the morphism $X_i \to F$ is a $\tau$-open immersion.
\end{defn}
\begin{example}Two interesting general categories of schemes that we have in mind are $\ttSch(\ttA, T,\tau)$ where $\ttA=\ttAff(\tC)$. First, the case where $\tau$ of Zariski open immersions and $T$ is the Zariski pre-topology. Second, the category $\tau$ of formal Zariski open immersions and $T$ is the formal Zariski pre-topology. 
\end{example}

\begin{example}
Let $k$ be any field and consider the category $\ttC= \ttVect_{k}$. Then $\ttComm(\ttC)$ is the category of $k$-algebras. Recall that a morphism $A \to B$ of $k$-algebras is of finite presentation in the usual sense when $B$ is finitely generated as an $A$ algebra and the ideal of relations is also finitely generated. Let us temporarily call a morphism TVfp if it satisfies the condition from Definition \ref{defn:FinitePres}. Consider the functor
\[\straightF : \ttSet \to A/\ttComm(\ttC) 
\]
which sends each set to the $A$-algebra freely generated by it. Then Lawvere's work on finitary algebraic theories and Corollary 3.13 and the remark following it in \cite{AR} show that $A \to B$ is TVfp if and only if there exists finite sets $S_g$ and $S_r$ and an isomorphism in $A/\ttComm(\ttC)$ of the form
\[\colim[\straightF S_r \rightrightarrows \straightF S_g ] \to B.
\]
So a morphism $f:A \to B$ of $k$-algebras is of finite presentation in the categorical sense if and only it is of finite presentation in terms of generators and relations. This fact also appears in the algebraic geometry literature. The implication that finite presentation in terms of generators and relations implies finite presentation in the categorical sense was shown in this case in Lemma III.8.8.2.3 of \cite{EGA4}. For the opposite implication see \cite{stacks-project}.  
Now \cite{EGA4} IV.17.9.1 tells us that a morphism of schemes is a flat monomorphism, locally of finite presentation if and only if it is an open immersion. Since a morphism of affine schemes $\spec(B) \to \spec(A)$ is locally of finite presentation if and only if the corresponding morphism $A \to B$ realizes $B$ as an $A$-algebra of finite presentation we can conclude that the Zariski open immersions are precisely the standard (Zariski) open immersions in algebraic geometry. The Zariski topology in the sense of relative algebraic geometry agrees with the Zariski topology in the standard sense in the case $\ttC= \ttVect_{k}$. We should remark that this is the only example in this article for which we can use Lawvere's theory and \cite{AR} in a straightforward way. Another way to characterize the Zariski open immersion is by replacing the flat epimorphism condition with a homotopy epimorphism condition, i.e. that the natural morphism in the derived category $B\otimes^{\mathbb{L}}_{A}B \to B$ is an isomorphism. It is this condition that we examine this article for the category of Banach spaces, and not the flatness condition which would give a different answer. We comment on the case of vector spaces again in Remark \ref{rem:recover}.
\end{example}
\begin{defn}\label{defn:PresheafZarOpIm} A morphism $f:F \to G$ in $\ttPr(\ttAff(\tC))$ is a Zariski open immersion if for every affine scheme $X$ and every morphism $X \to G$, the induced morphism $F \times_{G}X \to X$ is a monomorphism of presheaves and its image agrees with the image of the map of sheaves $\coprod_{i \in I} X_{i } \to X$ corresponding to a family of Zariski opens $X_{i} \to X$.  
\end{defn}
\begin{defn}\label{defn:AffToPreFlat}
Suppose  $ F \in \ttPr(\ttAff(\tC))$ and $X \in \ttAff(\ttC)$. A morphism $X \to F$ is flat if for every $Y \in \ttAff(\ttC)$ and every morphism $Y \to F$ there is a Zariski open cover $\coprod Z_{i} \to X \times_{F} Y$ such that the combined morphisms $Z_{i} \to  X \times_{F} Y \to Y$ are flat. 
\end{defn}
\begin{defn}\label{defn:PresheafFlatMor} A morphism $f:F \to G$ in $\ttPr(\ttAff(\tC))$ is flat if for every affine scheme $X$ and every morphism $X \to G$ and every flat morphism $W \to X \times_{G} F$ the composition $W \to X \times_{G} F \to X$ is flat.
\end{defn}

Consider the Grothendieck site $\ttAff(\tC)^{fpqc}$. The category of simplicial objects in $\ttPr(\ttAff(\tC))$ is denoted $\ttSPr(\ttAff(\tC)).$ This category comes with a (local) model structure as explained in \cite{T}.
\begin{defn}\cite{T} An object $F \in \ttSPr(\ttAff(\tC)^{fpqc})$ is called a pre-Stack. An object $F \in \ttSPr(\ttAff(\tC)^{fpqc})$ is called an fpqc stack if for any $X \in \ttAff(\tC)$ and any hypercovering $H_{\bullet} \to X,$ the natural morphism
\[F(X) \to \text{holim}_{[n] \in \Delta} F(H_n)
\]
is an equivalence of simplicial sets. The category $\ttSt(\ttAff(\ttC)^{fpqc}) =Ho(\ttSPr(\ttAff(\tC)^{fpqc}))$ is the category of stacks.
\end{defn}

An equivalent definition to the above is to define pre-stacks as functors $\ttAff(\ttC)^{op}\to \infty-\ttGpd$, where $\infty-\ttGpd$ is the category of infinity groupoids. This is a full subcategory of the category of $(\infty,1)$-categories for which one could use quasi-categories. $\infty$-groupoids in this model are Kan simplicial sets. Stacks would be pre-stacks which satisfy descent with respect to hypercovers \cite{DHI, L1,L2, TVe2,TVe3}. Similarly we could define (pre-)stacks valued in other categories, for instance a pre-stack in categories would be a functor $\ttAff(\ttC)^{op} \to (\infty,1)-\ttCat$. Note that the category of $(1,1)$-categories embeds into the category of $(\infty,1)$-categories. Using quasi-categories to model $(\infty,1)$-categories, the nerve of a $1$-category is a quasi-category. We can also view dg-categories as stable quasi-categories tensored over complexes \cite{Coh}. We will use this later to view categories of quasi-coherent $\mathcal{O}$-modules, and $\mathcal{D}$-modules (derived or underived) as pre-stacks valued in categories. There is an inductive definition of an $n$-algebraic stack for $n=0,1,2,\dots$. An algebraic fpqc stack on $\ttAff(\tC)$ is an fpqc stack on $\ttAff(\tC)$ which is $n$-algebraic for some $n$.

One can study schemes and stacks using with the fpqc or Zariski topology we have defined above. These could be useful in the analytic context as well when $\ttC$ is the category of Banach spaces and one can study faithfully flat descent in this context. However, the usual G-topology that is usually studied in non-Archimedean geometry as well as the classical metric topology of complex analytic geometry are finer topologies and have smaller ``open sets". The localizations from Definition \ref{defn:AffLocalization} are not flat with respect to the monoidal structure on $\ttBan_{k}$ and for this reason we need to introduce new abstractly defined topologies which will fit in well with these facts. To do so, we need to use quasi-abelian categories.

\section{Quasi-abelian categories}We review some of Schneiders' theory of quasi-abelian categories. These are special cases of Palamodov's semi-abelian categories and of pseudo-abelian categories. They also have the  structure of a (Quillen) exact category in one natural way.
The main reference for this section is \cite{SchneidersQA}.
\begin{defn}
Let $\mathcal{E}$ be an additive category with kernels and cokernels. A morphism $f:E\to F$ is $\mathcal{E}$ is called strict if the induced morphism 
\[\coim(f)\to \im(f)
\] is an isomorphism. 
\end{defn}
Here the image of $f$ is the kernel of the canonical map
$F \to \coker(f)$, and the coimage of $f$ is the cokernel of the canonical map $\ker (f)\to E$. 

\begin{defn}\label{defn:CartCoCart}
Let $\mathcal{E}$ be an additive category with kernels and cokernels. We say that $\mathcal{E}$ is quasi-abelian if it satisfies the following two conditions:
\begin{itemize}
\item In a cartesian square 
\begin{equation*}
\xymatrix{  E' \ar[r]^{f'} \ar[d] & F' \ar[d] \\
E \ar[r]_{f}& F}
\end{equation*}

If $f$ is a strict epimorphism then $f'$ is a strict epimorphism.

\item In a co-cartesian square
\begin{equation*}
\xymatrix{  E \ar[r]^{f} \ar[d] & F \ar[d] \\
E' \ar[r]_{f'}& F'}
\end{equation*}
If $f$ is a strict monomorphism then $f'$ is a strict monomorphism.

\end{itemize}
\end{defn}
\begin{rem}\label{rem:LeftRightInverse} Any morphism in a quasi-abelian category with a right inverse is a strict epimorphism.  Any morphism in a quasi-abelian category with a left inverse is a strict monomorphism. \end{rem}
\begin{defn}
Let $\ttE$ be a quasi-abelian category. 
Let $\xymatrix{ E'\ar[r]^{e'} & E\ar[r]^{e''} & E''}$ be a sequence of maps such that $e''\circ e'=0$. We call such a sequence strictly exact (resp. strictly coexact) if $e'$ (resp. $e''$) is strict and the canonical map 
$\im(e')\to \ker(e'')$ is an isomorphism. A complex 
$E_1\to \cdots \to E_n$ is strictly exact (resp. strictly coexact) if each subsequence $E_{i-1}\to E_i\to E_{i+1}$ is strictly exact (resp. strictly coexact).
\end{defn}

\begin{rem}\label{rem:threetermstrict}
Note that the sequence \begin{equation}\label{equation:threetermstrict}\xymatrix{ 0\ar[r] & E\ar[r]^{u}& F\ar[r]^{v}& G\ar[r]& 0}\end{equation} is strictly exact if and only if $u$ is the kernel of $v$ and $v$ is the cokernel of $u$. Any strict monomorphism or strict epimorphism can be completed to a strictly exact sequence in the form of Equation \ref{equation:threetermstrict}. 
This implies that such a sequence is strictly exact if and only if it is strictly coexact.
\end{rem}Let $\ttE$ be a closed symmetric monoidal quasi-abelian category with all finite limits and colimits. Then an object is flat if and only if tensoring with it preserves strict short exact sequences.

\begin{defn} Call a sequence  $\xymatrix{ E'\ar[r]^{e'} & E\ar[r]^{e''} & E''}$  exact (resp. coexact) if the canonical map 
$\im(e')\to \ker(e'')$ is an isomorphism. A sequence 
$E_1\to \cdots \to E_n$ is exact (resp. coexact) if each subsequence $E_{i-1}\to E_i\to E_{i+1}$ is exact (resp. coexact).
\end{defn}
\begin{rem}\label{rem:threeterm}
Note that the sequence \begin{equation}\label{equation:threeterm}\xymatrix{ 0\ar[r] & E\ar[r]^{u}& F\ar[r]^{v}& G\ar[r]& 0}\end{equation} is  exact if and only if $\ker(u)=0$, $\im(v)=G$ and $\im(u) \to \ker(v)$ is an isomorphism.  Any  monomorphism or epimorphism can be completed to a exact sequence in the form of Equation \ref{equation:threeterm}. 
\end{rem}

The following is remark $1.1.11$ in \cite{SchneidersQA}:
\begin{thm}
Let $\ttE$ be a quasi-abelian category. The class of strictly exact short exact sequences endows $\ttE$ with the structure of an exact category.
\end{thm}

\begin{defn}
Let $\ttE$ be a quasi-abelian category. Let $\straightK(\ttE)$ be its homotopy category. The derived category of $\ttE$ is $\straightD(\ttE)=\straightK(\ttE)/\straightN(\ttE)$ where $\straightN(\ttE)$ is the full subcategory of strictly exact sequences.
\end{defn}

\begin{defn}
Let $\ttE$ be a quasi-abelian category. Let $\straightK(\ttE)$ be its homotopy category. A morphism in 
$\straightK(\ttE)$ is called a strict quasi-isomorphism if its mapping cone is strictly exact. 
\end{defn}

The following is 1.2.17, 1.2.19, 1.2.20, 1.2.27 and 1.2.31 in \cite{SchneidersQA}: 
\begin{thm}
Let $\ttE$ be a quasi-abelian category.
\begin{enumerate}
\item $\straightD(\ttE)$ has a canonical t-structure (the left t-structure). A complex $E$ belongs to $D^{\leq 0}$ if and only if it is strictly exact in strictly positive degrees. $E$ belongs to $D^{\geq 0}$ if and only if it is strictly exact in strictly negative degrees.
\item The heart of this t-structure $\LH(\ttE)$, is equivalent to the localization of the full subcategory of $\straightK(\ttE)$ consisting of complexes E of the form 
\begin{equation}
\xymatrix{ 0\ar[r] & E\ar[r]^{u}& F\ar[r]& 0}
\end{equation}
where $u$ is a monomorphism and $F$ is in degree $0$, by the multiplicative system formed by morphisms which are both cartesian and cocartesian. 
\item There is a canonical fully faithful functor $\sI:\ttE\to LE(\ttE)$. A sequence $E'\to E\to E''$ is strictly exact in $\ttE$ if and only if $\sI(E')\to \sI(E)\to \sI(E'')$ is exact in $\LH(\ttE)$.
\item The functor $I$ induces an equivalence between $\straightD(\ttE)$ and $\straightD(\LH(\ttE))$. This equivalence sends the (left) t-structure on $\straightD(\ttE)$ to the standard t-structure on $\straightD(\LH(\ttE))$.
\end{enumerate}
\end{thm}

\begin{rem}
The embedding $\sI:\ttE\to \LH(\ttE)$ is universal 
in the sense that induces an equivalence for any abelian category $\mathcal{F}$ between left strictly exact functors from $\ttE$ to $\mathcal{F}$ and left exact functors 
from $\LH(\ttE)$ to $\mathcal{F}$. In this sense $\LH(\ttE)$ is the (left) abelian envelope of $\ttE$. See 1.2.33 in \cite{SchneidersQA}.
\end{rem}

\begin{defn}\label{defn:Projective}
Let $\ttE$ be an additive category with kernels and cokernels. An object $I$ is called injective (resp. strongly injective) if the functor $E\mapsto \Hom(E,I)$ is exact (resp. strongly exact), i.e. for any strict (resp. arbitrary) monomorphism $u:E\to F$, the induced map $\Hom(F,I)\to \Hom(E,I)$ is surjective. Dually, $P$ is called projective (resp. strongly projective) if the functor $E\mapsto \Hom(P,E)$ is exact (resp. strongly exact), i.e. for any strict (resp. arbitrary) epimorphism $u:E\to F$, the associated map $\Hom(P,E)\to \Hom(P,F)$ is surjective.     
\end{defn}

\begin{defn}\label{defn:enough}
A quasi-abelian category $\ttE$ has enough projectives if for any object $E$ there is a strict epimorphism $P\to E$ where $P$ is projective. A quasi-abelian category $\ttE$ has enough injectives if for any object $E$ there is a strict monomorphism $E \to I$ where $I$ is injective.
\end{defn}

The following is 1.3.24 in \cite{SchneidersQA}:
\begin{lem}
Let $\ttE$ be a quasi-abelian category. 
\begin{enumerate}
\item An object $P$ of $\ttE$ is projective if and only if $\sI(P)$ is projective in $\LH(\ttE)$.
\item $\ttE$ has enough projectives if and only if $\LH(\ttE)$ has enough projectives. In this case an object of $\LH(\ttE)$ is projective if it is isomorphic  to $\sI(P)$ where $P$ is projective in $\ttE$.
\end{enumerate}
\end{lem}

The following is 1.3.22 in \cite{SchneidersQA}:
\begin{thm}
Let $\ttE$ be a quasi-abelian category with enough projectives (resp. injectives). Let $\tP$ be the full additive subcategory of projective objects (resp. $\tI$ the category of injective objects). The canonical functor $\sK^-(\tP)\to \sD^-(\ttE)$ (resp. $\sK^+(\tI)\to \sD^+(\ttE)$) is an equivalence.
\end{thm}

\subsection{Closed symmetric monoidal quasi-abelian categories.}

The following is 1.5.1 in \cite{SchneidersQA}:
\begin{prop}\label{prop:StrictIFF}
Suppose that $\ttC$ is a closed symmetric monoidal quasi-abelian category with all finite limits and colimits. Suppose that $A \in \ttComm(\ttC)$.  The category $\ttMod(A)$ 
is quasi-abelian and the forgetful functor $\ttMod(A)\to \ttC$ preserves limits and colimits. A morphism in $\ttMod(A)$ is strict if and only if it is strict in $\ttC$. 
\end{prop}

\begin{lem}Suppose that $\ttC$ is a closed symmetric monoidal quasi-abelian category with all finite limits and colimits. Using Remark \ref{rem:LeftRightInverse} we see that for any $V  \in \Mod(A)$ the canonical morphism
\[V \ootimes A \to V
\]
is a strict epimorphism and the canonical morphism
\[V \to \uHom(A,V)
\]
is a strict monomorphism.
\end{lem}

\begin{defn}\label{defn:finiteAmodAbstract}
Suppose that $\ttC$ is a closed symmetric monoidal quasi-abelian category with all finite limits and colimits. An object $V$ is called finite if there is a strict epimorphism $\coprod_{i=1}^{n} \text{id}_{\ttC} \to V$ in $\ttC$ for some finite non-negative integer $n$. In the case that $\ttC=\ttMod(A)$ for $A$ a commutative monoid in a closed symmetric monoidal quasi-abelian category, we denote the full subcategory of finite objects by $\ttMod^{f}(A)$.
\end{defn}

\begin{lem}\label{lem:FreeProjCofreeInj} Suppose that $0 \to L \to M \to N \to 0$ is a strictly exact sequence in $\ttMod(A)$ and $F =A\ootimes P \in \ttMod(A)$ is free and $C=\uHom(A,I) \in \ttMod(A)$ is cofree. Then the sequences 
\[0 \to \Hom_{A}(F,L) \to  \Hom_{A}(F,M)\to  \Hom_{A}(F,N) \to 0
\]
and
\[0 \to \Hom_{A}(N,C) \to  \Hom_{A}(M,C)\to  \Hom_{A}(L,C) \to 0
\]
are exact. If the sequence $0 \to L \to M \to N \to 0$ is only exact and $F$ is strictly free and $C$ is strictly cofree we can make the same conclusion.
\end{lem}
{\bf Proof.}
Using Lemma \ref{lem:MeyerProperties} these sequences are isomorphic to the sequences
\[0 \to \Hom(P,L) \to  \Hom(P,M)\to  \Hom(P,N) \to 0
\]
and 
\[0 \to \Hom(N,I) \to  \Hom(M,I)\to  \Hom(L,I) \to 0
\]
which are exact by definition of projectivity and injectivity (or the strict versions).
\hfill $\Box$
\begin{lem}\label{lem:MaintainProjInj}If $P$ is projective in $\ttC$ then $P\ootimes A$ is projective in $\ttMod(A)$. Similarly, if $I$ is injective in $\ttC$ then $\uHom(A,I)$ is injective in $\ttMod(A)$.
\end{lem}
{\bf Proof.}
Both of these facts are immediately implied by Lemma \ref{lem:FreeProjCofreeInj} together with Remark \ref{rem:threetermstrict}.
\hfill $\Box$
\begin{lem} \label{lem:MonoEpiPreserved} The functor $E \mapsto E\ootimes A$ takes epimorphisms in $\ttC$ to epimorphisms in $\ttMod(A)$. The functor $E \mapsto \uHom(A,E)$ takes monomorphisms in $\ttC$ to monomorphisms in $\ttMod(A)$.
\end{lem}
{\bf Proof.}
This follows directly from the definitions and Lemma \ref{lem:MeyerProperties}.
\hfill $\Box$

\begin{defn}\label{defn:KerFlat} For $A \in \ttComm(\ttC)$ a module $M$ in $\ttMod(A)$ is called kernel flat if for any morphism $f:E \to F$ in $\ttMod(A)$  the natural morphism 
\begin{equation}\label{equation:want4flatAbstr}B\ootimes \ker(f) \to \ker(f_{M})
\end{equation}
is an isomorphism where $f_{M}$ is defined as $\text{id}_{M} \ootimes_{A}f:M \ootimes_{A}E \to M \ootimes_{A}F.$ A morphism $A \to B$ in $\ttComm(\ttC)$ is called kernel flat if it makes $B$ kernel flat over $A$.
\end{defn}

\begin{lem}\label{lem:FlatkerFlat}Suppose that $\ttC$ is a closed symmetric monoidal quasi-abelian category.  An object $V\in \ttC$ is kernel flat if and only it is flat (see Definition \ref{defn:AbstractFlat}). Therefore, a morphism of algebras $p: A \to B$ is kernel flat (see Definition \ref{defn:KerFlat}) if and only if it is flat (Definition \ref{defn:AbstractFlat}) in $\ttComm(\ttC)$. 

\end{lem}

{\bf Proof.}
First of all if $p$ is flat it is clearly kernel flat since a kernel is a type of limit. In the other direction suppose that $p$ is kernel flat. It means that tensoring with $B$ commutes with kernels. Note that every limit over a finite diagram can be written as a combination of finite products and kernels. Finite products are isomorphic to finite coproducts and the functor given by tensoring with $B$ commutes with coproducts and hence it commutes with finite products. Therefore, tensoring with $B$ commutes with finite limits and hence $p$ is flat. 

\ \hfill $\Box$
\begin{lem}\label{lem:HelpWithEnough} Suppose that $\ttC$ is a closed symmetric monoidal quasi-abelian category with all finite limits and colimits. Suppose that $A \in \ttComm(\ttC)$. If the category $\ttC$ has enough projectives 
then the category $\ttMod (A)$ has enough projectives. If the category $\ttC$ has enough injectives 
then the category $\ttMod(A)$ has enough injectives. 
\end{lem}
{\bf Proof.} Suppose that $\ttC$ has enough projectives. Suppose that $V \in \ttMod(A)$. Choose a strict epimorphism in $\ttC$ of the form $P \to V$ where $P$ is projective in $\ttC$.  Lemma \ref{lem:MaintainProjInj} implies that $P\ootimes A$ is projective in $\ttMod(A).$ Consider the morphism $P\ootimes A \to V.$ We need to show it is a strict epimorphism in $\ttMod(A).$ It factorizes as 
\begin{equation}
P \ootimes A \to V \ootimes A  \to V.
\end{equation}
The second morphism is a strict epimorphism because it admits a right inverse. The arrow $P \ootimes A \to V \ootimes A$  is an epimorphism by Lemma \ref{lem:MonoEpiPreserved} and in fact a strict epimorphism because the monoidal product with $A$ is a left adjoint functor and preserves cokernels. Therefore $\ttMod(A)$ has enough projectives.
Suppose that $\ttC$ has enough injectives. Choose a strict monomorphism in $\ttC$ of the form $V \to I$ where $I$ is injective in $\ttC$. Lemma \ref{lem:MaintainProjInj} implies that $\uHom(A,I)$ is injective in $\ttMod(A).$ Consider the morphism $V \to \uHom(A,I).$ We need to show that it is a strict monomorphism in $\ttMod(A).$ It factorizes as \begin{equation}
V \to \uHom(A,V) \to \uHom(A,I).
\end{equation}
Notice that here, we are considering $\uHom(A,V)$ and $\uHom(A,I)$ as elements of 
$\ttMod(A)$ using the action of $A$ on itself. The first arrow is a strict monomorphism because it admits a left inverse. Using Lemma \ref{lem:MonoEpiPreserved}, $\uHom(A,V) \to \uHom(A,I)$ is a monomorphism in $\ttMod(A)$ and in fact a strict monomorphism because the internal Hom from $A$ is a right adjoint functor and preserves kernels. Therefore $\ttMod(A)$ has enough injectives.
\hfill $\Box$
\begin{lem}\label{lem:BaseChangeInj}Suppose that $\sF: \ttC \to \tD$ is a functor which has a right adjoint $\sR$ and which 
preserves strict monomorphisms (preserves monomorphisms). Then an injective (strongly injective) in $\tD$ is an injective (strongly injective) when considered in $\ttC$ via $\sR$ .
\end{lem}
{\bf Proof.}

Suppose that $I \in \tD$ is injective (strongly injective). Then consider a strict monomorphism (monomorphism) $E\to F$ in $\ttC.$ Then $\sF(E) \to \sF(F)$ is a strict monomorphism (monomorphism) in $\ttC.$ We have a commutative diagram
\[\xymatrix{\ttC(F,\sR (I))  \ar[r]& \ttC (E,\sR (I)) \\
\tD(\sF (F),I) \ar[r] \ar[u]& \tD (\sF (E),I). \ar[u]
}
\]
Because the upwards arrows are isomorphisms and the lower horizontal arrow is surjective, the upper horizontal arrow is surjective as well. Therefore, $I$ is  injective (strongly injective) when considered as an object of $\ttC$. 
\hfill $\Box$

\begin{lem}\label{lem:BaseChangeProj}Suppose that $\sG: \ttC \to \tD$ is a functor which has a left adjoint $\sL$ and which preserves strict epimorphisms (preserves epimorphisms). Then a projective (strongly projective) in $\tD$  is projective (strongly projective) when considered in $\ttC$.
\end{lem}
{\bf Proof.}
Suppose that $P \in \tD$ is projective (strongly projective). Consider a strict epimorphism (epimorphism) $E\to F$ in $\tD(A).$ Then by $\sG(E) \to \sG(F)$ is a strict epimorphism (epimorphism) in $\tD.$ We have a commutative diagram
\[\xymatrix{\ttC(\sL (P),E)  \ar[r]& \ttC (\sL (P),F)\\
\tD(P,\sG (E))\ar[u] \ar[r]& \tD(P,\sG (F)) \ar[u].
}
\]
Because the upwards arrows are isomorphisms and the lower horizontal arrow is surjective, the upper horizontal arrow is surjective as well. Therefore, $P$ is projective (strongly projective) when considered in $\ttC$. 
\ \hfill $\Box$

Let $\ttE$ be a closed symmetric monoidal quasi-abelian category and let $A\in \ttComm(\ttE)$. 
The following is contained in 2.1.18 in \cite{SchneidersQA}: If $P$ is projective in $\ttE$ then $A\ootimes P$ is projective in $\ttMod(A)$.

\subsection{Derived Functors}
Let $\sF: \ttC \to \ttD$ be an additive functor betwen quasi-abelian categories $\ttC$ and $\ttD$. Schneiders gave the following definitions in 1.3.2 of \cite{SchneidersQA}
\begin{defn} A full additive subcategory $\ttP$ of $\ttC$ is called $\sF$-projective if:
\begin{enumerate}
\item for any object $V$ of $\ttC$ there is an object $P$ of $\ttP$ and a strict epimorphism $P\to V$
\item in any strictly exact sequence 
\[0 \to V' \to V \to V'' \to 0
\]
of $\ttC$ where $V$ and $V''$ are objects of $\ttP$, $V'$ is as well
\item for any strictly exact sequence 
\[0 \to V' \to V \to V'' \to 0
\]
of $\ttC$ where $V, V'$ and $V''$ are objects of $\ttP$, the sequence 
\[0 \to \sF(V') \to \sF(V) \to \sF(V'') \to 0
\]
is strictly exact in $\ttD.$
\end{enumerate}
 A full additive subcategory $\ttI$ of $\ttC$ is called $\sF$-injective if:
\begin{enumerate}
\item for any object $V$ of $\ttC$ there is an object $I$ of $\ttI$ and a strict monomorphism $V\to I$
\item in any strictly exact sequence 
\[0 \to V' \to V \to V'' \to 0
\]
of $\ttC$ where $V$ and $V''$ are objects of $\ttI$, $V'$ is as well
\item for any strictly exact sequence 
\[0 \to V' \to V \to V'' \to 0
\]
of $\ttC$ where $V, V'$ and $V''$ are objects of $\ttI$, the sequence 
\[0 \to \sF(V') \to \sF(V) \to \sF(V'') \to 0
\]
is strictly exact in $\ttD.$
\end{enumerate}
\end{defn}
Schneiders also includes the following (Lemma 1.3.3 \cite{SchneidersQA})
\begin{lem}
Let $\ttC$ be a quasi-abelian category and let $\ttP$ be a subset of the objects of $\ttC$. Assume that for any object $V$ of $\ttC$ there is a strict epimorphism $P \to V$ with $P \in \ttP$. Then for each object $V$ of $C^{-}(\ttC)$ there is a quasi-isomorphism $u:P \to V$ with $P$ in $C^{-}(\ttP)$ and such that each $u^{k}:P^{k} \to V^{k}$ is a strict epimorphism.
\end{lem}

From this we get (proposition 1.3.5 \cite{SchneidersQA}):
\begin{prop}\label{lem:derivable}
Let $\sF: \ttC \to \ttD$ be an additive functor between quasi-abelian categories $\ttC$ and $\ttD$.
\begin{enumerate}
\item Assume that $\ttC$ has an $\sF$-projective subcategory. Then $\sF$ has a left derived functor $L\sF:D^-(\ttC)\to D^-(\ttD)$.
\item Assume that $\ttC$ has an $\sF$-injective subcategory. Then $\sF$ has a right derived functor $R\sF:D^+(\ttC)\to D^+(\ttD)$.
\end{enumerate}
\end{prop}
\begin{defn}In the situations of Lemma \ref{lem:derivable}, $\sF$ is called explicitly left derivable or explicitly right derivable.
\end{defn}
Here, derived functors are defined as usual by their universal property.

\begin{rem}\label{rem:sse2sse}Note that if $\sF$ is exact (sends strict short exact sequences to strict short exact sequences) then the full subcategory $\ttC$ itself is an $\sF$-projective (and injective) subcategory. Hence exact functors are always derivable.   
\end{rem}

As in the abelian case, projective and injectives form $\sF$-projective and $\sF$-injective subcategories (remark 1.3.21 \cite{SchneidersQA}):
\begin{prop}\label{prop:ProjFProj}
Let $\sF: \ttC \to \ttD$ be an additive functor between quasi-abelian categories $\ttC$ and $\ttD$. 
\begin{enumerate}
\item Assume that $\ttC$ has enough projectives. Then the full subcategory of projective objects is a $\sF$-projective subcategory and therefore can be used to explicitly left derive the functor $\sF.$
\item Assume that $\ttC$ has enough injectives. Then the full subcategory of injective objects is a $\sF$-injective subcategory  and therefore can be used to explicitly right derive the functor $\sF.$
\end{enumerate}
\end{prop}
We also have the following (remark 1.3.7 \cite{SchneidersQA}):
\begin{lem}\label{lem:acyclic}
Let $\sF: \ttC \to \ttD$ be an additive functor between quasi-abelian categories $\ttC$ and $\ttD$. Assume that $\sF$ has a right derived functor $R\sF:D^+(\ttC)\to D^+(\ttD)$. Call an object $I$ $\sF$-acyclic if $R\sF(I)\cong \sF(I)$. Assume that for any object 
$A$, there is an $\sF$-acyclic object $I$ and a monomorphism 
$A\to I$. Then the $\sF$-acyclic objects form a $\sF$-injective subcategory.  Assume that $\sF$ has a left derived functor $L\sF:D^-(\ttC)\to D^-(\ttD)$. Call an object $P$ $\sF$-acyclic if $L\sF(P)\cong \sF(P)$. Assume that for any object 
$A$, there is an $\sF$-acyclic object $P$ and a epimorphism 
$P \to A$. Then the $\sF$-acyclic objects form a $\sF$-projective subcategory. 
\end{lem}
\begin{defn}Let $\ttC$ be a closed symmetric monoidal quasi-abelian category with monoidal structure $\ootimes$. An object $V$ of $\ttC$ is called $\ootimes$-acyclic if $V$ is $\sF$-acyclic for all of the functors $\sF:\ttC \to \ttC$ given by $U \mapsto U \ootimes W$ for any object $W$ in $\ttC$.
\end{defn}
\subsection{Topologies based on homological algebra}
Using the homological algebra in this section, we now introduce some more classes of morphisms and Grothendieck topologies on a closed symmetric monoidal quasi-abelian category $\tC$ with all finite limits and colimits.
\begin{lem}For any morphism $p:\spec(B)\to \spec(A)$ in $\ttAff(\ttC),$ the induced morphism $p_{*}:\ttMod(B)\to \ttMod(A)$ is derivable to a functor $\sD^{-}(B) \to \sD^{-}(A).$ 
\end{lem}
{\bf Proof.} This functor sends strict exact sequences to strict exact sequences so this follows from Remark \ref{rem:sse2sse}.
\ \hfill $\Box$
\begin{defn}\label{defn:homotopyEpi}A morphism $p:\spec(B)\to \spec(A)$ in $\ttAff(\ttC)$ is called a homotopy monomorphism in $\ttAff(\ttC)$ if the induced morphism $p_{*}:\sD^{-}(B) \to \sD^{-}(A)$ is fully faithful.
\end{defn}
Notice that by considering $i=0$ in Definition \ref{defn:homotopyEpi} we see that a homotopy epimorphism in $\ttComm(\ttC)$ is in particular an epimorphism in $\ttComm(\ttC).$

\begin{lem}\label{lem:ComposHom}The composition of homotopy monomorphisms in $\ttAff(\ttC)$ is a homotopy monomorphism in $\ttAff(\ttC)$.
\end{lem}
{\bf Proof.}
This follows from the fact that the composition of fully faithful functors is fully faithful.
\ \hfill $\Box$
\begin{lem}

Assume that $p:\spec(B)\to \spec(A)$ in $\ttAff(\ttC)$ and that the functor $p^{*}:\ttMod(A) \to \ttMod(B)$ given by tensoring with $B$ over $A$ is explicitly left derivable to a functor $\mathbb{L}p^{*}:\sD^{-}(A)\to \sD^{-}(B)$. Then $p$ is homotopy monomorphism if and only if the natural morphism of functors $\mathbb{L}p^{*}p_{*}\to  id_{\sD^{-}(B)}$ is an isomorphism.

\end{lem}
{\bf Proof.}
We have natural isomorphisms for any objects $M,N \in \sD^{-}(B)$ 
\[\Hom_{\sD^{-}(B)}(\mathbb{L}p^{*}p_{*}M,N)   \cong \Hom_{\sD^{-}(A)}(p_{*}M,p_{*}N) .
\]
Therefore, if $\mathbb{L}p^{*}p_{*}\to  id_{\sD^{-}(B)}$ is a isomorphism then $p$ is a homotopy epimorphism. The converse follows from a simple application of the Yoneda lemma.
\ \hfill $\Box$
\begin{lem}\label{lem:HomotopyMon}
 Assume that $p:\spec(B)\to \spec(A)$ is a morphism in $\ttAff(\ttC)$ and that the functor $\ttMod(A) \to \ttMod(B)$ given by tensoring with $B$ over $A$ is explicilty left derivable to a functor $\sD^{-}(A)\to \sD^{-}(B)$. Then $p$ is homotopy monomorphism if and only if $B\ootimes^{\mathbb{L}}_{A}B\cong B$.
\end{lem}
{\bf Proof.} 
For any object $M$ of $\sD^{-}(B)$ we have 
\begin{equation*}
M\ootimes^{\mathbb{L}}_{A}B\cong M\ootimes^{\mathbb{L}}_{B}(B\ootimes^{\mathbb{L}}_{A}B).
\end{equation*}

Hence $\mathbb{L}p^{*}p_{*}\to  id_{\sD^{-}(B)}$ is an isomorphism if and only if  we have natural isomorphisms $M\ootimes^{\mathbb{L}}_{A}B\cong M$ for any $M \in \sD^{-}(B)$ which happens if and only if $B\ootimes^{\mathbb{L}}_{A}B\cong B$.
\ \hfill $\Box$

\begin{defn}\label{defn:fhTVZ} Let $\ttC$ be a closed, symmetric monoidal quasi-abelian category with enough projectives.
 The morphism $\spec(B) \to \spec(A)$ of $\ttAff(\ttC)$ is called a homotopy formal Zariski open immersion if the corresponding morphism $A \to B$ in $\ttComm(\ttC)$ is a homotopy epimorphism.
\end{defn}
\begin{defn}\label{defn:hTVZ} Let $\ttC$ be a closed, symmetric monoidal quasi-abelian category with enough projectives.
 The morphism $\spec(B) \to \spec(A)$ of $\ttAff(\ttC)$ is called a homotopy Zariski open immersion if the corresponding morphism $A \to B$ in $\ttComm(\ttC)$ is a homotopy epimorphism of finite presentation.
\end{defn}

\begin{defn}\label{defn:Amitsur}The Amitsur complex of a morphism $f:A \to B$ in $ \ttComm(\ttC)$ is the complex $\mathscr{A}(f)$ given by 
\[0 \to A \to B \to B\ootimes_{A}B \to  B\ootimes_{A}B\ootimes_{A}B \to \cdots
\]
where the morphism $B^{\ootimes_{A}^{m}}\to B^{\ootimes_{A}^{m+1}}$ is defined by 
\[d(b_{1} \ootimes b_{2} \ootimes \cdots \ootimes b_{m})= \sum_{i=1}^{m+1}(-1)^{i}b_{1} \ootimes \cdots  \ootimes b_{i-1}\otimes 1 \ootimes b_{i} \ootimes \cdots \ootimes b_{m}.
\]
In the case that we can chose a decomposition $f= \prod_{i=1}^{n} f_{i}:A \to \prod_{i=1}^{n} B_{i}=B$ there is a strictly included subcomplex $\mathscr{A}^{a}(f)\to \mathscr{A}(f)$ where $\mathscr{A}^{a}(f)$ is given by the finite complex 
\[0 \to A \to \prod_{1 \leq i_1 \leq n} B_{i}  \to \prod_{1\leq i_1<i_2\leq n} B_{i_1}\ootimes_{A} B_{i_2} 
\to \cdots \to B_{1} \ootimes_{A} \cdots \ootimes_{A} B_{n} \to 0\]
with the induced differentials. 
\end{defn}
\begin{defn}\label{defn:hZt} Let $\ttC$ be a closed, symmetric monoidal quasi-abelian category with enough projectives. We call a full subcategory $\ttA \subset \ttAff(\ttC)$ \emph{homotopy Zariski transversal} if it is closed under fiber products and for any homotopy monomorphism $\spec (B) \to \spec (A)$ in $\ttA$ and for any morphism $\spec (C) \to \spec (A)$ in $\ttA$ 
\[\text{the natural morphism} \ \ B \ootimes_{A}^{\mathbb{L}}C \to  B \ootimes_{A}C \ \ \text{is an isomorphism.}
\]
\end{defn}
\begin{lem}\label{lem:IfThen} If $\ttA \subset \ttAff(\ttC)$ homotopy Zariski transversal subcategory then the base change of a homotopy monomorphism is a homotopy monomorphism.
\end{lem}
{\bf Proof.} If $A \to B$ is a homotopy epimorphism in $\ttA$ and $A \to C$ is any morphism in $\ttA$ then 
\[(B \ootimes_{A} C) \ootimes^{\mathbb{L}}_{C} (B \ootimes_{A} C) \cong (B \ootimes^{\mathbb{L}}_{A} C) \ootimes^{\mathbb{L}}_{C} (B \ootimes^{\mathbb{L}}_{A} C) \cong B \ootimes^{\mathbb{L}}_{A} (B \ootimes^{\mathbb{L}}_{A} C)\cong (B \ootimes^{\mathbb{L}}_{A} B) \ootimes^{\mathbb{L}}_{A} C \cong B \ootimes^{\mathbb{L}}_{A} C\cong B \ootimes_{A} C.
\]
\ \hfill $\Box$

The following definition comes from work \cite{RR} of Ramis-Ruget on quasi-coherent sheaves in complex analytic geometry. Based on their work, quasi-coherent modules in the complex analytic context were discussed in the book \cite{EP}. Our definition is inspired by that one. 
\begin{defn}\label{defn:RRqc}Let $\ttA \subset \ttAff(\ttC)$ be a homotopy Zariski transversal subcategory. For $\spec(A) \in \ttA$ define $\ttMod_{\ttA}^{RR}(A)$ to be the full subcategory of $\ttMod(A)$ consisting of modules $M$ such that $M$ is transversal to all homotopy epimorphisms in $\ttA$. That is, we consider modules $M$ such that the natural morphism \[M \ootimes_{A}^{\mathbb{L}}B \to  M \ootimes_{A}B \] is an isomorphism for all homotopy epimorphisms $A \to B$. We call these RR-quasi-coherent modules.
\end{defn}
\begin{lem}Let $\ttA \subset \ttAff(\ttC)$ be a homotopy Zariski transversal subcategory. If $\spec(C) \to \spec(A)$ is any morphism in $\ttA$ then the morphism $A \to C$ gives $C$ the structure of an object of $\ttMod_{\ttA}^{RR}(A)$.
\end{lem}
{\bf Proof.} This is an obvious consequence of Definition \ref{defn:hZt}.
\ \hfill $\Box$
\begin{lem}\label{lem:RRpreserve} If $\ttA \subset \ttAff(\ttC)$ homotopy Zariski transversal subcategory then pushforwards by morphisms in $\ttA$ preserve the category of RR-quasi-coherent modules. If a morphism is transverse to a RR-quasi-coherent module, then the pull-back of this module by that morphism is also RR-quasi-coherent. Therefore, pullbacks by homotopy monomorphisms or flat morphisms in $\ttA$ also preserve this category. 
\end{lem}
{\bf Proof.} Say we are given a morphism $ \spec(C) \to \spec(A)$ in  $\ttA$. Given any object $M$ of $\ttMod_{\ttA}^{RR}(C)$ and a homotopy monomorphism $\spec(B) \to \spec(A)$, we have 
\[M \ootimes^{\mathbb{L}}_{A} B \cong M \ootimes^{\mathbb{L}}_{C} (C \ootimes_{A}^{\mathbb{L}}B) \cong M \ootimes^{\mathbb{L}}_{C} (C \ootimes_{A}B) \cong M \ootimes_{C} (C \ootimes_{A}B) \cong M \ootimes_{A} B.
\]
This proves the statement about push-forwards. For the statements about pullbacks, fix a morphism $\spec(D) \to \spec(C)$ in $\ttA$ (which we will pullback with) which is transverse to $M$ and a homotopy monomorphism $\spec(E) \to \spec(D)$ in $\ttA$. Then
\[(M \ootimes_{C} D)\ootimes^{\mathbb{L}}_{D}E \cong (M \ootimes^{\mathbb{L}}_{C} D)\ootimes^{\mathbb{L}}_{D}E \cong M \ootimes^{\mathbb{L}}_{C} E \cong M \ootimes_{C} E \cong (M \ootimes_{C} D)\ootimes_{D}E.
\]
\ \hfill $\Box$

As a corollary, we can get a helpful analogue of Lemma \ref{lem:BaseChangeConserv}: 
\begin{cor}\label{cor:BaseChangeConservHzAR}Let $\ttC$ be a closed symmetric monoidal quasi-abelian category. Suppose that $\ttA \subset \ttAff(\ttC)$ homotopy Zariski transversal subcategory. Suppose that a family $\{p_{i}:X_{i}\to X\}$ in $\ttA$ of homotopy monomorphisms is such that the family $\{p^{*}_{i}:\ttMod_{\ttA}^{RR}(X) \to \ttMod_{\ttA}^{RR}(X_i)\}$ has a finite conservative subfamily. Then any pull-back family $\{p_{i}:X_{i}\times_{X} Y\to Y\}$ coming from a base change $Y\to X$ has the same property.
\end{cor}
{\bf Proof.} The proof of Lemma \ref{lem:BaseChangeConserv} uses only the functors $q_{*}$, $q'_{*}$, $p_{i}^{*}$ and ${p'}_{i}^{*}$. The first three preserve the categories of RR-quasi-coherent modules by Lemma \ref{lem:RRpreserve} and the last one  does by Lemmas \ref{lem:IfThen} and \ref{lem:RRpreserve}. Therefore, that same proof works here.
\ \hfill $\Box$

\begin{prop}\label{hfTV1}Let $\ttC$ be a closed, symmetric monoidal quasi-abelian category with enough projectives. Let $\ttA$ be a homotopy transversal subcategory of $\ttAff(\ttC)$.  Consider the families $\{p_{i}:X_{i}\to X\}_{i \in I}$ in $\ttA$ such that the family $\{p^{*}_{i}:\ttMod_{\ttA}^{RR}(X) \to \ttMod_{\ttA}^{RR}(X_i)\}_{i \in I}$ has a finite conservative subfamily and that each $p_{i}$ is a homotopy monomorphism. These families define a pretopology on $\ttA$. 
\end{prop}
{\bf Proof.} In order to show the base change property, consider $q:Y \to X$ in $\ttA$ and let  $q'_{i}$, $p'_i$ and $p_i$ play the role of $q'$, $p'$ and $p$ in diagram (\ref{basechange}). Lemma \ref{lem:BaseChangeConserv} implies that the family $\{p^{*}_{i}\}$ has a finite conservative subfamily. The fact that the $p'_{i}$ are homotopy monomorphisms follows from the assumption.

\ \hfill $\Box$

Let $\ttC$ be a closed symmetric monoidal quasi-abelain category with enough projectives. Let $\ttA \subset \ttAff(\ttC)$ be a homotopy Zariski transversal subcategory. Say that $q:\spec(C) \to \spec(A)$ is arbitrary. Say we are given $A, B, C \in \ttComm(\ttsC)$.  Consider a Cartesian diagram 
\begin{equation}\label{Derivedbasechange}\xymatrix{  \spec(C\otimes^{\mathbb{L}}_{A}B) \ar[r]^{q'} \ar[d]_{p'} & \spec(B) \ar[d]^{p} \\
\spec(C) \ar[r]_{q}& \spec(A). 
}
\end{equation}

Using the notation of diagram (\ref{Derivedbasechange}) there is a natural equivalence 
\begin{equation}\label{eqn:basechangeequation}(\mathbb{L}p^{*})q_{*} \Longrightarrow q'_{*}(\mathbb{L}p'^{*})
\end{equation}
called base change. 

\begin{lem}\label{lem:DerivedBaseChangeConserv}Let $\ttC$ be a closed symmetric monoidal quasi-abelain category with enough projectives. Let $\ttA \subset \ttAff(\ttC)$ be a homotopy Zariski transversal subcategory. Suppose that a family $\{p_{i}:X_{i}\to X\}$ of homotopy monomorphisms in $\ttA$ is such that the family $\{\mathbb{L}p^{*}_{i}:D^{-}(X) \to D^{-}(X_i)\}$ has a finite conservative subfamily. Then any pull-back family $\{p_{i}:X_{i}\times_{X} Y\to Y\}$ coming from a base change $Y\to X$ in $\ttA$ has the same property.
\end{lem}
{\bf Proof.}  In order to show the base change property, consider $q:Y \to X$. Choose a finite set $J \subset I$ such that $\prod_{i \in J} \mathbb{L}p^{*}_{i}$ is conservative. Consider the functor $\prod_{i \in J} \mathbb{L}{p'}^{*}_{i}$ where $q'_{i}$, $p'_i$ and $p_i$ play the role of $q'$, $p'$ and $p$ in diagram (\ref{basechange}). In order to show it is conservative, its enough to show that $\prod_{i \in J} q'_{i*}\mathbb{L}{p'}^{*}_{i}$ is conservative. By the definition of homotopy Zariski transversal, the natural morphism $X_{i}\times_{X} Y \to  X_{i}\times^{h}_{X} Y$ is an isomorphism. Therefore, we can use  (\ref{eqn:basechangeequation}) to conclude that $\prod_{i \in J} q'_{i*}\mathbb{L}{p'}^{*}_{i}$ is natrually equivalent to $(\prod_{i \in J} \mathbb{L}p^{*}_{i}) q_{*}$ which is conservative since $q_{*}$ is conservative.
\ \hfill $\Box$
\begin{rem}Notice that the {\it derived} base change of a homotopy monomorphism is always a homotopy monomorphism but there is no reason that such a thing would work for ordinary base changes. In order to construct a Gronthendieck topology, we need the admissible opens to be preserved by base change. Because the current article is not about derived geometry which will be discussed in \cite{BeKr2}, we defined our Gronthendieck topology on a homotopy transversal subcategory $\ttA \subset \ttAff(\ttC)$. The main application in this article is the case where $\ttC$ is the category $\ttBan_{k}$ for a non-Archimedean valued field $k$ and $\ttA$ is the opposite to the category of affinoid algebras. In the preprints \cite{BaBe} and \cite{BaBeKr} we take $\ttC$ to be the category of Ind-Banach spaces and $\ttA$ to be categories of dagger affinoid and Stein algebras which we define. In these upcoming articles we do not assume the field is non-Archimedean so they apply to complex analytic geometry as well. In the article on derived analytic geometry \cite{BeKr2}, we take $\ttC$ to be the monoidal model category of simplicial Ind-Ban spaces.
\end{rem}

\begin{defn}\label{defn:fhZar}Let $\ttC$ be a closed, symmetric monoidal quasi-abelian category with enough projectives.  Let $\ttA$ be homotopy Zariski transversal subcategory of $\ttAff(\ttC)$ . The topology coming from Proposition \ref{hfTV1} is called the formal homotopy Zariski topology on $\ttA$. When equipped with this topology, we denote this category by $\ttA^{fhZar}.$ The category of sheaves of sets is denoted $\ttSh(\ttA^{fhZar}).$ The category of schemes is denoted by $\ttSch(\ttA^{fhZar}).$
\end{defn}

\begin{prop}\label{hTV1}
Consider the families $\{p_{i}:X_{i}\to X\}_{i \in I}$ in $\ttA$ such that the family $\{p^{*}_{i}:\ttMod^{RR}(X) \to \ttMod^{RR}(X_i)\}_{i \in I}$ has a finite conservative subfamily and that each $p_{i}$ is a homotopy monomorphism of finite presentation. These families define a pretopology on $\ttA$. 
\end{prop}
{\bf Proof.} In order to show the base change property, consider $q:Y \to X$ in $\ttA$ and let  $q'_{i}$, $p'_i$ and $p_i$ play the role of $q'$, $p'$ and $p$ in diagram (\ref{basechange}). Lemma \ref{lem:BaseChangeConserv} implies that the family $\{p^{*}_{i}\}$ has a finite conservative subfamily. The fact that the $p'_{i}$ are homotopy monomorphisms follows from the assumption.

\ \hfill $\Box$

\begin{defn}\label{defn:hZar}Let $\ttC$ be a closed, symmetric monoidal quasi-abelian category with enough projectives. Let $\ttA$ be homotopy Zariski transversal subcategory of $\ttAff(\ttC)$. The topology coming from Proposition \ref{hTV1} is called the homotopy Zariski topology on $\ttA$. When equipped with this topology, we denote this category by $\ttA^{hZar}.$ The category of sheaves of sets is denoted $\ttSh(\ttA^{hZar}).$ The category of schemes is denoted by $\ttSch(\ttA^{hZar}).$
\end{defn}


\section{Main Theorems}
\subsection{From Berkovich geometry to Banach algebraic geometry}
 Let $k$ be a non-Archimedean valuation field. 
We now introduce the full subcategory $\ttAfnd_{k} \subset \ttComm(\ttBan_{k})$. The objects in this category are the $k$-affinoid algebras from the literature on Berkovich analytic spaces.
\begin{defn}\label{defn:NiceAlg} For any finite ordered set of positive real numbers $r=(r_1, \dots, r_n)$, let \[k\{r_{1}^{-1}x_{1}, \dots, r_{n}^{-1} x_{n}\} \in \ttComm(\ttBan_{k})\] be the completion of $k[x_{1}, \dots,  x_{n}]$ with respect to the norm 
\[\|\sum_{I \in \mathbb{Z}^{n}_{\geq 0}} a_{I} x^{I}\|_{r}=\text{max}_{I}\{|a_{I}|r^{I}\}
\] 
where for each $I= (i_1, \dots, i_n) \in \mathbb{Z}^{n}_{\geq 0}$,  $a_{I}= a_{i_1, \dots, i_n}$ and $x^{I}=x^{i_1}_{1} \cdots x^{i_n}_{n}$. 
\end{defn}
\begin{lem} More concretely, we have
\begin{equation}\label{equation:concrete}
k \{r_{1}^{-1}x_{1}, \dots, r_{n}^{-1} x_{n}\} = \{ \sum_{} a_{I}x^{I} \in k[[x_1, \dots, x_n]] \ \ | \ \ \lim_{|I| \to \infty} |a_{I}| r^{I} = 0 \}
\end{equation}
where $|I|= i_{1} + \cdots +i_{n}$ and we equip the right hand side with the norm 
\[\|\sum_{I \in \mathbb{Z}^{n}_{\geq 0}} a_{I} x^{I}\|_{r}=\sup_{I}\{|a_{I}|r^{I}\}
\]
\end{lem}
{\bf Proof.} First, notice that the right hand side is complete: given a Cauchy sequence $\{f^{(j)}\}$ in the right hand side of (\ref{equation:concrete}), note that for each $I$, the sequence of coefficients $f^{(j)}_{I}$ of $x^I$ must be Cauchy, and hence converge in $k$, to some $f_{I}.$ Define $f=\sum_{I \in \mathbb{Z}^{n}_{\geq 0}} f_{I} x^{I}.$ It is easy to check that $f$ lies in the right hand side using the inequality
\[|f_{I}|r^{I} \leq \max \{|f_{I} - f^{(j)}_{I}|r^{I},|f^{(j)}_{I}|r^{I} \}.
\]
 
If $f$ is an element of the right hand side, let $f^{(j)} = \sum_{I \in \mathbb{Z}^{n}_{\geq 0}, |I| \leq j} f_{I} x^{I}.$ These polynomials $f^{(j)}$ converge to $f$ so the polynomials are dense in the right hand side.
\ \hfill $\Box$
\begin{rem}Notice by Lemma \ref{lem:sumsConverge} that the right hand side consists precisely of the formal power series which can be evaluated on elements $(x_1, \dots, x_n) \in k^{n}$ such that $|x_{i}| \leq r_i$ for $i=1, \dots, n$.
\end{rem}

\begin{defn} A $k$-affinoid algebra $\mathcal{A}$ is an object in $\ttComm(\ttBan_{k})$ which admits an admissible surjection 
\[k\{r_{1}^{-1}x_{1}, \dots, r_{n}^{-1} x_{n}\} \to \mathcal{A}.
\]
whose multiplication is contracting (see Definition \ref{defn:nonexpanding}). The full category of $\ttComm(\ttBan_{k})$ consisting of such objects is denoted $\ttAfnd_{k}$.
\end{defn}
 Notice that the completed symmetric algebra in $\ttBan^{\leq 1}_{k}$ on $k_{r_1} \oplus k_{r_2} \oplus \cdots \oplus k_{r_n}$ is just $k\{r_{1}^{-1}x_{1}, \dots, r_{n}^{-1} x_{n}\}$ from Definition \ref{defn:NiceAlg}. So these algebras are precisely the free commutative monoids in $\ttBan^{\leq 1}_{k}$. In fact, we have
\begin{lem}Affinoid algebras are precisely the finitely presented objects in $\ttComm(\ttBan^{\leq 1}_{k})$. This means that they are exactly those objects which are the cokernel of a morphism in $\ttComm(\ttBan^{\leq 1}_{k})$ of free commutative algebra objects: 
\[k\{s_{1}^{-1}y_{1}, \dots, s_{m}^{-1} y_{m}\} \to k\{r_{1}^{-1}x_{1}, \dots, r_{n}^{-1} x_{n}\}
\]
\end{lem}
{\bf Proof.} Let $\mA = k\{r_{1}^{-1}x_{1}, \dots, r_{n}^{-1} x_{n}\}/I$ be an affinoid algebra. It is known \cite{Te} that $k\{r_{1}^{-1}x_{1}, \dots, r_{n}^{-1} x_{n}\}$ is Noetherian. Therefore, $I$ must be finitely generated by some finite set of non-zero elements $f_1, \dots, f_m \in k\{r_{1}^{-1}x_{1}, \dots, r_{n}^{-1} x_{n}\}$. Let $s_{j} = \|f_{j}\|$ for $j= 1, \dots, m$. Consider the morphism of algebras
\[k\{s_{1}^{-1}y_{1}, \dots, s_{m}^{-1} y_{m}\} \to k\{r_{1}^{-1}x_{1}, \dots, r_{n}^{-1} x_{n}\}
\]
determined by the (contracting) morphism of Banach spaces 
\[k_{s_1} \oplus \cdots \oplus k_{s_m} \to k\{r_{1}^{-1}x_{1}, \dots, r_{n}^{-1} x_{n}\}
\]
determined by the maps sending $1 \in k_{s_j}$ to $f_j$. Its image is the ideal generated by $f_1, \dots, f_m$.

\ \hfill $\Box$

\begin{defn}\label{defn:AffLocalization} A $k$-affinoid localization is a morphism $\mathcal{A} \to \mathcal{D}$ of $k$-affinoid algebras such that the morphism $|\mathcal{M}(\mathcal{D})|\to |\mathcal{M}(\mathcal{A})|$ is injective, the image of $|\mathcal{M}(\mathcal{D})|$ is closed in $|\mathcal{M}(\mathcal{A})|$ and any morphism of $k$-affinoid algebras $\mathcal{A} \to \mathcal{B}$ such that $|\mathcal{M}(\mathcal{B})|$ lands in the image of  $|\mathcal{M}(\mathcal{D})|$ factors as $\mathcal{A} \to \mathcal{D} \to \mathcal{B}$. 
\end{defn}
Note that $k$-affinoid localizations correspond to the subspaces known as affinoid domain embeddings. They are usually written as $\mA \to \mA_{V}$ where $V$ is the closed image and they satisfy the property (proven in 2.2.2 (iv) of \cite{Ber1990}) that 
\[\mA_{V_1 \cap V_2} \cong \mA_{V_1} \wotimes_{\mA} \mA_{V_2} 
\]
for affinoid domains $V_1$ and $V_2$.

Three important examples are the localizations corresponding to the rational, Weirstrass and Laurent domains. Any affinoid domain is a union of a finite collection of rational domains.

\begin{defn}\label{defn:RatDom}A rational localization is a morphism of affinoid algebras of the form 
\[\mA \to \mA \{\frac{T_1}{r_1}, \dots, \frac{T_m}{r_m} \}/(g T_1 - f_1, \dots, gT_m - f_m)
\]
for some $g, f_j \in \mA$ such that $(f_1,\dots, f_m,g)=1$.
\end{defn}
Notice that the norm of this morphism is one. 

\begin{rem}

Proposition 2.2.4 (ii) of \cite{Ber1990} tells us that $\mA_{V}$ is a flat $\mA$-algebra with respect to the algebraic (non-completed) tensor product, however it is not flat with respect to $\wotimes_{\mA}$. Note that since completions commute with colimits, a rational localization $\mA_{V}$ is the completion of 
$\mA[T_1, \dots, T_n]/(g T_1 - f_1, \dots, gT_n - f_n)$ with respect to the residue semi-norm coming from the semi-norm in Definition \ref{def:NDiskrelA}. We have an 
isomorphism of $\mA$ modules 
\[\mA_g= \mA[S]/(g S-1) \to \mA[T_1, \dots, T_n]/(g T_1 - f_1, \dots, gT_m - f_m)
\]
given by $S = b+\sum_{i=1}^{n}a_i T_i$ where $bg + \sum_{i=1}^{n}a_if_i = 1$. The inverse is given by $T_i=f_iS$. Therefore $\mA[T_1, \dots, T_n]/(g T_1 - f_1, \dots, gT_n - f_n)$ is flat over $\mA$ in the algebraic sense. Given any $\mM \in \ttMod(\mA)$ we have a morphism $\mM_{g} \to \mM\wotimes_{\mA} \mA_V$ and $\mM\wotimes_{\mA} \mA_V$ is the completion of $\mM_{g}=\mM[T_1, \dots, T_n]/(g T_1 - f_1, \dots, gT_m - f_m)$ with the residue semi-norm of the algebraic tensor product $\mM[T_1, \dots, T_n]$.
\end{rem}
\begin{defn}\label{defn:WDom}A Weirstrass localization is a rational localization for which $g=1$.
\end{defn}
\begin{defn}\label{defn:LDom}A Laurent loaclization is a localization of the form 
\[\mA\to \mA\{\frac{T_1}{p_1}, \dots, \frac{T_n}{p_n}, q_1S_1, \dots, q_mS_m\}(T_1 - f_1, \dots, T_n-f_n, g_1S_1-1, \dots, g_mS_m-1)
\]
where the $p_i$ and $q_j$ are positive reals and $f_i,g_j \in \mA$.
\end{defn}
\begin{defn}\label{defn:CAcomplex}
Fix a system $\mA \to \mA_{V_i}$ which give a cover of $\mM(\mA)$ by a finite collection of affinoid domains $V_i$. The \v{C}ech-Amitsur complex is the complex 
\begin{equation}\label{equation:CA} 0 \to \mM \to \prod_{i} \mM \wotimes_{\mA} \mA_{V_i} \to \prod_{i,j}\mM \wotimes_{\mA} \mA_{V_i} \wotimes_{\mA} \mA_{V_j} \to \cdots.
\end{equation}
To any morphism $\mM \to \mN$ we have the obvious morphism of \v{C}ech-Amitsur complexes. 
The complex written here differs from the standard long exact sequence defined in section 8.2 of \cite{BGR} in that we consider the completed tensor products everywhere whereas they only complete the tensor products between the various $\mA_{V_i}$ terms. The standard complex mentioned is exact as proven in the Acyclicity Theorem (Proposition 2.2.5 of \cite{Ber1990}). The above complex is not always exact, however.
\end{defn}
\begin{lem}\label{lem:CAcplxISexact} The complex in definition \ref{defn:CAcomplex} is strictly exact for modules $\mM$ such that the natural morphism 
\[\mM\wotimes^{\mathbb{L}}_{\mA}\prod_{i}\mA_{V_i}\to \mM\wotimes_{\mA}\prod_{i}\mA_{V_i}
\] is an isomorphism.
\end{lem}

{\bf Proof.}

This follows from an easy double complex argument using free resolutions.
\ \hfill $\Box$

\begin{lem}\label{lem:easyLW} Let $\mA$ be an affinoid algebra over a non-Archimedean  valuation field $k$. Let $\mA_{V}$ be the localization $\mA_{V}= \mA\{\frac{T}{s}\}/(T-f)$ for some $f\in \mA$ or $\mA_{V}= \mA\{\frac{T}{s}\}/(gT-1)$ for some $g\in \mA$. Let $\mB$ be an affinoid $\mA$-algebra. Then the natural morphism
\[\mB \wotimes^{\mathbb{L}}_{\mA} \mA_{V} \to   \mB \wotimes_{\mA} \mA_{V}\] is an isomorphism in $\sD^{-}(\mA)$. In particular,  by taking $\mB=\mA_{V}$ we see that the morphisms $\mA\to \mA\{\frac{T}{s}\}/(gT-1)$ and  the morphism $\mA\to \mA\{\frac{T}{s}\}/(T-f)$ homotopy epimorphisms.
\end{lem}
{\bf Proof.} 
Notice that it enough to show that for any affinoid algebra $\mC$ and any element $f\in \mC$ the morphism 
\begin{equation}\label{mono1}\mC\{\frac{T}{s}\} \stackrel{T-f}\to \mC\{\frac{T}{s}\}
\end{equation}
is a strict monomorphism and for any $g \in \mC$ the morphism
\begin{equation}\label{mono2}\mC\{\frac{T}{s}\} \stackrel{gT-1}\to \mC\{\frac{T}{s}\}
\end{equation}
is a strict monomorphism. Indeed, we can use in the case $\mC=\mA$
\[\mA\{\frac{T}{s}\} \stackrel{T-f}\to \mA\{\frac{T}{s}\}
\]
or
\[\mA\{\frac{T}{s}\} \stackrel{gT-1}\to \mA\{\frac{T}{s}\}
\]
as resolutions of $\mA_{V}$ whose terms are projective and $\wotimes_{\mA}$-acyclic. Then by taking the completed tensor product over $\mA$ with $\mB$ we find a representative for $\mB \wotimes^{\mathbb{L}}_{\mA} \mA_{V}$ which looks like 
\[\mB\{\frac{T}{s}\} \stackrel{T-f}\to \mB\{\frac{T}{s}\}
\]
or
\[\mB\{\frac{T}{s}\} \stackrel{gT-1}\to \mB\{\frac{T}{s}\}
\]
where the $f$ and $g$ here are the images of the original ones in $\mB$. These complexes are immediately recognized as also being special cases of (\ref{mono1}) and (\ref{mono2}) in the case $\mC=\mB.$ This shows that the operation of taking the completed tensor product over $\mA$ with $\mB$ is (strictly) exact.

We start by showing that (\ref{mono1}) and (\ref{mono2}) are monomorphisms which using Lemma \ref{lem:BanProps} simply means to show that they are injective. Multiplication by $gT-1$ is clearly injective as can be seen by looking at the lowest order term in $T$. Consider the ascending sequence of ideals of $\mC$ given by the kernel of the morphisms of $\mC$ to itself defined by multiplication by $f^{i}$: \[\ker(f) \subset \ker(f^{2}) \subset \ker(f^{3}) \subset \cdots. \] Since $\mC$ is Noetherian by Proposition 2.1.3 of \cite{Ber1990} this sequence must terminate at some $\ker(f^{N}).$ Suppose that \[(T-f)\left(\sum_{j=0}^{\infty} a_{j} T^{j} \right)= 0\] with $\sum_{j=0}^{\infty} a_{j} T^{j} \neq 0$ and let $a_{i}$ be the first non-zero coeficient. Then $f^{N}a_{i+N}= a_{i}$ and so $f^{N+1}a_{i+N}= fa_{i}=0$ and so $a_{i+N} \in \ker(f^{N+1})= \ker(f^{N}).$ This shows that $a_{i}=0,$ a contradiction. Therefore, multiplication by $T-f$ is injective.
Let us assume for a moment that $k$ is non-trivially valued. To show that these morphisms are strict one must simply show that the set-theoretic image is closed. However, this set-theoretic image in both cases is an ideal and all ideals in affinoid algebras are closed by Proposition 2.1.3 of \cite{Ber1990}. Therefore, they are strict monomorphisms in this case. Now if $k$ is arbitrary we find that these two morphisms become strict after tensoring with the non-trivially valued field introduced in Proposition 2.1.2 of \cite{Ber1990}. Therefore, by this same proposition, they were strict monomorphisms all along.
\ \hfill $\Box$

\begin{lem}\label{lem:LCSrHoEpis}Let $\mA_{V}$ be a rational localization of an affinoid $k$-algebra $\mA$. Let $\mB$ be an affinoid $k$-algebra. Then the natural morphism
\[\mA_{V}\wotimes^{\mathbb{L}}_{\mA}\mB \to \mA_{V}\wotimes_{\mA}\mB \] 
is an isomorphism in $\sD^{-}(\mA)$. In particular taking $\mB=\mA_{V}$, any rational, Weirstrass, or Laurent localization $\mA\to \mA_{V}$ is a homotopy epimorphism.
\end{lem}

{\bf Proof.}

Assume now that $k$ is non-trivially valued and consider the rational localization 
\[\mA_{V} = \mA \{\frac{T_1}{r_1}, \dots, \frac{T_m}{r_m} \}/(g T_1 - f_1, \dots, gT_m - f_m)
\]
where $r_{i}>0$ and $(f_1, \dots, f_m,g)=1.$
Notice that following Proposition 1 of 7.2.4 of \cite{BGR} $|g|$ cannot be arbitrarily small for $|\;| \in \mM(\mA_{V})$ and in fact we can realize the rational localization as an Laurent localization of a Weierstrass localization. That is, there is an $\epsilon>0$ such that
\[\mA_{V}
\cong \left(\mathcal{A}\{ \frac{S}{\epsilon^{-1}}\}/(gS-1)\right)\{\frac{T_1}{r_1}, \dots, \frac{T_m}{r_m} \}/( T_1 - \frac{f'_1}{g'}, \dots, T_m -\frac{f'_m}{g'})
\]
where the $f'_i$ and $g'$ are the image of $f_i$ and $g$ in $\mathcal{A}\{\epsilon S\}/(gS-1).$ This corresponds to the fact that we can choose $\epsilon >0$ such that the conditions $|f_{i}|\leq r_{i}|g|$ and $(f_1, \dots, f_m,g)=1$ are equivalent to $|g|>\epsilon$ and $|\frac{f_{i}}{g}|\leq r_{i}$. Therefore Lemma \ref{lem:easyLW} implies that the natural morphism
\[\mA_{V}\wotimes^{\mathbb{L}}_{\mA}\mB \to \mA_{V}\wotimes_{\mA}\mB \] 
is an isomorphism in $\sD^{-}(\mA)$. 
Now if $k$ is trivially valued, consider, using Proposition 2.1.2 of \cite{Ber1990}, a field $K_r$  containing $k$ which is flat over $k$ with respect to the completed tensor product and non-trivially valued. When we take the tensor product of the canonical morphism \[\mA_{V}\wotimes^{\mathbb{L}}_{\mA}\mB\to \mA_{V}\wotimes_{\mA}\mB\] with $K_{r}$ over $k$ we get, using the flatness of $K_r$, that
\[(\mA_{V} \wotimes_{k} K_r)\wotimes^{\mathbb{L}}_{\mA  \wotimes_{k} K_r}(\mB \wotimes_{k} K_r) \to \mA_{V} \wotimes_{\mA}\mB \wotimes_{k} K_r
\] is an isomorphism in $\sD^{-}(\mA\wotimes_{k}K_r)$ and hence the original morphism $\mA_{V}\wotimes^{\mathbb{L}}_{\mA}\mA_{V}\to \mA_{V} \wotimes_{\mA}\mB$ is an isomorphism in $\sD^{-}(\mA).$

\ \hfill $\Box$
\begin{lem}Let $\mA_{W_1}$ and $\mA_{W_2}$ be affinoid localizations of an affinoid $k$-algebra $\mA$ corresponding to subdomains $W_1$ and $W_2$. Assume also that $W_1 \cup W_2$ is an affinoid subdomain. Let $\mB$ be an affinoid $k$-algebra. Assume that the morphisms 
\[\mA_{W_i}\wotimes^{\mathbb{L}}_{\mA}\mB \to \mA_{W_i}\wotimes_{\mA}\mB
\]
are isomorphisms for $i=1,2$. Assume also that the morphism
\[\mA_{W_1 \cap W_2}\wotimes^{\mathbb{L}}_{\mA}\mB \to \mA_{W_1 \cap W_2}\wotimes_{\mA}\mB
\]
is an isomorphism. Then the morphism 
\[\mA_{W_1 \cup W_2}\wotimes^{\mathbb{L}}_{\mA}\mB \to \mA_{W_1 \cup W_2}\wotimes_{\mA}\mB
\]
is an isomorphism.
\end{lem}
{\bf Proof.}
This follows immediately from considering the strict short exact sequence 
\[0 \to \mA_{W_1 \cup W_2} \to \mA_{W_1} \times \mA_{W_2} \to \mA_{W_1 \cap W_2} \to 0
.\]

\ \hfill $\Box$
Finally, we are able to show in the following theorem that affinoid subdomains of affinoids give examples homotopy monomorphisms of affine schemes in the abstract sense.
\begin{thm} \label{thm:TensLocs}Let $\mA_{V}$ be an affinoid localization of an affinoid $k$-algebra $\mA$. Let $\mB$ be an affinoid $k$-algebra. Then the natural morphism  
\[\mA_{V}\wotimes^{\mathbb{L}}_{\mA}\mB \to \mA_{V}\wotimes_{\mA}\mB \]
is an isomorphism in $\sD^{-}(\mA)$.
In particular let $\mA_{W_1}$ and $\mA_{W_2}$ be affinoid localizations of an affinoid algebra $\mA$. Then $\mA_{W_1} \wotimes^{\mathbb{L}}_{\mA} \mA_{W_2} \cong \mA_{W_1} \wotimes_{\mA} \mA_{W_2}.$ Therefore taking $W_1=W_2$, any affinoid localization is a homotopy epimorphism.
\end{thm}

{\bf Proof.}
For Weirstrass or Laurent localizations this follows immediately from using induction on Lemma \ref{lem:easyLW} or the fact that it holds for the more general set of rational localizations which we now consider. When $\mA_{W_1}$ and $\mA_{W_2}$ are rational localizations this has been shown already in Lemma \ref{lem:LCSrHoEpis}. Suppose now that $V=W_1$ is an rational domain and $W=W_2=V_1 \cup \cdots \cup V_{N}$ is an affinoid domain written as a union of rational domains. For $N=1$ the claim is true. The induction step follows from considering the derived tensor product of $\mA_{V}$ over $\mA$ with the short exact sequence 
\begin{equation}\label{equation:ind}0\to \mA_{W} \to \mA_{V_1 \cup \cdots \cup V_{N-1}} \times \mA_{V_N} \to \mA_{V_1 \cup \cdots \cup V_{N-1}} \wotimes_{\mA} \mA_{V_N} \to 0.
\end{equation}
Assume by induction that the for some $N>2$ the derived tensor product of the rational localization $\mA_V$ and the affinoid algebra corresponding to the union of $N-1$ or fewer rational domains is equivalent to the ordinary tensor product. This implies that the derived tensor product of $\mA_{V}$ with $\mA_{V_1 \cup \cdots \cup V_{N-1}} \times \mA_{V_N}$ and with $\mA_{V_1 \cup \cdots \cup V_{N-1}} \wotimes_{\mA} \mA_{V_N} =\mA_{(V_1 \cap V_{N})\cup \dots \cup (V_{N-1} \cap V_{N})}$ are equivalent to the ordinary tensor products. Hence the same holds with $\mA_{W}$. Finally, one takes $W$ as $W_1$ in (\ref{equation:ind}) and considers the derived tensor product with $\mA_{W_2}$ and uses a similar induction to do the general case.

\ \hfill $\Box$

\begin{rem}\label{rem:recover}This is analogous to Lemma 2.1.4 (1) of \cite{TVe5} where it is shown that homotopy Zariski open immersions in the category of affine schemes relative to the closed symmetric monoidal category of abelian groups give precisely the ordinary notion of a Zariski open imemrsions.  In that reference first a special case (inverting a single element of the ring) was shown and the general case follows by their descent formalism.
\end{rem}
\begin{defn}
For any topological space $X$, a {\it quasi-net} is a set $T$ of subsets of $X$ such that any point $x \in T$ has a neighborhood of the form $\cup_{i=1}^{n}V_{i}$ with $x \in V_{i} \in T$ for $1 \leq i \leq n$. A {\it net} is a quasi-net $T$ such that such that for every $U,V \in T$ the set $\{W \in T | W \subset U \cap V\}$ is a quasi-net of subsets of $U \cap V$.  A $k$-analytic space is a locally Hausdorff topological space $X$, a net $\tau_{0}$ on $X$, a functor $\phi: \tau_{0} \to \ttAfnd_{k},$ and an invertible natural transformation $\ttTop \Longrightarrow \ttTop \circ \phi$. 
\end{defn}
\begin{defn} 
For any $k$-affinoid algebra $\mathcal{A}$, the topological space $|\mathcal{M}(\mathcal{A})|$ is defined to be the set of non-archimedean bounded semivaluations $| \ |$ on $\mA$ equipped with the weakest topology such that for each $f \in \mA$, the maps $|\mM(\mA)| \to \mathbb{R}_{+}$ defined by sending $| \ |$ to $|f|$ is continuous.
\end{defn}

\begin{defn}A {\it k-affinoid space} is a locally ringed space of the form $\mathcal{M}(\mathcal{A})=(|\mathcal{M}(\mathcal{A})|,\mathcal{O}_{\mathcal{M}(\mathcal{A})})$ where $\mathcal{A}$ is a $k$-affinoid algebra and $\mathcal{O}_{\mathcal{M}(\mathcal{A})}(U)$ is the limit over $\mathcal{A}_{V}$, where $V \subset U$ is a finite union of affinoid domains. The category of $k$-affinoid spaces defined to be a full subcategory of the category of locally ringed spaces of the given form. 
\end{defn}
The category of $k$-affinoid spaces is equivalent to the category  $\ttAfnd_{k}^{op}$. So we treat $\mathcal{M}$ as a functor giving this equivalence from $\ttAfnd_{k}^{op}$ to the category of  $k$-affinoid spaces.
\begin{defn}
A {\it k-analytic space} consists of a triple $(X,\tau,\mathcal{A})$ where $X$ is a locally Hausdorff topological space, $\tau$ is a net on $X$, and for each $V \in \tau$, $\mathcal{A}(V)$ is a k-affinoid algebra along with a homeomorphism $|\mathcal{M}(\mathcal{A}(V))| \cong V$ (functorially assigned to the elements of $\tau$) such that if $V, V' \in \tau$ and $V' \subset V$ then $V'$ is an affinoid subdomain of $\mathcal{M}(\mathcal{A}(V))$ with coordinate ring $\mathcal{A}(V')= \mathcal{A}(V)_{V'}$.  In the event that for every $U \in \tau_{2},$ $\tau_{1}$ restricted to $g^{-1}(U)$ is an atlas of $g^{-1}(U),$ a morphism  $(X_{1},\tau_{1},\mathcal{A}_{1}) \to (X_{2},\tau_{2},\mathcal{A}_{2})$ consists of a continuous and $G$-continuous map $g:X_{1} \to X_{2}$ along with bounded homomorphisms $g^{\#}_{U,V}:\mathcal{A}_{2}(U) \to \mathcal{A}_{1}(V)$ for every $U \in \tau_{2}, V  \in \tau_{1}$ with $g(V) \subset U$ such that for every $V, V' \in \tau_{1}$ with $V'\subset V$ and $U, U' \in \tau_{2}$ with $U' \subset U$ such that $g(V) \subset U$ and $g(V') \subset U'$ the diagram 
\begin{equation}
\xymatrix{ \mathcal{A}_{2}(U) \ar[d] \ar[r] & \mathcal{A}_{1}(V) \ar[d] \\  \mathcal{A}_{2}(U') \ar[r] &  \mathcal{A}_{1}(V') }
\end{equation}
commutes. We use the terms $k$-analytic space and Berkovich analytic spaces interchangeably. Let $\ttAn_{k}$ denote the category of $k$-analytic spaces.
\end{defn}

\begin{defn}\label{quasinet}A quasi-net on a topological space $X$ is a collection $T$ of subsets of $X$ such that for every $x \in X$ there is a subset $T_x \subset T$ with $|T_x| < \infty$ such that $x \in \cap_{V \in T_x}V$ and there is an open set $U \subset X$ such that $x \in U \subset \cup_{V \in T_x} V.$
\end{defn}
\begin{lem}\label{easyquasinet}Any finite set $T=\{V_1,V_2, \dots V_m \}$ of closed subsets of a topological space $X$ which cover $X$ is a quasi-net.
\end{lem}
{\bf Proof.} Given $x \in X,$ consider the subset $T_x \subset T$ defined by those subsets in $T$ which contain $x$.  Then $x \in \cap_{V \in T_x}V$. Let $U=X-\cup_{V \in T-T_x}V,$ this satisfies the required property.
\ \hfill $\Box$

Consider $\ttAfnd^{op}_{k}/X$, the category of affinoid $k$-analytic spaces over $X$. This means that the objects are pairs $(\mathcal{M}(\mathcal{A}),f)$ where $\mathcal{M}(\mathcal{A})$ is an affinoid $k$-analytic space and $f:\mathcal{M}(\mathcal{A}) \to X$ is a morphism of $k$-analytic spaces.  The morphisms from $(\mathcal{M}(\mathcal{A}_1),f_1)$ to $(\mathcal{M}(\mathcal{A}_2),f_2)$ are morphisms to $X$ commuting with the $f_i$.    
\begin{lem}\label{EverythingColimit} Any $k$-analytic space $X$ is a colimit of the category $\ttAfnd^{op}_{k}/X$ when considered as a subcategory of $\ttAn_{k}$.
\end{lem}
{\bf Proof.}
Consider the family $\hat{\tau}$ of all affinoid domains in $|X|$. It is a net and $|X|$ has a maximal $k$-affinoid atlas $\hat{\mathcal{A}}$. For each $V \in \hat{\tau}$, $\hat{\mathcal{A}}$ assigns $\mathcal{M}(\mathcal{A}_{V})\to X.$  In particular, $\hat{\mathcal{A}}$ assigns homeomorphisms $|\mathcal{M}(\mathcal{A}_{V})| \cong V\subset|X|$ and such that these homeomorphisms satisfy obvious compatibilities. For any other $k$-analytic space $X'$ we have by Exercise 3.2.2 of \cite{Ber2009} an isomorphism
\begin{equation}
\Hom(X,X') \rightarrow \eq[\prod_{V \in  \hat{\tau}} \Hom(\mathcal{M}(\mathcal{A}_{V}), X')\rightrightarrows \prod_{(V,W) \in \hat{\tau}^{2}}  \Hom(\mathcal{M}(\mathcal{A}_{V}\widehat{\otimes_{k}}\mathcal{A}_{W}), X')
].\end{equation}
Together with the factorization of this isomorphism as
\footnotesize
\[\Hom(X,X') \to \lim_{\mathcal{M} \in \ttAfnd^{op}_{k}/X}\Hom(\mathcal{M},X') \hookrightarrow  \eq[\prod_{V \in  \hat{\tau}} \Hom(\mathcal{M}(\mathcal{A}_{V}), X')\rightrightarrows \prod_{(V,W) \in \hat{\tau}^{2}}  \Hom(\mathcal{M}(\mathcal{A}_{V}\widehat{\otimes_{k}}\mathcal{A}_{W}),X')]
\]
\normalsize
this implies that the natural morphism 
\[\Hom(X,X') \to \lim_{\mathcal{M} \in \ttAfnd^{op}_{k}/X}\Hom(\mathcal{M},X')
\]
is an isomorphism. Therefore, $X=\colim_{\mathcal{M} \in \ttAfnd^{op}_{k}/X} \mathcal{M}.$
\ \hfill $\Box$
\subsection{From Banach algebraic geometry to Berkovich geometry}\label{BSSOCBA}
Consider the category $\tC=\ttBan_{k}$ for some valuation field $k$. We have shown in Section \ref{BanachSpaces} that it is a closed symmetric monoidal quasi-abelian categories with $\uHom=\uHom_{k}$ and $\ootimes =\wotimes_{k}$ with all finite limits and colimits and enough projectives so that we can do algebraic geometry relative to $\ttBan_{k}$. In particular the categories of affine schemes over it has certain distinguished morphisms and topologies and we have notions of (Archimedean/non-Archimedean) Banach schemes Banach (infinity) stacks and $n$-algebraic Banach stacks over an (Archimedean/non-Archimedean) valuation field $k$.  
In the Archimedean case we could compare this geometry to the geometry of complex varieties covered by Stein compact subsets. However, we focus here on the non-Archimedean case. In this section $k$ will be a non-Archimedean valuation field.

\begin{defn}\label{def:NDiskrelA}
Let $\mathcal{A}$ be a $k$-affinoid algebra. Given $r_{1}, \cdots, r_{n} \in \mathbb{R}$, we can define an $\mA$ algebra 
\begin{equation}\mA\{r_1^{-1}T_1, \dots, r_n^{-1}T_n\}
\end{equation}
as the completion of $\mA[T_1, \dots, T_n]$ with respect to the norm 
\[\|\sum a_{I} T^{I}\|_{r}=\text{max}_{I}\{\|a_{I}\|_{\mA}r^{I}\}.
\] 
\end{defn}
Let $r$ be a real number greater than zero, denote by $\mA_{r}$ the $\mA$ module with norm $\|a\|= r\|a\|_{\mA}$.
\begin{lem}\label{lem:DiskAlg}$\mA\{r_1^{-1}T_1, \dots, r_n^{-1}T_n\}$ is the symmetric algebra (see subsection \ref{Symmetric}) on $\mV=\mA_{r_1} \oplus \cdots \oplus \mA_{r_n}$ in $\ttMod^{\leq 1}(\mA)$. It can also be seen as the filtered colimit in $\ttMod^{\leq 1}(\mA)$ of 
\[\sS^{0}(\mV) \hookrightarrow \sS^{0}(\mV) \oplus \sS^{1}(\mV)  \hookrightarrow \sS^{0}(\mV) \oplus \sS^{1}(\mV) \oplus \sS^{2 }(\mV)\hookrightarrow \cdots 
\]
where 
\[\sS^{m}(\mV) = \{ \sum_{|I|=m} a_{I} T^{I} | a_{I} \in \mA\}\]
equipped with the norm 
\[\|\sum_{|I|=m} a_{I} T^{I}  \|= \max_{|I|=m} \|a_I\| r^{I}.
\]
\end{lem}
{\bf Proof.} Left to the reader.
\ \hfill $\Box$

\begin{lem}\label{lem:FinTypePres} 

A morphism $p: \mathcal{A} \to \mathcal{B}$ in  $\ttComm(\ttBan^{\leq 1}_{k})$ induces a presentation 
\begin{equation}\label{eqn:NiceForm} \mB \cong \mathcal{A}\{\frac{T_1}{r_1}, \dots, \frac{T_g}{r_g}\}/(P_1, \dots, P_r) 
\end{equation}
where $P_{i} \in \mathcal{A}\{\frac{T_1}{r_1}, \dots, \frac{T_g}{r_g}\}$
if and only if $p$ is of finite presentation in $\ttComm(\ttBan^{\leq 1}_{k})$ as defined in Definition \ref{defn:FinitePres}.
\end{lem}

\ \hfill $\Box$

The following two lemmas are technical results that will be used only in the proof of Theorem \ref{thm:localForm}.
\begin{lem}\label{SplitProdStr2}Let $\mathcal{A}, \mathcal{C}$ be $k$-affinoid algebras considered as objects in $\ttComm(\ttBan_{k})$ and let $f:\mathcal{A} \to \mathcal{C}$ be morphism in $\ttComm(\ttBan_{k})$ which is a strict epimorphism when considered in $\ttBan_{k}$ such that there is a $k$-affinoid algebra $\mathcal{B}$ and the post-composition of $f$ with some homotopy epimorphism $g:\mathcal{C} \to \mathcal{B}$ in $\ttComm(\ttBan_{k})$ is a homotopy epimorphism $h:\mA\to \mB$ in $\ttComm(\ttBan_{k})$. Then there is a $k$-affinoid algebra $\mA'$ and a isomorphism $\mathcal{A} \cong \mC \times \mA'$ such that the projection to $\mC$ corresponds to $f$ under this isomorphism.\end{lem}
{\bf Proof.} 
The morphisms 
\[\mA\to \mC \to \mB
\]
induce morphisms on the derived categories 
\[\sD^{-}(\mB) \to \sD^{-}(\mC) \to \sD^{-}(\mA).
\]
Since the composition is fully faithful and the first morphism is as well, the second morphism must be fully faithful and so by Lemma \ref{lem:HomotopyMon} we have $\mC \wotimes^{\mathbb{L}}_{\mA} \mC \cong \mC$. Let $I=\ker (f)$. There is a strict, short exact sequence 
\begin{equation}\label{eqn:idealAC} 0 \to I \to \mA \to \mC \to 0.
\end{equation}
If we consider the derived completed tensor product of (\ref{eqn:idealAC}) over $\mA$ with $\mC$ we find an exact triangle 
\[I \wotimes^{\mathbb{L}}_{\mA} \mC \to \mC \to \mC \wotimes^{\mathbb{L}}_{\mA} \mC
\]
and because the second morphism is an isomorphism, we see that 
$I \wotimes^{\mathbb{L}}_{\mA} \mC$ is isomorphic to $0$. If we now consider the derived completed tensor product over $\mA$ of (\ref{eqn:idealAC}) with $I$ we an exact triangle 
\[I \wotimes^{\mathbb{L}}_{\mA} I \to I \to I \wotimes^{\mathbb{L}}_{\mA} \mC
\]
and so we get an isomorphism $I \wotimes^{\mathbb{L}}_{\mA}I \to I$.
So we have $I=\text{image}[I \wotimes_{\mA}I \to I]$ and in fact this implies that $I= I^2:= \text{image}[I \otimes_{\mA}I \to I]$. Therefore, there exists an element $e \in \mA$ such 
that $e^{2}=e$ and $e\mA =I$. This gives the structure of a $k$-affinoid algebra 
to $I$, which we denote 
by $\mA'=\mA/(1-e)\mA$. Now because $f$ is 
a strict epimorphism, there is a strict short exact sequence 
\[0 \to I \to \mA \stackrel{f}\to \mC \to 0
\]
which in fact is split by the morphism of algebras $e:\mA \to \mA'.$ Therefore, \[(e,f):\mA \to \mA' \times \mC\] is 
an isomorphism.
\ \hfill $\Box$
\begin{lem}\label{CoverOfAff}Let $\mathcal{A}, \mathcal{B}$ be $k$-affinoid algebras and let $f:\mathcal{A} \to \mathcal{B}$ be a morphism in the category of $k$-affinoid algebras with the property that $|\mM(\mA)|$ has a finite covering by affinoid domains $V_{j}$ corresponding to affinoid domain embeddings $\mM(\mA_{V_j})\to \mM(\mA).$ Suppose also that morphisms $\mM(\mA_{V_j}\widehat{\otimes}_{\mA}\mB) \to \mM(\mA_{V_j})$ are affinoid domain embeddings. Then the morphism $\mM(\mB) \to \mM(\mA)$ is an affinoid domain embedding.
\end{lem}
{\bf Proof.}
Let us denote by $U$ the image of $|\mM(\mB)|$ inside $|\mM(\mA)|.$ We have $\overline{U} \cap V_{j} = \overline{U \cap V_j} $ for all $j$. Since $U \cap V_j$ is closed in $V_j$ and $V_j$ is closed in $|\mM(\mA)|$ we see that  $U \cap V_j$ is closed in $|\mM(\mA)|$.  Therefore $\overline{U} \cap V_{j} = U \cap V_{j} $ for all $j$. Hence $\overline{U}=U$ and so $U$ is closed inside $|\mM(\mA)|.$
Let $\mA \to \mC$ be a bounded homomorphism of affinoid $k$-algebras such that the image of $|\mM(\mC)|$ lies in $U$.  We wish to show that the morphism $\mM(\mC) \to \mM(\mA)$ factors through a morphism $\mM(\mC) \to \mM(\mB)$. Notice that $\mA_{V_j} \to \mC \widehat{\otimes}_{\mA}\mA_{V_j} $ is a bounded homomorphism of affinoid $k$-algebras such that the image of $\mM(\mC \widehat{\otimes}_{\mA}\mA_{V_j} )$ lies in $U\cap V_j$. Therefore, the morphisms $\mM(\mC \widehat{\otimes}_{\mA}\mA_{V_j} ) \to \mM(\mA_{V_j} )$ factor in a unique way through morphisms  $\mM(\mC \widehat{\otimes}_{\mA}\mA_{V_j} )  \to \mM(\mA_{V_j}\widehat{\otimes}_{\mA}\mB).$ When thought of as morphisms   $\mM(\mC \widehat{\otimes}_{\mA}\mA_{V_j} )  \to \mM(\mB)$ they agree when pulled back to $\mM(\mC \widehat{\otimes}_{\mA}\mA_{V_j \cap V_k} ). $ The preimages of $V_j$ in $|\mM(\mC)|$ are analytic domains by \cite{Te} Exercise 3.2.2 (v). These preimages are the pullback of a quasi-net and therefore form a quasi-net by Lemma \ref{easyquasinet} and therefore by Exercise 3.2.2 (v) of \cite{Ber2009} we have a unique morphism $\mM(\mC)  \to \mM(\mB)$ which restricts to the morphisms $\mM(\mC \widehat{\otimes}_{\mA}\mA_{V_j} )  \to \mM(\mB).$ 
Indeed, this follows from the commutative diagram of exact sequences 
\footnotesize
\begin{equation}
\xymatrix{
0 \ar[r] & \Hom(\mM(\mC),\mM(\mB)) \ar[d] \ar[r]  & \prod_{j}\Hom(\mM(\mC\widehat{\otimes}_{\mA}\mA_{V_j} ),\mM(\mB)) \ar[r] \ar[d] & 
\prod_{j,k}\Hom(\mM(\mC\widehat{\otimes}_{\mA}\mA_{V_j\cap V_k} ),\mM(\mB)) \ar[d] \\ 
0 \ar[r] &\Hom(\mM(\mC),\mM(\mA)) \ar[r] & \prod_{j}\Hom(\mM(\mC\widehat{\otimes}_{\mA}\mA_{V_j} ),\mM(\mA)) \ar[r] &  \prod_{j,k}\Hom(\mM(\mC\widehat{\otimes}_{\mA}\mA_{V_j\cap V_k} ),\mM(\mA)). }
\end{equation}
\normalsize
This clearly provides the required factorisation. 
\ \hfill $\Box$

If also, the $\mM(\mA_{V_j})$ are rational in $\mM(\mA)$ and  $\mM(\mA_{V_j}\widehat{\otimes}_{\mA}\mB)$ is rational in $\mM(\mA_{V_j})$ notice that $\mM(\mB)$ is a union of the rational domains $\mM(\mA_{V_j}\widehat{\otimes}_{\mA}\mB)$ in $\mM(\mA).$
 
\begin{lem}\label{lem:NonExpFP} Let $\mathcal{A}$ be a $k$-affinoid algebra and suppose there is a morphism $f: \mA \to \mB$ of finite presentation in $\ttComm(\ttBan^{\leq 1}_{k}).$ Then $\mB$ is a $k$-affinoid algebra.
\end{lem}
{\bf Proof.}
By combining a presentation for $\mB$ over $\mA$ and a presentation for $\mA$ over $k$ one can write $\mB$ as a finite colimit of objects of finite presentation in $\ttComm(\ttBan^{\leq 1}_{k})$. Therefore, $\mB$ has finite presentation in $\ttComm(\ttBan^{\leq 1}_{k})$.
\ \hfill $\Box$
\begin{thm}\label{thm:localForm}Let $\mathcal{A}, \mathcal{B}$ be $k$-affinoid algebras and let $f:\mathcal{A} \to \mathcal{B}$ be a morphism in the category of $k$-affinoid algebras. Assume that $f$ is a homotopy epimorphism (see Definition \ref{defn:homotopyEpi}) when considered in the category $\ttComm(\ttBan_{k})$. Then the morphism $\mM(\mB) \to \mM(\mA)$ corresponding to $f$ is an affinoid domain embedding.\end{thm}
{\bf Proof.} 
We refer here to Temkin's proof \cite{Te2} of the Gerritzen-Grauert Theorem for morphisms of affinoid algebras.  This theorem, assuming only the epimorphism condition on $f$, produces a finite collection of morphisms of $k$-affinoid algebas $\mA \to \mA_{V_i}$ corresponding to rational domain embeddings $\mM(\mA_{V_i}) \to \mM(\mA)$ covering $|\mM(A)|$ with the images $V_i.$  The theorem further ensures that the morphisms $\mA_{V_i} \to \mB \widehat{\otimes}_{\mA} \mA_{V_i}$ induced from $f$ admit factorzations \[\mA_{V_i} \twoheadrightarrow \mC_{i} \hookrightarrow \mB \widehat{\otimes}_{\mA} \mA_{V_i}.\] These factorizations correspond to the composition of the morphism of $k$-affinoid algebras $\mC_{i} \to (\mC_{i})_{W_i}=\mB \widehat{\otimes}_{\mA} \mA_{V_i}$ corresponding to Weierstrass domain embeddings  with the surjective morphisms of affinoid $k$-algebras $\mA_{V_i} \to \mC_{i}$ corresponding to closed immersions. 
Therefore, by Lemma \ref{lem:LCSrHoEpis} the morphism  $\mC_{i} \to \mB \widehat{\otimes}_{\mA} \mA_{V_i}$ is a homotopy epimorphism in the category $\ttComm(\ttBan_{k})$. Notice that 
\[(\mB\wotimes^{\mathbb{L}}_{\mA}\mA_{V_i})\wotimes^{\mathbb{L}}_{\mA_{V_i}}(\mB\wotimes^{\mathbb{L}}_{\mA}\mA_{V_i}) \cong \mB\wotimes^{\mathbb{L}}_{\mA}\mA_{V_i}
\]
because homotopy epimorphisms are closed under derived base change. However, applying Lemma \ref{lem:LCSrHoEpis} we have 
\[\mB\wotimes^{\mathbb{L}}_{\mA}\mA_{V_i} \cong \mB\wotimes_{\mA}\mA_{V_i}
\]
and so we see that 
\[(\mB\wotimes_{\mA}\mA_{V_i})\wotimes^{\mathbb{L}}_{\mA_{V_i}}(\mB\wotimes_{\mA}\mA_{V_i}) \cong \mB\wotimes^{\mathbb{L}}_{\mA}\mA_{V_i}
\]
and so by Lemma \ref{lem:HomotopyMon} the morphisms $\mA_{V_i} \to \mB \widehat{\otimes}_{\mA} \mA_{V_i}$ are homotopy epimorphisms. Also, the morphisms of affinoid algebras corresponding to Weierstrass domain embeddings are injective. Therefore, Lemma \ref{SplitProdStr2} can applied by choosing the $f$ from that lemma to be the morphism $\mA_{V_i} \to \mC_{i}$ and $g$ from that lemma to be the morphism $\mC_{i} \to \mB \widehat{\otimes}_{\mA} \mA_{V_i}.$ The lemma then tells us that the morphism $\mM(\mC_i) \to \mM(\mA_{V_i})$ is simply the inclusion of a connected component in a disjoint union of affinoids. Therefore,  $\mM(\mC_i) \to \mM(\mA_{V_i})$ is an affinoid domain embedding. Because the composition of affinoid domain embeddings is an affinoid domain embedding, we conclude that the morphisms $\mM(\mB \widehat{\otimes}_{\mA} \mA_{V_i})\to \mM(\mA_{V_i})$ are affinoid domain embeddings as well. By Lemma \ref{CoverOfAff} we conclude that the original morphism gives an affinoid domain embedding $\mM(\mB) \to \mM(\mA)$.

\ \hfill $\Box$


From now on we write $\ttMod^{RR}(\mathcal{A})$ in place of $\ttMod_{\ttAfnd_{k}^{op}}^{RR}(\mathcal{A})$.
\begin{lem}\label{lem:ConsImpliesSur}
Let $\mathcal{A}$ be a $k$-affinoid algebra. Let $\{f_{i}:\mathcal{A} \to \mathcal{A}_{V_i}\}_{i \in I}$ be a family of affinoid localizations such that for some finite set $J \subset I$ the corresponding family of functors 
\[\ttMod^{RR}(\mA) \to \ttMod^{RR}(\mA_{V_i})
\]
for $i \in J$ is conservative. Then the morphism $\coprod_{i\in J}\mM(\mA_{V_i}) \to \mM(\mA)$ is surjective.
\end{lem}
{\bf Proof.}
We argue by contradiction.  First assume that $k$ is non-trivially valued and $\mathcal{A}$ is strictly affinoid. Suppose that the family of functors is conservative and some point $x\in \mM(\mA)$ is not in the image. By Proposition 2.1.15 of \cite{Ber1990} the subset of points of $y\in \mM(\mA)$ such that $\ker(|\;|_{y})$ is a maximal ideal is a dense subset of $\mM(\mA)$. Therefore, since the image is the closed set $\cup_{i\in J}V_{i}$ we may assume (by changing the point $x$) that $x\in \mM(\mA)$ is not in the image and $\ker(|\;|_{x})$ is a maximal ideal. Chose using Proposition 2.2.3 (iii) of \cite{Ber1990} an affinoid subdomain $W$ of $\mM(\mA)$ such that $x \in W$ and $W \cap V_{i}$ is empty for all $i \in J$.  Consider the 
 morphism $0 \to \mA_W$ of $\ttMod^{RR}(\mA)$. It is not an isomorphism but for each $i \in J$, the pullback to each $\spec(\mA_{V_i})$ is the isomorphism $0 \to 0=\mA_{W} \wotimes_{A} \mA_{V_i}$ of $\ttMod^{RR}(\mA_{V_i})$. This gives a contradiction. For the general case, choose using Proposition 2.1.2 of \cite{Ber1990}, a valuation field extension $k \to K$ such that the valuation on $K$ is non-trivial and $\mA\wotimes_{k}K$ is a strictly $K$-affinoid algebra. Notice that the conservativity assumption on the original family implies by Lemma \ref{cor:BaseChangeConservHzAR} applied to the base change $\spec(\mA\wotimes_{k}K) \to \spec(\mA)$ that the family of functors $\{\Mod^{RR}(\mA\wotimes_{k}K) \to \Mod^{RR}(\mA_{V_i}\wotimes_{k}K)\}_{i \in J}$ is also conservative. The morphism $\coprod_{i\in J}\mM(\mA_{V_i}\wotimes_{k} K) \to \mM(\mA\wotimes_{k}K)$ cannot be surjective because in the commutative diagram, 
\[
\xymatrix{ \coprod_{i\in J}\mM(\mA_{V_i}\wotimes_{k}K) \ar[d] \ar[r] & \coprod_{i\in J}\mM(\mA_{V_i}) \ar[d] \\  \mM(\mA\wotimes_{k}K) \ar[r] &  \mM(\mA) }\]
the horizonal arrows are surjective.  Therefore, we have reduced to the previous case and so the proof is complete.
\ \hfill $\Box$

\begin{lem}\label{lem:AlternatingTate}Consider a (surjective) cover of $X= \mM(\mA)$ by a finite collection of affinoid domains $V_i= \mM(A_{V_{i}})$. Then the complex 
\[
0 \to \mA \to \prod_{i_1} \mA_{V_{i_1}} \to \prod_{i_1 < i_2} \mA_{V_{i_1}} \wotimes_{\mA} \mA_{V_{i_2}} \to \cdots \to  \mA_{V_{1}} \wotimes_{\mA} \mA_{V_{2}}\wotimes_{\mA} \cdots \wotimes_{\mA}  \mA_{V_{n}} \to 0. 
\]
is strictly exact.
\end{lem}
{\bf Proof.} By Proposition 1, section 8.1 of \cite{BGR}, the inclusion (which is strict) of alternating cochains inside all cochains is a quasi-isomorphism. On the other hand, the complex of all cochains is strictly exact by Proposition 2.2.5 of \cite{Ber1990}. The conclusion follows.
\ \hfill $\Box$
\begin{lem}\label{lem:CoversRconservative}Consider a (surjective) cover of $X= \mM(\mA)$ by a finite collection of affinoid domains $V_i= \mM(A_{V_{i}})$. Then the corresponding family of functors $\Mod^{RR}(\mA) \to \Mod^{RR}(\mA_{V_i})$ is conservative.
\end{lem}
{\bf Proof.}
Let $f: \mM \to \mN$ in $\ttMod^{RR}(\mA)$ be any morphism such that $f_{i}: \mM \wotimes_{\mA}\mA_{V_i} \to \mN \wotimes_{\mA}\mA_{V_i} $ are isomorphisms for all $i$. 
The alternating version of the \v{C}ech-Amitsur complex (see Definition \ref{defn:Amitsur}) corresponding to the morphism $\mA \to \prod_{i=1}^{n} \mA_{V_i}$ is a strictly exact bounded above complex by \ref{lem:AlternatingTate} and so defines an element of $D^{-}(\mA)$ (in fact the $0$ element!).
\[0 \to \mA \to \prod_{i_1} \mA_{V_{i_1}} \to \prod_{i_1 < i_2} \mA_{V_{i_1}} \wotimes_{\mA} \mA_{V_{i_2}} \to \cdots \to  \mA_{V_{1}} \wotimes_{\mA} \mA_{V_{2}}\wotimes_{\mA} \cdots \wotimes_{\mA}  \mA_{V_{n}} \to 0. 
\]
Each object in this complex is acyclic for the functor $\mM \wotimes_{\mA}(-)$ since $\mM$ (being RR-quasi-coherent) is transversal to localizations of $\mA$. Therefore by Proposition \ref{prop:ProjFProj}, if we apply the derived functor $\mM \wotimes_{\mA}^{\mathbb{L}}(-)$ we are left with a strictly exact complex 
\[0 \to \mM \to \prod_{i_1} \mM \wotimes_{\mA} \mA_{V_{i_1}} \to \prod_{i_1<i_2}\mM \wotimes_{\mA} \mA_{V_{i_1}} \wotimes_{\mA} \mA_{V_{i_2}} \to \cdots \to \mM \wotimes_{\mA}\mA_{V_{1}} \wotimes_{\mA} \mA_{V_{2}}\wotimes_{\mA} \cdots \wotimes_{\mA}  \mA_{V_{n}} \to 0.
\]
We can do the same thing for $\mN$. The $f_{i}$ extend uniquely to morphisms of the complexes resolving $\mM$ and $\mN$. Therefore $f$ is an isomorphism. Conversely, if $f$ is an isomorphism, the $f_i$ obviously are as well.
\ \hfill $\Box$
\begin{rem}\label{rem:TateAcylicity} In the proof of Lemma \ref{lem:CoversRconservative} we showed the interesting fact that for any $\mM \in \ttMod^{RR}(\mA)$, and any finite cover $\coprod_{i} \mM(\mA_{V_i}) \to \mM(\mA)$ the complex  
\[0 \to \mM \to \prod_{i_1} \mM \wotimes_{\mA} \mA_{V_{i_1}} \to \prod_{i_1<i_2}\mM \wotimes_{\mA} \mA_{V_{i_1}} \wotimes_{\mA} \mA_{V_{i_2}} \to \cdots \to \mM \wotimes_{\mA}\mA_{V_{1}} \wotimes_{\mA} \mA_{V_{2}}\wotimes_{\mA} \cdots \wotimes_{\mA}  \mA_{V_{n}} \to 0.
\]
is strictly exact (a version of Tate acyclicity). 
\end{rem}
\begin{example}This example was suggested via a correspondence with V. Berkovich. The analogue of Corollary 4.48 is false for finite modules since they are not preserved by push-forward. On the other hand,  RR-quasi-coherent modules are preserved by push-forward (see Lemma \ref{lem:RRpreserve}) and do detect surjectivity. If $k$ is a field equipped with the trivial valuation and $\mA= k\{\frac{x}{r}\}$ for $0<r<1$ and $\mA_{V}= k\{\frac{x}{r'}\}$ for $0<r'<r$ then $\mM(\mA_{V})\to \mM(\mA)$ is not surjective while the natural morphism $\mA \to \mA_{V}$ induces an equivalence $\ttMod^{f}(\mA)\to \ttMod^{f}(\mA_{V})$, as both categories are equivalent to the category of finite modules in the algebraic sense over $k[[x]]$.  On the other hand, it is easy to find objects of $\ttMod^{RR}(\mA)$ which go to zero in $\ttMod^{RR}(\mA_{V})$. For instance, choose a valued extension $K$ of $k$ such that $K$ is non-trivialy valued and $\mA\wotimes_{k}K$ is strictly affinoid. Chose a point $x \in \mM(\mA\wotimes_{k}K) - \mM(\mA_{V}\wotimes_{k}K)$ such that $\ker(|\ |_{x})$ is a closed maximal ideal of $\mA\wotimes_{k}K$. Choose using Proposition 2.2.3 (iii) of \cite{Ber1990} an affinoid subdomain $W$ of $\mM(\mA\wotimes_{k}K)$ such that $W$ does not intersect $\mM(\mA_{V}\wotimes_{k}K)$ and $x \in W$. Then we have the non-zero object $(\mA\wotimes_{k}K)_{W}$ of $\ttMod^{RR}(\mA)$. It is the push-forward of $(\mA\wotimes_{k}K)_{W} \in \ttMod^{RR}(\mA\wotimes_{k}K)$ along $\spec(\mA\wotimes_{k}K) \to \spec(\mA)$ and therefore RR-quasicoherent by Lemma \ref{lem:RRpreserve}.  It goes to zero after applying the functor $(-)\wotimes_{\mA}A_{V}$ because using (\ref{eqn:basechangeequation}) we have \[(\mA\wotimes_{k}K)_{W}\wotimes_{\mA}A_{V} \cong (\mA\wotimes_{k}K)_{W} \wotimes_{\mA\wotimes_{k}K}(\mA_{V}\wotimes_{k}K)\cong 0.
\]
\end{example}
\subsection{Topologies in the Banach algebraic geometry setting}

\begin{thm}\label{thm:covers}Consider the Berkovich space $\mathcal{M}(\mA)$ for $\mA\in \ttAfnd_{k}$. The covers in the weak $G$-topology on $\mathcal{M}(\mA)$, as defined on page 30 of \cite{Ber1990}, coincide precisely with the covers of $\spec(\mA)$ by affinoids in the formal homotopy Zariski topology on the homotopy transversal subcategory $\ttAfnd_{k}^{op} \subset\ttAff(\ttBan_{k})$ (using the terminology of Definition \ref{defn:fhTVZ} and Proposition  \ref{hfTV1}).
\end{thm}
{\bf Proof.}
Say that we are given a cover of $\mM(\mA)$ by rational domains which has a finite subcover.  Lemma \ref{lem:LCSrHoEpis} tells us that the affinoid domains correspond to homotopy epimorphisms $\mathcal{A} \to \mathcal{B}_{i}$ for some $i \in I$ in $\ttAfnd^{op}_{k}$. Lemma \ref{lem:CoversRconservative} tells us that the family of functors indexed by the finite subset $J \subset I$ corresponding to the finite subcover is conservative. In the other direction suppose we are given a cover of $X$ in the formal homotopy Zariski topology. It has a finite conservative sub-cover which must be surjective by Lemma \ref{lem:ConsImpliesSur}. Every element of it must be a subset of $\mM(\mA)$  because every morphism in the cover is a monomorphism. In fact, Theorem \ref{thm:localForm} implies that every morphism in the cover is an affinoid domain embedding. This concludes the proof.
\ \hfill $\Box$

\begin{lem}\label{lem:IsSheaf}
Let $\mathcal{A}, \mathcal{B}_{i}$ be $k$-affinoid algebras and let $\{f_{i}:\mathcal{A} \to \mathcal{B}_{i}\}_{i\in I}$ be a family of morphisms of $k$-affinoid algebras which is a formal homotopy Zariski open cover in the category $\ttAfnd_{k}^{op}$. Chose a finite subset $J \subset I$ with the conservative property. Then for any $k$-analytic space $X$ we have 
\begin{equation}\label{equation:gluingdiagram}\Hom(\mM(\mA),X)=\eq [\prod_{i \in J}\Hom(\mM(\mathcal{B}_{i}),X) \rightrightarrows \prod_{i \in J,j \in J}\Hom(\mM(\mB_{i} \widehat{\otimes}_{\mA} \mB_{j}), X)].
\end{equation}
Therefore, the presheaf $\straighth_{X}$ is a sheaf in the 
formal homotopy Zariski topology on the homotopy transversal subcategory $\ttAfnd_{k}^{op} \subset\ttAff(\ttBan_{k})$.
\end{lem}
{\bf Proof.}
Lemma \ref{lem:ConsImpliesSur} and Theorem \ref{thm:localForm} imply that the family of morphisms indexed by $J$ corresponds to a finite cover by affinoid domains. Lemma \ref{easyquasinet} implies that it is a quasinet. Exercise 3.2.2 (v) of \cite{Ber2009} now tells us that Equation \ref{equation:gluingdiagram} is valid.

\ \hfill $\Box$
\noindent

Unlike the rest of this article, the following theorem is only significant given the results in \cite{BeKr2}.
\begin{thm}\label{thm:BigRestriction} The topology of Definition \ref{defn:TVtopologyIntro} restricted to the category $\ttAfnd_{k}^{op}$ agrees with the weak G-topology. 
\end{thm}

{\bf Proof.}

In light of Theorems \ref{thm:TensLocs} and \ref{thm:localForm} we only need to understand why condition (2) of Definition \ref{defn:TVtopologyIntro} is equivalent to a finite collection of affinoid subdomains $V_i$ covering every point of an affinoid $\mM(\mA)$. Suppose first that the union of the $V_i$ is all of $\mM(\mA)$. Say we have a morphism $\mM \to \mN$ in $D^{-}(\mA)$ and we know that the induced morphism $\mM\wotimes^{\mathbb{L}}_{\mA}\mA_{V_i}\to \mN\wotimes^{\mathbb{L}}_{\mA}\mA_{V_i}$ is an isomorphism in  $D^{-}(\mA)$ for all $i$. For any $i_1 < i_2 < \cdots < i_p$ we get an isomorphism $\mM\wotimes^{\mathbb{L}}_{\mA}\mA_{V_{i_1, \dots, i_p}}\to \mN\wotimes^{\mathbb{L}}_{\mA}\mA_{V_{i_1, \dots, i_p}}$ where \[\mA_{V_{i_1, \dots, i_p}}= \mA_{V_{i_1}} \wotimes_{\mA} \cdots  \wotimes_{\mA} \mA_{V_{i_1}}= \mA_{V_{i_1}} \wotimes^{\mathbb{L}}_{\mA} \cdots  \wotimes^{\mathbb{L}}_{\mA} \mA_{V_{i_1}}.\] Therefore, we get an isomorphism from 
\begin{equation}\label{eqn:Mres}\mM\wotimes^{\mathbb{L}}_{\mA}( \prod_{i_1 } \mA_{V_{i_1}} \to \prod_{i_1 < i_2} \mA_{V_{i_1, i_2}}  \to \cdots \to  \mA_{V_{1,2,...,n}}  \to 0) 
\end{equation}
to
\begin{equation}\label{eqn:Nres}\mN\wotimes^{\mathbb{L}}_{\mA}( \prod_{i_1 } \mA_{V_{i_1}} \to \prod_{i_1 < i_2} \mA_{V_{i_1, i_2}}  \to \cdots \to  \mA_{V_{1,2,...,n}}  \to 0) 
\end{equation}
in $D^{-}(\mA)$. Since by Lemma \ref{lem:AlternatingTate} the natural morphism from $\mA$ 
to 
\[
\prod_{i_1 } \mA_{V_{i_1}} \to \prod_{i_1 < i_2} \mA_{V_{i_1, i_2}}  \to \cdots \to  \mA_{V_{1,2,...,n}}  \to 0
\]
is an isomorphism in $D^{-}(\mA)$,  (\ref{eqn:Mres}) and (\ref{eqn:Nres}) are equivalent to $\mM$ and $\mN$ respectively. Therefore, we can conclude that the original morphism $\mM \to \mN$ is an isomorphism. Conversely, suppose that there is a point $x\in \mM(\mA)$ which is not in any of the $V_{i}$. Suppose first that $k$ is non-trivially valued and that $\mA$ is strictly $k$-affinoid. Chose as in Lemma \ref{lem:ConsImpliesSur} a subdomain $W$ of $\mM(\mA)$ such that $x \in W$ and $W$ does not intersect any of the $V_i$. Then $\mA_{W}\to 0$ is a morphism in $D^{-}(\mA)$ but for each $V_{i}$ the derived pullback gives an isomorphism $\mA_{W} \wotimes_{A}^{\mathbb{L}}\mA_{V_{i}} \cong  \mA_{W} \wotimes_{A}\mA_{V_{i}} \cong 0 \to 0$. For the general case, choose using Proposition 2.1.2 of \cite{Ber1990}, a valuation field extension $k \to K$ such that the valuation on $K$ is non-trivial and $\mA\wotimes_{k}K$ is a strictly $K$-affinoid algebra. Notice that the conservativity assumption on the original family implies by Lemma \ref{lem:DerivedBaseChangeConserv} applied to the base change $\spec(\mA\wotimes_{k}K) \to \spec(\mA)$ that the family of functors $\{D^{-}(\mA\wotimes_{k}K) \to D^{-}(\mA_{V_i}\wotimes_{k}K)\}_{i \in J}$ is also conservative. The morphism $\coprod_{i\in J}\mM(\mA_{V_i}\wotimes_{k} K) \to \mM(\mA\wotimes_{k}K)$ cannot be surjective because in the commutative diagram, 
\[
\xymatrix{ \coprod_{i\in J}\mM(\mA_{V_i}\wotimes_{k}K) \ar[d] \ar[r] & \coprod_{i\in J}\mM(\mA_{V_i}) \ar[d] \\  \mM(\mA\wotimes_{k}K) \ar[r] &  \mM(\mA) }\]
the horizonal arrows are surjective.  Therefore, we have reduced to the previous case and so the proof is complete. 
\ \hfill $\Box$

For the next lemma, we will need the notion of continuous and cocontinuous functors. These are defined in the appendix in Definition \ref{defn:ContCocont}. 
\begin{lem}\label{lem:inclprops}The inclusion functor 

\[\ttAfnd^{op}_{k} \to \ttAff(\ttBan_{k}) 
\]
is fully faithful and continuous with respect to the Zariski, formal Zariski, homotopy Zariski, formal homotopy Zariski and fpqc topologies.
\end{lem}
{\bf Proof.}
Most of these properties are more or less obvious at this point. The fiber products correspond to completed tensor products of affinoid algebras or commutative algebra objects in $\ttBan_{k}$. The fact that these functors take covers to covers in the various topologies follows from the definition of covers on the domain of the functors. The fully faithful property is discussed in \ref{lem:extension} of the appendix.
  
\hfill $\Box$

Recall that $\ttAn_{k}$ is the category of Berkovich spaces over a non-Archimedean valuation field $k$.
\begin{thm}\label{thm:embedding} There is a fully faithful embedding from $\ttAn_{k}$ to the categories of presheaves of sets on $\ttAfnd^{op}_{k}$ and $\ttComm(\ttBan_{k})^{op}.$
This embedding induces a fully faithful embedding 
\[\ttAn_{k} \hookrightarrow  \ttSch((\ttAfnd^{op}_{k})^{fhZar}) 
\].
For any affinoid algebras $\mA$ and $\mB$, the scheme assigned to
$\mathcal{M}(\mA)$ evaluates on $\spec(\mB)$ to $\Hom_{\ttAfnd_{k}}(\mA,\mB)$.
\end{thm}
{\bf Proof.}
Let $X$ be a $k$-analytic space. Consider the contravariant functors 
\[\straighth_{X}:\ttAfnd^{op}_{k} \to \ttSet
\]
defined by 
\[\straighth_{X}(\spec(\mA))=\Hom_{\ttAn_{k}}(\mathcal{M}(\mA),X)
\]
for $\mA \in \ttAfnd_{k}$.
They form a functor
\begin{equation}\label{equation:AnToPre}\ttAn_{k} \stackrel{\straighth}\to \ttPr(\ttAfnd^{op}_{k})
\end{equation}
which can be seen the Yoneda embedding for $\ttAn_{k}$ followed by the restriction functor $\ttPr(\ttAn_{k}) \to \ttPr(\ttAfnd^{op}_{k})$. 
We would like to know that for every $X_{1},X_{2} \in \ttAn_{k}$ that 
\[\Hom_{\ttAn_{k}}(X_{1},X_{2}) \to \Hom_{\ttPr(\ttAfnd^{op}_{k})}(\straighth_{X_1},\straighth_{X_2})
\]
is an isomorphism. The fact that the functor $\straighth$ is fully faithfull follows immediately from Lemma \ref{EverythingColimit} which says that $X$ is the final object in the category of morphisms $\spec(\mA) \to X$ where $\mA \in \ttAfnd_{k}$ and the fact that for any full subcategory $\ttC \subset \ttD$ such that every object $d$ in $\ttD$ is a colimit of objects in $\ttC$ mapping to $d$, the functor $\straighth$ embeds $\ttC$ fully and faithfully into the category of presheaves of sets on $\ttD$.  This fact can be found in \cite{SGA4}, Exp I, Prop 7.2. 
Consider the fully faithful, continuous and cocontinuous functor from Lemma \ref{lem:inclprops}
\[\su: \ttAfnd^{op}_{k} \to \ttAff(\ttBan_{k}). 
\]

which is left adjoint to the restriction functor. 

It is clear by the definition of Berkovich analytic spaces that we actually get a fully faithful embedding 
\[\ttAn_{k} \hookrightarrow  \ttSch((\ttAfnd^{op}_{k})^{fhZar}) 
\]
and when we restrict the scheme to the subcategory $\ttAfnd_{k}^{op}$ this assignment agrees with the standard functor of points in the category $\ttAn_{k}$.
\ \hfill $\Box$

\begin{rem}
Note that exactly the same argument shows that the category of rigid analytic spaces (which contains $k$-analytic spaces as a full subcategory) embeds fully faithfully into $\ttSch((\ttAfnd^{op}_{k})^{fhZar}).$ 
\end{rem}

\begin{rem}
In order to embed $k$-analytic spaces and rigid analytic spaces into a catgeory of Banach schemes, $\ttSch(\ttBan_{k})$, we need to work in the derived setting as the homotopy Zariski topology is only well defined on simplicial Banach algebras. As $\ttAfnd_{k}$ is homotopy Zariski transversal this topology restricts to the non-derived site of affinoids. The derived approach is pursued in \cite{BeKr2}.  
\end{rem}

\subsection{Huber Points}

Let $X$ be an object of $\ttSch((\ttAfnd_{k})^{op})^{hZar})$. Define a category $\hZar(X)$ as the full sub-category of $\ttSh((\ttAfnd_{k})^{op})^{hZar})/X$ whose objects are $u:Y\to X$ with $u$ a homotopy Zariski open immersion. This category 
is a locale and inherits a topology: a family $\{Y_i\to Y\}$ is a covering if and only if $\coprod Y_i\to Y$ is an epimorphism of sheaves. We have the full subcategory $\hZarAff(X)$ whose objects are $u:Y\to X$ where $Y$ is affine. 
The continuous inclusion $\hZarAff(X)\to \hZar(X)$ induces an equivalence of categories of sheaves:
$\ttSh(\hZar(X))\cong \ttSh(\hZarAff(X))$. 

The fact that $\hZar(X)$ is a locale and that the topology on it 
is generated by a quasi-compact  pre-topology (as covering families in $\ttSh(\hZarAff(X))$ are finite), implies that 
$\hZar(X)$ is equivalent as a locale to the locale of open subsets of a topological space $|X|$. Using this we can view $X$ a $\ttComm(\ttBan_{k})$-valued ringed space:
$(|X|,\mathcal{O}_X)$, where $\mathcal{O}_X$ is a sheaf valued in $\ttComm(\ttBan_{k}).$

\begin{lem}
Let $X$ be a Berkovich space. $\hZarAff(X)$ is the rigid analytic site of $X$.
\end{lem}

{\bf Proof.} This is a rephrasing of Theorem \ref{thm:covers}.
\hfill $\Box$

\begin{thm}
Let $X$ be a Berkovich space. The space $|X|$ is the space of Huber points.
\end{thm}
{\bf Proof.} From the previous lemma we know that the locale
corresponding to $X$ is the rigid analytic one. From Huber \cite{Hu}, we see that the topological space corresponding to this locale is the space of Huber points of $X$.
\hfill $\Box$

\section{Work in progress}
In future work, we intend to embed rigid analytic spaces and Huber spaces into the categories of schemes over the opposite category to Banach algebras. We also intend develop a variant which will work with complete, convex bornological algebras, instead of Banach algebras or Fr\'echet algebras or other types of topological vector spaces. Similarly to Paugam \cite{Pa} we would like to handle the Archimedean and non-Archimedean cases with a single language. Also with the same goal, Bambozzi has defined and studied affinoid bornological algebras and dagger algebras in his thesis \cite{Bam}. The categories based on convex bornological spaces are nicer than those based on $\ttBan_{k}$ in that they have arbitrary limits and colimits. The bornological categories contain the categories based on Banach spaces or Fr\'echet spaces via fully faithful embeddings. A lot of work in this direction of bornological geometry and representation theory has been done by Meyer \cite{M} in the Archimedean case over $\mathbb{C}$ and by Houzel and others \cite{H,H2} in the case of a valuation field where the valuation is not the trivial valuation. The extension to the case of a trivially valued field is more difficult, especially in the definition of completeness. Therefore, we could have done our constructions for convex bornological vector spaces over an arbitrary field or for complete convex bornological vector spaces over a non-discretely valued field.  In another direction we think that categories of topological spaces, topological manifolds, differentiable manifolds and real analytic manifolds could also be described in terms of this style of algebraic geometry. This could fit in well with the algebraic analysis of Kashiwara, Schapira, Saito and others. Non-commutative versions may be possible as well along the lines of \cite{So,KR}. It would also be nice to be able to work over rings instead of fields, for instance the rings of $p$-adic intgers $\mathbb{Z}_{p}$ and the integers equipped with either the discrete or the standard norm and do types of analytic geometry over these rings. We are developing a formalism of derived analytic stacks along the lines of following To\"{e}n and Vezzosi's approach \cite{TVe} or Lurie's approach \cite{L1,L2}. It will work simultaneously in the Archimedean and non-Archimedean setting. In particular we need to define a cotangent complex and the notion of smooth and \'{e}tale morphisms and compare these notions with the notions defined by Berkovich.
We also believe that a non-Archimedean versions of both Lurie's Tanakian duality theorem \cite{L3} and (after developing differential graded categories of modules) the algebrization theorem \cite{TV2} of To\"{e}n and Vaqui\'{e} should be provable in our framework. We plan to prove a version of the Beilinson-Bernstein localization theorem over certain non-Archimedean fields in our context which would be a generalisation of the work of \cite{AW} over $\mathbb{Q}_{p}$. Finally, we believe that when working over a finite field with trivial valuation, our formalism will be useful for studying $G$-bundles over an algebraic surface or Kac-Moody bundles over an algebraic curve with an eye towards representation theory as in \cite{BrKa}, \cite{K}, \cite{GP}. This is based on the relationship between $G$-bundles on a punctured formal neighborhood of a curve in surface and twisted $G((t))$-bundles on the curve. 
\appendix
\section{Semi-normed spaces and Banach spaces}\label{SNrmSpaces}
Let $k$ be any field.  We use $\ttVect_{k}$ to denote the category of $k$-modules and $k$-linear morphisms. This is a closed symmetric monoidal category which has all limits and colimits. The set of morphisms will be denoted $\Hom_{k}( \cdot ,\cdot )$ and we use the same notation when considering it as a $k$-module in the natural way. The tensor product will be denoted $\otimes_{k}.$ We will often consider a valuation field: a field equipped with a multiplicitive map (called a norm) $|\;|:k \to \mathbb{R}_{\geq 0}$ such that $|a+b| \leq |a|+|b|$ for all $a,b \in k$ and $|a|=0$ if and only if $a=0$ and such that the field is complete with respect to the (metric defined by the) norm. Such a field together with its norm is called a valuation field. We consider two types of valuation fields. 
\begin{defn}
A non-Archimedean field is a valuation field $k$ such that $|c_1+c_2| \leq \max\{|c_1|, |c_2|\}$ for any $c_1,c_2 \in k$.  Any valuation field that which is not non-Archimedean is called an Archimdedean field. 
\end{defn}

\begin{rem}Any field can be considered a non-Archimedean field by equipping it with the trivial valuation which is defined by $|k-\{0\}|=1$. In the Archimedean case it is know that the only examples are $\mathbb{R}$ or $\mathbb{C}$ equipped with the norm $|\;|^{\epsilon}$ where $0<\epsilon \leq 1$. 
\end{rem}

\begin{defn}\label{defn:ArchimedeanNormedSpace}
An Archimedean semi-normed space over an Archimedean field $k$ is a $k$-module $V$ together with a map $\| \|:V \to \mathbb{R}_{\geq 0}$ which satisfies 

\begin{itemize}
\item $\|cv\|=|c|\|v\|$
\item $\|v+w\| \leq \|v\|+ \|w\|$
\end{itemize}
for all $c \in k$ and $v,w \in V$. An Archimedean normed space is an Archimedean semi-normed space such that $\|v\|=0$ if and only if $v=0$.
\end{defn}
\begin{defn}\label{defn:nonArchimedeanNormedSpace}
A non-Archimedean semi-normed space over a non-Archimedean field $k$ is a $k$-module $V$ together with a map $\| \|:V \to \mathbb{R}_{\geq 0}$ which satisfies 

\begin{itemize}
\item $\|cv\|=|c|\|v\|$
\item $\|v+w\| \leq \max\{\|v\|, \|w\|\}$
\end{itemize}
for all $c \in k$ and $v,w \in V$. A non-Archimedean normed space is a non-Archimedean semi-normed space such that $\|v\|=0$ if and only if $v=0$.
\end{defn}
\begin{rem}\label{rem:amazing}
If $V$ is a non-Archimedean normed space over a non-Archimedean field $k$ and $v,w\in V$ then if $\|v\|\neq \|w\|$ then $\|v+w\|=\max\{\|v\|,\|w\|\}.$
\end{rem}
\begin{rem}When we speak of a normed space over a valuation field $k$ we are talking about one of the above situations depending on whether $k$ is Archimedean or non-Archimedean.
\end{rem}
\begin{defn}
The category of Archimedean (resp. non-Archimedean) semi-normed spaces is defined as the category with the objects given by Archimedean (resp. non-Archimedean) normed spaces and the morphisms given by those $k$-linear maps $f:V \to W$ for which there exists a constant $C \in \mathbb{R}$ such that 
$\|f(v)\| \leq C\|v\|$ for all $v \in V$. We denote this category $\ttSNrm_{k}$. 
\end{defn}
\begin{rem}
Notice that a subspace $V$ of a semi-normed space $W$ inherits a natural semi-norm. If $V$ is a subspace of $W$ then the quotient vector space $W/V$ can be considered a semi-normed space when equipped with the semi-norm $\|[w]\|= \inf_{v\in V} \|w-v\|$. The quotient map $\pi:W \to W/V$ is bounded and in fact $\|\pi \| \leq 1$.
\end{rem}

\begin{defn}Let $V$ be a semi-normed space over a valuation field $k$. A subset $S \subset V$ is called bounded if $\sup_{s,t\in S} \|s-t\|< \infty$.
\end{defn}
\begin{defn}\label{defn:NormOfMap}
If $V, W \in \ttSNrm_{k}$, we can define a map 
$\Hom_{k}(V,W) \to \mathbb{R}_{+}$ which sends $T$ to $\|T\|$, the map is defined by 
\begin{equation}\label{equation:OperatorNorm}\|T\|=\inf \{C \in \mathbb{R} \ \ | \ \ 
\|Tv\| \leq C \|v\| 
\ \ \forall v \in V 
\}.
\end{equation}
\end{defn}
\begin{defn}\label{defn:nonexpanding}A morphism $f:V \to W$ in $\ttSNrm_{k}$ is called contracting if $\|f\| \leq 1$.
\end{defn}

\begin{lem} For any valuation field, the category $\ttSNrm_{k}$ is quasi-abelian. 
\end{lem}
{\bf Proof.} We simply note that the proof in \cite{SchneidersQA} works fine in this more (Archimedean or non-Archimdean) general  context. 
\ \hfill $\Box$

\subsection{The non-expanding semi-normed category}
Let $k$ be a valuation field.
\begin{defn}The category $\ttSNrm^{\leq 1}_{k}$ 
is defined to have the same objects as $\ttSNrm_{k}$.
The morphisms are the linear maps with norm less than or equal to one (they are contracting).
\end{defn}In the Archimedean case the product $\prod_{i\in I} V_{i}$ of a collection $\{V_{i}\}_{i \in I}$ in $\ttSNrm^{\leq 1}_{k}$ is given by 
\[\{
(v_i)_{i \in I} \in \times_{i \in I}V_{i} \ \ | \ \ \sup_{i \in I} \|v_{i}\| < \infty \}  
\] equipped with the norm 
\[\|(v_i)_{i \in I} \| =\sup_{i \in I} \|v_{i}\|.
\]
In the Archimedean case the coproduct $\coprod_{i\in I} V_{i}$ of a collection $\{V_{i}\}_{i \in I}$ in $\ttSNrm^{\leq 1}_{k}$ is given by \[\bigoplus_{i \in I} V_{i}
\] equipped with the norm 
\[\|(v_i)_{i \in I} \| = \sum_{i \in I} \|v_{i}\|.
\]
In the non-Archimedean case, the product $\prod_{i\in I} V_{i}$ of a collection $V_{i}$ in $\ttSNrm^{\leq 1}_{k}$ is given by a subspace of the product of vector spaces
\[\{
(v_i)_{i \in I} \in \times_{i \in I}V_{i} \ \ | \ \ \sup_{i \in I} \|v_{i}\| < \infty \}  
\] equipped with the semi-norm 
\[\|(v_i)_{i \in I} \| =\sup_{i \in I} \|v_{i}\|.
\]
In the non-Archimedean case, the coproduct $\coprod_{i\in I} V_{i}$ of a collection $V_{i}$ in $\ttSNrm^{\leq 1}_{k}$ is given by 
\[\bigoplus_{i \in I} V_{i}
\]
equipped with the semi-norm 
\[\|(v_i)_{i \in I} \| = \sup_{i \in I} \|v_{i}\|.
\]

Kernels and cokernels in $\ttSNrm^{\leq 1}_{k}$ are the same as those in $\ttSNrm_{k}$  which will be described in Lemma \ref{lem:SNrmLimColim}. The category $\ttSNrm^{\leq 1}_{k}$ has all limits and colimits.

\subsection{Morphisms and the closed structure}

\begin{defn}\label{defn:IntHomSNrm}
The structure defined in \ref{defn:NormOfMap} gives functors
\[\ttSNrm_{k}^{op} \times \ttSNrm_{k} \to \ttSNrm_{k}
\]
and 
\[(\ttSNrm^{\leq 1}_{k})^{op} \times \ttSNrm^{\leq 1}_{k} \to \ttSNrm^{\leq 1}_{k}
\]
which will be denoted
\[(V,W) \mapsto \underline{\Hom}_{k}(V,W).
\]
\end{defn}
\begin{defn}\label{defn:TensSnrm}
The symmetric monoidal structure which we will use assigns to two objects $V,W \in \ttSNrm_{k}$ or $\ttSNrm_{k}^{\leq 1}$ their \emph{projective} tensor product.  In the Archimedean case it is given by the algebraic tensor product $V\otimes_{k}W$ equipped the semi-norm 
\begin{equation}\label{equation:TensArch}
\|u\|=\inf\{\sum_{i=1}^{n}\|v_{i}\|\|w_{i}\|\ \  | \ \ u=\sum_{i=1}^{n} v_{i}\otimes w_{i}\}.
\end{equation}
In the non-Archimedean case it is given by the algebraic tensor product $V\otimes_{k}W$ equipped with the semi-norm
\begin{equation}\label{equation:TensNArch}\|u\|=\inf\{\max \{\|v_{i}\|\|w_{i}\| \ \ | \ \ i=1 , \dots ,  n \ \ \} \ \  | \ \ u=\sum_{i=1}^{n} v_{i}\otimes w_{i}\}
.\end{equation}
 It will be denoted simply $V \otimes_{k} W$ 
 and it defines bi-functors 
 \[\ttSNrm_{k} \times \ttSNrm_{k} \to \ttSNrm_{k}
 \]
 and 
  \[\ttSNrm^{\leq 1}_{k} \times \ttSNrm^{\leq 1}_{k} \to \ttSNrm^{\leq 1}_{k}
 \]
 which are exact in each variable.
 \end{defn}
 The following result is appears in many forms in the literature but we include it anyway to make sure it works when the field has trivial valuation.
\begin{lem}\label{lem:SNrmAdj}Let $k$ be a valuation field and let $U,V,W$ be semi-normed spaces over $k$. The natural equivalence of functors 
\[\ttVect^{op}_{k} \times \ttVect^{op}_{k} \times \ttVect_{k} \to \ttVect_{k}
\]
given by
\[\Hom_{k}(U \otimes_{k} V, W) \cong \Hom_{k}(U, \Hom_{k}(V,W)) 
\]
induces a natural equivalence of functors
\[\ttSNrm^{op}_{k} \times \ttSNrm^{op}_{k} \times \ttSNrm_{k} \to \ttSNrm_{k}
\]
given by morphisms of norm $1$. 
\begin{equation}\label{equation:show}\uHom_{k}(U \otimes_{k} V, W) \cong \uHom_{k}(U, \uHom_{k}(V,W)) 
\end{equation}
where $U \otimes_{k} V$ was defined in Definition \ref{defn:TensSnrm}. Therefore, 
\[\ttSNrm_{k}(U \otimes_{k} V, W) \cong \ttSNrm_{k}(U, \uHom_{k}(V,W))\] 
and
\[\ttSNrm^{\leq 1}_{k}(U \otimes_{k} V, W) \cong \ttSNrm^{\leq 1}_{k}(U, \uHom_{k}(V,W))\] 
showing that $U \mapsto U \otimes_{k} V$ is left adjoint to $W \mapsto \uHom_{k}(V,W)$ in $\ttSNrm_{k}$ and $\ttSNrm^{\leq 1}_{k}$.
\end{lem}
{\bf Proof.}
We check this in the non-Archimedean case. It will be obvious how to adapt it to the Archimedean case. Consider corresponding vectors $\phi \in \Hom_{k}(U,\Hom_{k}(V,W))$ and $\psi \in \Hom_{k}(U\otimes_{k} V, W)$. This means that for any $u\in U$ and $v\in V$ that $\psi(u\otimes v)=\phi(u)(v)$.

If $\psi$ is bounded and so can be considered as an element of $\uHom_{k}(U \otimes_{k} V, W)$ we have for any $u \in U$
\begin{equation}
\begin{split}\|\phi(u) \| & = \inf \{C \in \mathbb{R} \ \ | \ \ 
\|\psi(u \otimes v)\| \leq C \|v\| 
\ \ \forall v \in V \} \\
& \leq \inf \{C \in \mathbb{R} \ \ | \ \ 
\|\psi \| \|u \otimes v\| \leq C \|v\| 
\ \ \forall v \in V \} \\
&=\|\psi\| \inf \{C \in \mathbb{R} \ \ | \ \ 
 \|u \otimes v\| \leq C \|v\| 
\ \ \forall v \in V \} \\
&\leq \|\psi\| \inf \{C \in \mathbb{R} \ \ | \ \ 
 \|u \| \| v\| \leq C \|v\| 
\ \ \forall v \in V \} \\
& = \|\psi\| \|u\|
\end{split}
\end{equation}
and so $\phi(u)$ is bounded for all elements $u$ of $U$ and also $\phi$ is bounded so $\phi \in \uHom_{k}(U, \uHom_{k}(V,W))$ and $\|\phi\| \leq \|\psi\|$. In the other direction if $\phi(u)$ is bounded for all $u \in U$ and $\phi \in \uHom_{k}(U, \uHom_{k}(V,W))$ then for any $y \in U \otimes_{k} V$ and for all real $\epsilon >0$ there is a collection $u_1,v_1,u_2,v_2, \dots, u_n,v_n$ such that $y=\sum_{i=1}^{n}u_{i}\otimes v_{i}$ and $\|y\|+ \epsilon \leq \max_{i=1, \dots, n}\|u_i\|\|v_i\|$. Then we have
\begin{equation}
\begin{split}
\|\psi(y)\| & \leq \|\sum_{i=1}^{n}\phi(u_i)(v_i)\|  \leq \max_{i=1, \dots, n}\|\phi(u_i)(v_i)\| \\
&
\leq \max_{i=1, \dots, n}\|\phi(u_i)\|\|v_i\| \leq \max_{i=1, \dots, n}\|\phi\|\|u_i\|\|v_i\| \\
& \leq  \|\phi\| (\|y\|+\epsilon).
\end{split}
\end{equation}
Since $\epsilon$ was arbitrary, we conclude that $\psi$ is bounded so can be considered as an element of $\uHom_{k}(U \otimes_{k} V, W)$ and $\|\phi\| \leq \|\psi\|$. Thus we have shown that the natural bijection induces an isomorphism as in \ref{equation:show} and it is an isometry in the sense that $\|\phi\| = \|\psi\|$. In the Archimedean case, the only difference is that $\max_{i=1}^{n}$ is replaced by $\sum_{i=1}^{n}$ in the four places it appears in the proof.
\ \hfill $\Box$
\begin{lem}\label{lem:SNrmLimColim} For any valuation field $k$, the category $\ttSNrm_{k}$ has all finite limits and finite colimits.
\end{lem}
{\bf Proof.}
First of all, it is easy to see that given a finite collection of semi-normed spaces $V_1, \dots, V_n$ with norms $\| \|_1, \dots, \| \|_n$, the product $\prod_{i=1}^{n}V_{i}$ as vector spaces equipped with the norm $\|(v_1,\dots, v_n) \| =\max_{i=1}^{n}\|v_i \|_i$ is a product in $\ttSNrm_{k}$. In the Archimedean case, the sum $\bigoplus_{i=1}^{n}V_{i}$ of vector spaces equipped with the norm $\|(v_1,\dots, v_n) \| =\sum_{i=1}^{n}\|v_i \|_i$  is a co-product in $\ttSNrm_{k}$.  In the non-Archimedean case, the sum $\bigoplus_{i=1}^{n}V_{i}$ of vector spaces equipped with the norm $\|(v_1,\dots, v_n) \| =\max_{i=1}^{n}\|v_i \|_i$  is a co-product in $\ttSNrm_{k}$. The category of semi-normed spaces over $k$ has kernels: the kernel of a morphism $f:V \to W$ in $\ttSNrm_{k}$ is just the kernel in $\ttVect_{k}$ equipped with the restriction of $\| \|_{V}$. The category of semi-normed spaces over $k$ also has has cokernels: The cokernel of the morphism $f$ is just $V/\text{im}(f)$. By combining these operations, it is easy to form the limit of a finite diagram of semi-normed spaces as an appropriate kernel of the product of objects and the colimit of a finite diagram of semi-normed spaces as an  appropriate cokernel of the direct sum of objects in a finite category. 
\ \hfill $\Box$

\begin{lem}\label{lem:ProdCoprod}
For any finite collection of semi-normed spaces $V_1 , \dots , V_n$ the natural morphism $\bigoplus_{i=1}^{n}V_{i} \to \prod_{i=1}^{n}V_{i}$ is an isomorphism in $\ttSNrm_{k}$ when these spaces are equipped with the semi-norms described in Lemma \ref{lem:SNrmLimColim}.
\end{lem}
{\bf Proof.}
In the non-Archimedean case this is obvious and in the Archimedean case, it follows immediately from the fact that 
\[\max_{i=1 , \dots, n} \|\|_{i} \leq \sum_{i=1}^{n} \|\|_{i} \leq n \max_{i=1 , \dots, n} \|\|_{i}
\]
and so the identifications of the underlying vector spaces are bounded in both directions.
\ \hfill $\Box$
\begin{defn}\label{defn:Complete} A semi-normed space $V$ is complete if every Cauchy sequence in $V$ converges.
\end{defn}

\subsection{Banach Spaces}\label{BanachSpaces}
\begin{defn}
The category of Archimedean (resp. non-Archimedean) Banach spaces is defined as the full subcategory of category of of Archimedean (resp. non-Archimedean) normed spaces whose objects are complete (defined in Definition \ref{defn:Complete}). We use $\ttBan_{k}$ to denote this category.
\end{defn}

\noindent
The case of a trivially valued field presents some interesting behavior, as we show in the following example.

\begin{example}
Consider a field $k$ equipped with the trivial valuation. Consider the vector space $k[[t]]$ over $k$. Let $r$ be a real number satisfying $0 <r<1$. We can consider two non-Archimedean Banach spaces over $k$: given by $(k[[t]],\|\|_{triv})$ and $(k[[t]],\|\|_{adic})$ where $\|\|_{triv}$ is the trivial norm and $\|\|_{adic}$ is the norm which sends a series $\sum_{i=0}^{\infty} a_{i}t^{n_i}$ with $a_i \neq 0$ and $n_i\geq 0$ strictly increasing integers to $r^{n_0}$. Then the identity map 
\[(k[[t]],\|\|_{triv}) \to (k[[t]],\|\|_{adic})
\] is bounded but these spaces have no isomorphism between them due to the fact that in the adic space there is a non-zero sequence $t, t^{2}, t^{3}, \dots$ converging to zero and in the trivial valuation case there is no such sequence.  The normed space $k[t]$ when equipped with the adic norm is not complete, whereas it is clearly complete when equipped with the trivial norm.
\end{example}
\begin{rem}If $V \in \ttBan_{k}$ and $W$ is a closed subspace of $V$ then the induced semi-norm on $W$ is complete and the semi-norm on the quotient $V/W$ is complete. 

\end{rem}
\begin{lem}\label{lem:BanLimColim} For any valuation field $k$, the category $\ttBan_{k}$ has all finite limits and finite colimits. As in Lemma \ref{lem:ProdCoprod}, finite products and finite coproducts agree.
\end{lem}
{\bf Proof.}
The finite products and finite coproducts in $\ttBan_{k}$ are inherited automatically from those in $\ttSNrm_{k}$ which were discussed in the proof of Lemma \ref{lem:SNrmLimColim}.  The category of Banach spaces over $k$ has kernels: the kernel of a morphism $f:V \to W$ in $\ttBan_{k}$ is just the kernel (which is always closed) in $\ttVect_{k}$ equipped with the restriction of $\| \|_{V}$. The cokernel of the morphism $f$ is just $W/\overline{\text{im}(f)}$. By combining these operations, it is easy to form the limit as an appropriate kernel of the product of objects in a finite category and the colimit as an appropriate cokernel of the direct sum of objects in a finite category. 

\ \hfill $\Box$

\begin{defn}Suppose that $k$ is a non-Archimedean valuation field and $V \in \ttBan_{k}$. Given a collection of elements $v_{i} \in V$ indexed by some set $I$, a  limit of them is an element $v \in V$ such that for all $\epsilon >0$ there is a finite subset $J \subset I$ such that if $i$ is not in $J$ then $\|v_{i} - v\| < \epsilon$. If it exists then it is unique and written as $\lim_{i\in I}v_{i}$. The notation $\sum_{i \in I} v_{i}$ refers to a limit of the elements $\sum_{i \in J} v_{i}$ where $J$ runs over all finite subsets of $I$.
\end{defn}
The following lemma is well known
\begin{lem}\label{lem:sumsConverge}Suppose that $k$ is a non-Archimedean valuation field and $V \in \ttBan_{k}$. A sum $\sum_{i \in I} v_{i}$ of elements $v_{i} \in V$ converges if and only if $\lim_{i\in I}v_{i}=0$, i.e. if and only if $\lim_{i\in I}\|v_{i}\|=0$.
\end{lem}

\ \hfill $\Box$

\begin{lem}\label{lem:NoProdNoCoprod} The category $\ttBan_{k}$ has no infinite product or coproduct of any collection of non-zero objects. 
\end{lem}
{\bf Proof.}
Suppose that $C$ is the coproduct in $\ttBan_{k}$ of the infinite collection $V_i$ for $i \in I$ where no $V_i$ is $0$. Let $\tilde{C}$ be the coproduct of the $V_i$ in $\ttBan^{\leq 1}_{k}$. Consider a collection of morphisms $\text{id}_{i} \in \ttBan_{k}(V_i,\tilde{V_i})$ where $\tilde{V_{i}}$ is the vector space $V_{i}$ with its norm multiplied by $n_i>0$ where $n_i$ take arbitrarily high values. Then notice that $\|\text{id}_{i}\|$ is unbounded. Notice that there is a canonical morphism in $\ttBan_{k}$ from $C$ to $\tilde{C}$ commuting with the morphisms from the $V_i$ to both $C$ and $\tilde{C}$. By construction of $\tilde{C}$, the morphisms from $V_{i}$ to $\tilde{C}$ have norm $1$. Let $\tilde{\tilde{C}}$ be the coproduct of the $\tilde{V_{i}}$ in the category $\ttBan^{\leq 1}_{k}$. This results in the commutative diagram in $\ttBan_{k}$
\[\xymatrix{C  \ar[r] & \tilde{C} \\
\tilde{\tilde{C}} & V_{i}. \ar[u]\ar[ul] \ar[l]}
\]
Because the upper triangle is commutative, the norms of the diagonal maps must be bounded as we range over $i \in I$. By the universal property of $C$, there is a morphism $C \to \tilde{\tilde{C}}$ making the diagram commute. However this is a contradiction because the lower horizontal morphism has norms that go to infinity as we range over $i \in I$, while the norms of the diagonal maps remain bounded. The proof in the case of infinite products is similar.
\ \hfill $\Box$
\begin{rem}\label{rem:nofiltcolims}
Suppose that some filtered colimit exists in $\ttBan_{k}$. By passing to the (closure of the) images of the terms in this colimit it turns out that one necessarily has only gets finitely many distinct sub-objects of the colimit this way. In this case, the original filtered diagram can be replaced by a finite subcategory with finitely many objects and morphisms.
\end{rem}
The coimage of a morphism need not be isomorphic to the image.
\begin{defn}\label{defn:equivalentNorms} Two norms $\|\|$ and $\|\|'$ on a vector space are equivalent if there exist real numbers $C,C' >0$ such that $C \|\| \leq \|\|' \leq C' \|\|.$
\end{defn}

The following facts are easily proven the category $\ttBan_{k}$ in complete generality following along the lines of the proofs in \cite{Pr}. The last item uses Banach's open mapping theorem which  is valid only in the non-trivially valued setting and can be found in Chapter 1 section 3 of \cite{TVS}.
\begin{lem}\label{lem:BanProps}
For any morphism $u:E \to F$ in $\ttBan_{k}$ 
\begin{enumerate}
\item $u$ is continuous (Conversely, continuous linear map between Banach spaces will be bounded as long as the valuation is non-trivial)
\item $u^{-1}(0)$ is closed
\item $\ker (u) \cong u^{-1}(0)$ with the induced norm from $E$
\item $ \coker (u) \cong F/\overline{u(E)}$ with the quotient norm
\item $\im  (u) \cong \overline{u(E)}$ with the induced norm from $F$
\item $\coim (u) \cong u(E)$ where $\|f\|_{u(E)} = \inf_{u(e)=f} \|e\|_{E}$ 
\item $u$ is a monomorphism if and only if it injective
\item $u$ is an epimorphism if and only if $u(E)$ is dense in $F$
\item $u$ is strict if and only if $u(E)$ is closed and there exists a constant $C>0$ such that 
\[\inf_{e \in u^{-1}f}\|e\| \leq C\|f\|.
\]
for all $f \in u(E)$.
\item $u$ is a strict monomorphism if and only if there is a constant $C>0$ such that 
\[\|e\| \leq C \|u(e)\|
\]
for all $e \in E$.
\item If $u$ is a strict epimorphism then $u$ is surjective.
\item If the valuation on $k$ is non-trivial then $u$ is strict if and only if $u(E)$ is closed which in turn happens if and only if there exists a constant $C>0$ such that 
\[\inf_{e \in u^{-1}f}\|e\| \leq C\|f\|.
\]
for all $f \in u(E)$. Therefore, if the valuation on $k$ is non-trivial, $u$ is a strict epimorphism if and only if $u$ is surjective.
\end{enumerate}
\end{lem}
The following is based on Proposition 3.1.7 of \cite{Pr}.
\begin{lem} For any valuation field, the category $\ttBan_{k}$ is quasi-abelian. 
\end{lem}
{\bf Proof.} The proof in \cite{Pr} works fine in the general Archimedean setting, so we simply check that it can be adapted to the non-Archimedean setting.  The additivity holds just because finite products and coproducts are clearly isomorphic. Suppose that $u$ is a strict epimorphism in the cartesian diagram \begin{equation*}
\xymatrix{  E' \ar[r]^{u'} \ar[d]^{v'} & F' \ar[d]^v \\
E \ar[r]_{u}& F}
\end{equation*}
in $\ttBan_{k}$. Then $E' = \ker [(u, -v): E \oplus F' \to F].$ For any $(e,f') \in E'$ with $u'(e,f')= f'$ we have $\|f'\| \leq \|(e,f')\|$ and so 
\[\inf \{\| (e,f')\|  |u'(e,f')= f'\} \geq \| f'\|
\]
For any $\epsilon>0$ and any $f' \in F',$ using that $u$ is strict epimorphism, we can pick $e \in E$ with $u(e) = v(f')$ and $\|e\|<  C_{u}\|u(e)\| + \epsilon.$
Then since $\|v(f')\| \leq \|v\|\|f'\|$ we have 
\[ \| (e,f')\| = \max\{\|e\|,\|f'\| \} \leq \max \{ C_{u}\|u(e)\| + \epsilon, \|f'\|\} \leq \max \{ C_{u}\|v\|\|f'\| + \epsilon, \|f'\|.\} 
\]
Since $\epsilon$ was arbitrary, we conclude that
\[\inf \{\| (e,f')\|  |u'(e,f')= f'\} \leq  \max \{C_{u}\|v\|, 1\} \| f'\|.
\] 
Because $u$ is surjective, $u'$ is surjective as well and so $u'(E')$ is closed and so we can conclude that $u'$ is a strict epimpormism.  Suppose that $u$ is a strict monomorphism in the co-cartesian diagram \begin{equation*}
\xymatrix{  E' \ar[r]^{u'}  & F'  \\
E \ar[u]^{v} \ar[r]_{u}& F \ar[u]^{v'}}
\end{equation*}
in $\ttBan_{k}$. Then $\|e\| \leq C_{u}\|u(e)\| \leq C_{u} \|(v,-u)(e)\|$ and therefore $(v,-u)(E)$ is closed and so we have $F'= (E' \oplus F)/((v,-u)E).$ Fix $f' \in F'$ and suppose that $u'(e') = f'$.
In the quotient, we have 
\[\|(e',0)\| = \inf \{\|(e'',f) \| | f=-u(e), e''-e'=v(e)\}.
\]
Also, if $f=-u(e), e''-e'=v(e)$ we have $\|(e'',f) \|=\max\{\|e''\|,\|f\|\}$ and so
\begin{equation}\begin{split}\max \{C_{u}\|v\|,1\} \|(e'',f) \| & \geq \max( \|e''\|, C_{u}\|v\|\|f\|)=\max( \|e'+v(e)\|, C_{u}\|v\|\|u(e)\|) \\
& \geq  \max(\|e'+v(e)\|, \|v\|\|e\|)  \geq \max( \|e'+v(e)\|, \|v(e)\|) \geq \|e'\|.
\end{split}
\end{equation}
Therefore, $\max \{C_{u}\|v\|,1\}\|(e',0)\| \geq \|e'\|.$ 
Therefore, $u'$ is a strict monomorphism.
\ \hfill $\Box$

\begin{obs}\label{obs:IntHomBan}
The category $\ttBan_{k}$ is closed: the internal Hom functor defined on semi-normed spaces over $k$ in Definition \ref{defn:IntHomSNrm} preserves the property of being complete. Therefore for Banach spaces, the internal Hom functor between them is defined by treating them as semi-normed spaces.
\end{obs}
\begin{defn}\label{defn:ProjTensBan}
The symmetric monoidal structure which use on $\ttBan_{k}$ assigns to two objects $V,W \in \ttBan_{k}$ their \emph{projective} tensor product.  In the Archimedean case it is given (see Chapter 2 of \cite{Ryan}) by the completion of the algebraic tensor product $V\otimes_{k}W$ with respect to the norm 
\[\|u\|=\inf\{\sum_{i=1}^{n}\|v_{i}\|\|w_{i}\|\ \  | \ \ u=\sum_{i=1}^{n} v_{i}\otimes w_{i}\}.
\]
In the non-Archimedean case it is given (see Section 3.2 of \cite{Gruson}) by the  completion of the algebraic tensor product $V\otimes_{k}W$ with respect to the 
norm
\[\|u\|=\inf\{\max \{\|v_{i}\|\|w_{i}\| \ \ i=1 \dots n \ \ \} \ \  | \ \ u=\sum_{i=1}^{n} v_{i}\otimes w_{i}\}
.\]
 It will be denoted $V \wotimes_{k} W.$ 
  \end{defn}  
The completed tensor product is the bi-functor 
 \[\ttBan_{k} \times \ttBan_{k} \to \ttBan_{k}
 \]
 given on objects by 
 \[(V,W) \mapsto V \wotimes_{k} W.
 \] 

\begin{defn} 
 The internal Hom in this category 
\[\ttBan^{op}_{k} \times \ttBan_{k} \to \ttBan_{k}.
\] 
is denoted by $\uHom_{k}(V,W)$ and given by the Banach space whose underlying vector space is 
\[\{T \in \Hom_{k}(V,W)| \|T\| < \infty \}
\]
with norm given by  $\|T\|=\sup_{v \in V, v\neq0}\frac{\|T(v)\|}{\|v\|}$. We write $V'$ for $\uHom(V,k) \in \ttBan_{k}$.
\end{defn}

\begin{lem}\label{lem:TensStrictEpi} Let $k$ be a valuation field. For any $W \in \ttBan_{k}$ and any strict epimorphism $p:E \to V$ in $\ttBan_{k}$, the induced morphism  $\tilde{p}:E \wotimes_{k}W \to V \wotimes_{k}W$ is a strict epimorphism in $\ttBan_{k}$.
\end{lem}
{\bf Proof.} We prove this first in the non-Archimedean setting and then comment on how to alter the proof in the Archimedean case. Because $p$ is surjective, it is clear that $\tilde{p}(E \wotimes_{k} W)$ is dense, so we only need to show strictness. Using 1.1.9, Corollary 6 of \cite{BGR} it is enough to show that the morphism $E \otimes_{k} W \to V \otimes_{k}W$ is strict, where we use the same definition of strictness in the Banach and semi-normed settings. Here, the tensor products $E \otimes_{k} W $ and $V \otimes_{k}W$ carry the semi-norms from Equation (\ref{equation:TensNArch}).  Choose $C$ such that $\inf_{e \in p^{-1}(v)}\|e\| \leq C\|v\|$ for all $v \in V$. Fix $f\in  V \otimes_{k}W$ and $\epsilon>0$. Suppose that $f= \sum_{i=1}^{n}v_{i}\otimes w_{i}$. Using surjectivity of $p$ and the strictness of $p$, choose $e_{i}\in E$ such that $\|e_{i}\|  < C\|v_{i}\|+ \frac{\epsilon}{\|w_i\|}$ and $p(e_i)=v_i$. Then 
and  $\tilde{p}(\sum_{i=1}^{n}e_{i}\otimes w_{i}) = f$ and

\[\|\sum_{i=1}^{n}e_{i}\otimes w_{i} \|  \leq \max_{i=1,\dots, n}\|e_i\|\|w_i\|  < \max_{i=1,\dots, n}C\|v_i\|\|w_i\| + \epsilon.
\]
Therefore, 
\[\inf_{\tilde{e} \in \tilde{p}^{-1}(f)}\|\tilde{e}\| < \max_{i=1, \dots, n}C\|v_i\|\|w_i\|+\epsilon.
\]
Since $\epsilon$ was arbitrary, we conclude that 
\[\inf_{\tilde{e} \in \tilde{p}^{-1}(f)}\|\tilde{e}\| \leq \max_{i=1, \dots, n}C\|v_i\|\|w_i\|.
\]
Because this holds for all presentations $f= \sum_{i=1}^{n}v_{i}\otimes w_{i}$ we can conclude that 
\[\inf_{\tilde{e} \in \tilde{p}^{-1}(f)}\|\tilde{e}\| \leq C\|f\|.
\]
The only difference in the proof for the Archimedean setting is that we use the tensor product from Equation (\ref{equation:TensArch}) and $\sum_{i=1}^{n}$ replaces $\max_{i=1, \dots, n}$.
\ \hfill $\Box$
\begin{lem}\label{lem:HomStrictMon} Let $k$ be a valuation field. For any $W \in \ttBan_{k}$ and any strict monomorphism $j:V \to E$ in $\ttBan_{k}$, the induced morphism  $\tilde{j}:\uHom_{k}(W,V) \to \uHom_{k}(W,E)$ is a strict monomorphism in $\ttBan_{k}$.
\end{lem}
{\bf Proof.}
Because $j$ is a strict monomorphism we can choose $C$ such that $\|v\| \leq C \|j(v)\|$ for all $v\in V$. Then for any $w \in W$ where $w \neq 0$ and any $\phi \in \uHom_{k}(W,V)$ we have \[\frac{\|\phi(w)\|}{\|w\|} \leq C\frac{\|j(\phi(w))\|}{\|w\|}.\] By passing to the supremum over all $w \in W$ where $w \neq 0$ we conclude that \[\|\phi\| \leq C \|j \circ \phi \|= C\| \tilde{j}(\phi)\|\] for all $\phi \in \uHom_{k}(W,V)$. Therefore, $\tilde{j}$ is a strict monomorphism.
\ \hfill $\Box$
\subsection{The non-expanding Banach category}\label{NexBan}
Let $k$ be a valuation field. In this subsection we consider the category of Banach spaces over $k$ with non-expanding morphisms. We do not consider geometry relative to this category in the present article, but rather will use the good properties of this category to construct and characterize objects in the category $\ttBan_{k}$.
\begin{defn}The category $\ttBan^{\leq 1}_{k}$ 
is defined to have the same objects as $\ttBan_{k}$.
The morphisms are the linear maps with norm less than or equal to one (they are contracting).
\end{defn}
Product and coproducts in $\ttBan^{\leq 1}_{k}$ exist. In the Archimedean case (see page 63 of \cite{Hel}) the product $\prod_{i\in I} V_{i}$ of a collection $\{V_{i}\}_{i \in I}$ in $\ttBan^{\leq 1}_{k}$ is given by 
\[\{
(v_i)_{i \in I} \in \times_{i \in I}V_{i} \ \ | \ \ \sup_{i \in I} \|v_{i}\| < \infty \}  
\] equipped with the norm 
\[\|(v_i)_{i \in I} \| =\sup_{i \in I} \|v_{i}\|
\]
while the coproduct $\coprod_{i\in I} V_{i}$ of a collection $\{V_{i}\}_{i \in I}$ in $\ttBan^{\leq 1}_{k}$ is given by \[\{
(v_i)_{i \in I} \in \times_{i \in I}V_{i} \ \ | \ \sum_{i \in I} \|v_{i}\| <\infty \}  
\] equipped with the norm 
\[\|(v_i)_{i \in I} \| = \sum_{i \in I} \|v_{i}\|.
\]
In the non-Archimedean case they can be computed as in \cite{Gruson}: the product $\prod_{i\in I} V_{i}$ of a collection $\{V_{i}\}_{i \in I}$ in $\ttBan^{\leq 1}_{k}$ is given by 
\[\{
(v_i)_{i \in I} \in \times_{i \in I}V_{i} \ \ | \ \ \sup_{i \in I} \|v_{i}\| < \infty \}  
\] equipped with the norm 
\[\|(v_i)_{i \in I} \| =\sup_{i \in I} \|v_{i}\|
\]
while the coproduct $\coprod_{i\in I} V_{i}$ of a collection $\{V_{i}\}_{i \in I}$ in $\ttBan^{\leq 1}_{k}$ is given by \[\{
(v_i)_{i \in I} \in \times_{i \in I}V_{i} \ \ | \ \lim_{i \in I} \|v_{i}\| =0 \}  
\] equipped with the norm 
\[\|(v_i)_{i \in I} \| = \sup_{i \in I} \|v_{i}\|.
\]
Suppose we are given a collection $\{f_{i}:V_{i} \to W_{i}\}_{i \in I}$ in $\ttBan^{\leq 1}_{k}.$ Then observe that the natural morphism
\[\coprod_{i\in I} \ker(f_i) \rightarrow \ker [\coprod_{i\in I} V_{i} \rightarrow  \coprod_{i\in I} W_{i}] 
\]
is an isomorphism. Similarly, if $V_{i} \subset V$ and $W_{i} \subset W$ are an increasing union of subspaces with union $V$ and $W$ respectively then the natural map 
\[\cup_{i\in I} \ker(f_i) \rightarrow \ker [V \rightarrow W] 
\]
is an isomorphism. The closed symmetric monoidal structure on  $\ttBan^{\leq 1}_{k}$ is defined in precisely the same way as in Observation \ref{obs:IntHomBan} and Definition \ref{defn:ProjTensBan}.
Suppose now that $k$ is non-Archimedean. Then $\ttBan^{\leq 1}_{k}$ is an additive category. 
Kernels and cokernels in $\ttBan^{\leq 1}_{k}$ are the same as those in $\ttBan_{k}$ which were described in Lemma \ref{lem:BanLimColim}. They commute. The category $\ttBan^{\leq 1}_{k}$ has all limits and colimits and the product of a finite collection agrees with the coproduct of a finite collection. Although we don't use this, if $k$ is non-Archimedean, the category $\ttBan^{\leq 1}_{k}$ is quasi-abelian.

\subsection{Enough projectives and injectives}

The proof in \cite{Pr} that the quasi-abelian category $\ttBan_{k}$ has enough projectives does not work in the case where the valuation is trivial. We give a more general proof to make sure that the category $\ttBan_{k}$ always has enough projectives.
\begin{defn}Suppose that $k$ is non-Archimedean and let $V \in \ttBan_{k}$. Define
\[c_{0}(V)= \{(c_{v})_{v \in V-\{0\}} | c_{v} \in k, \lim_{v \in V-\{0\}} \|c_{v}v\|=0 \}
\]
equipped with the norm 
\[\|(c_{v})_{v \in V-\{0\}}\| = \sup_{v \in V-\{0\}} \| c_{v} v\|.
\]
This is a Banach space because it is the coproduct of the collection $k_{\|v\|}$ over all $v \in V-\{0\}$ in the category $\ttBan^{\leq 1}_{k}$.
\end{defn}
\begin{lem}Suppose that $k$ is a non-Archimedean valuation field. Suppose that $V \in \ttBan_{k}$.  Then $c_{0}(V)$ is projective in $\ttBan_{k}$.
\end{lem}
{\bf Proof.}

Suppose that $\mu:E \to F$ is a strict epimorphism, and consider a bounded linear map $\nu: c_{0}(V)\to F.$ Define $1_{w}=(\delta_{v,w})_{v \in V-\{0\}} \in c_{0}(V)$. For each $v \in V-\{0\}$ we have by Lemma \ref{lem:BanProps} 
\[\inf_{\mu(e)=\nu(1_{v})} \|e\|_{E} \leq C_{\mu}\|\nu(1_{v})\|_{F} \leq C_{\mu} \|\nu\| \|1_{v}\|_{c_{0}(V)}= C_{\mu} \|\nu\| \|v\|
\]
Choose $\epsilon>0$. For each $v \in V-\{0\}$ we can use the surjectivity of $\mu$ from Lemma \ref{lem:BanProps} to choose $e_{v} \in E$ such that $\mu(e_{v})=\nu(1_{v})$ and 
\[\|e_{v}\|_{E} < C_{\mu} \|\nu\| \|v\|+\|v\|= (C_{\mu} \|\nu\|+\epsilon)\|v\|.
\]

Notice that 
\[\lim_{v \in V-\{0\}} \|c_{v}e_{v}\|= \lim_{v \in V-\{0\}} |c_{v}| \|e_{v}\| \leq \lim_{v \in V-\{0\}}(C_{\mu} \|\nu\|+\epsilon)|c_{v}|\|v\| = (C_{\mu} \|\nu\|+\epsilon)\lim_{v \in V-\{0\}}\|c_{v}v\| = 0.
\]
Therefore by Lemma \ref{lem:sumsConverge}, we can form the sum $\sum_{v \in V-\{0\}}c_{v}e_{v}.$ Notice that 
\[\|\sum_{v \in V-\{0\}}c_{v}e_{v}\| \leq \sup_{v \in V-\{0\}}\|c_{v}e_{v}\|\leq \sup_{v \in V-\{0\}}(C_{\mu} \|\nu\|+\epsilon)\|c_{v}v\| = (C_{\mu} \|\nu\|+\epsilon)\|(c_{v})_{v \in V-\{0\}}\|.
\]
This means that the linear map 
\[\nu':c_{0}(V)\to E
\]
given by
\[(c_{v})_{v \in V-\{0\}} \mapsto \sum_{v \in V-\{0\}}c_{v}e_{v}
\]
is bounded. It is easy to see that $\mu \circ \nu'=\nu$. Therefore, $c_{0}(V)$ is projective in $\ttBan_{k}$. 
\ \hfill $\Box$
\begin{lem}\label{lem:BanEnoughProj}There is a strict epimorphism $c_{0}(V)\to V$ in $\ttBan_{k}$. For any valuation field, the category $\ttBan_{k}$ 
has 
enough projectives. 
\end{lem}

{\bf Proof.}
For the general Archimedean case one can use a variation on the proof over $\mathbb{C}$ in \cite{Pr}. Suppose first that $k$ is a 
non-Archimedean field and let $V \in \ttBan_{k}$.

Notice that for every element $(c_{v})_{v \in V-\{0\}} \in c_{0}(V)$ there is a well defined element \[\sum_{v \in V-\{0\}} c_{v} v \in V\] because $\lim_{v \in V-\{0\}}\|c_{v} v\| =0$. Using Lemma \ref{lem:sumsConverge} we may define a linear morphism 
\begin{equation}\label{equation:EnoughSur}
\kappa_{V}:c_{0}(V)\to V
\end{equation}
given by
\[(c_{v})_{v \in V-\{0\}} \mapsto \sum_{v \in V-\{0\}}c_{v}v.
\]
It is bounded and in fact $\|\kappa_{V}\| \leq 1$ because 
\[\|\sum_{v \in V-\{0\}}c_{v}v\| \leq \sup_{v \in V-\{0\}} \|c_{v}v\| = \|(c_{v})_{v \in V-\{0\}}\|. 
\]
Also, it is surjective because $\kappa_{V}(0)=0$ and for any $w \in V-\{0\},$ $\kappa_{V}(1_{w})=w$.  In fact it is a strict epimorphism in $\ttBan_{k}$ because by the above if $\kappa_{V}((c_{v})_{v \in V-\{0\}})=w$ then $\|w\| \leq  \|(c_{v})_{v \in V-\{0\}}\|$ and so $\|w\| \leq \inf_{\kappa_{V}((c_{v})_{v \in V-\{0\}})=w}\|(c_{v})_{v \in V-\{0\}}\|$ and also 
\[\inf_{\kappa_{V}((c_{v})_{v \in V-\{0\}})=w}\|(c_{v})_{v \in V-\{0\}}\| \leq \|1_{w}\| = \|w\|
\]
for all $w \in V-\{0\}$ and of course $\inf_{\kappa_{V}((c_{v})_{v \in V-\{0\}})=0}\|(c_{v})_{v \in V-\{0\}}\| = 0$.


\ \hfill $\Box$

\begin{defn}\label{defn:ellInfty}
Suppose that $k$ is non-Archimedean valuation field and $V \in \ttBan_{k}$.

Define
\[\ell_{\infty}(V)= \{(c_{v})_{v \in V-\{0\}} | c_{v} \in k,  \sup_{v \in V-\{0\}} \frac{|c_v|}{\|  v\|} < \infty \}
\]
equipped with the norm 
\[\|(c_{v})_{v \in V-\{0\}}\| = \sup_{v \in V-\{0\}}  \frac{|c_v|}{\|  v\|} .
\]
This is a Banach space because it is the product of the collection $k_{\|v\|^{-1}}$ over all $v \in V-\{0\}$ in the category $\ttBan^{\leq 1}_{k}$.  For any morphism $f:V\to W$ in $\ttBan_{k}$, let \[\ell_{\infty}(f):\ell_{\infty}(V') \to \ell_{\infty}(W')\] be the morphism $(f(c_{\alpha})_{\alpha \in V'-\{0\}})_{\beta}=c_{\beta \circ f}.$ 
\end{defn}

\begin{lem}\label{lem:InjObj} Suppose that $k$ is non-Archimedean with a non-trivial valuation and also that $k$ is spherically complete. Suppose that $V \in \ttBan_{k}$.  Then $\ell_{\infty}(V)$ is injective in the quasi-abelian category $\ttBan_{k}$.
\end{lem}
{\bf Proof.} 
Suppose that $\mu:E \to F$ is a strict monomorphism in $\ttBan_{k}$ and $\nu:E \to \ell_{\infty}(V)$ is any morphism in $\ttBan_{k}$. Define, for each $v\in V-\{0\},$
$\nu_{v}$ as the composition $E \to \ell_{\infty}(V) \to k$ where we project onto the coordinate labeled by $v$.  We have
\[|\nu_{v}(e)|\leq \|v\| \|\nu(e)\| \leq \|v\| \|\nu\| \|e\|
\]
and therefore the linear map $\nu_{v}$ is bounded with norm less than or equal to $\|\nu\|\|v\|$. Using the Hahn-Banach theorem (valid in this spherically complete context by \cite{Ing}) we can extend each $\nu_{v}$ to a morphism $\nu'_{v}:F \to k$ without increasing the norm. We have $|\nu_{v}'(f)| \leq \|\nu\| \|v\| \|f\|$ for any $f\in F$. Finally, the sought after factorization is provided by 
\[\nu':F \to \ell_{\infty}(V)
\]
defined by 
\[f \mapsto (\nu_{v}'(f))_{v \in V-\{0\}}.
\]
Notice that this is bounded and in fact $\|\nu'\| \leq \|\nu\|$ because for any $f\in F$ we have
\[\| (\nu_{v}'(f))_{v \in V-\{0\}} \| =  \sup_{v \in V-\{0\}} \frac{ | \nu_{v}'(f) |}{\|  v\| }\leq \|\nu\| \|f\|.
\]
Therefore, $\ell_{\infty}(V)$ is injective in $\ttBan_{k}$.
\ \hfill $\Box$

\begin{lem}\label{lem:BanEnoughInj}Let $k$ be a valuation field. The quasi-abelian category $\ttBan_{k}$ has enough injectives.
\end{lem}
As far as the first statement is concerned, for the general Archimedean case or the non-trivially valued non-Archimedean case one can use a variation on the proof over $\mathbb{C}$ in \cite{Pr}. So we assume that $k$ is non-Archimedean. 

{\bf Proof.}
First, assume that $k$ is non-trivially valued and spherically complete. Consider the morphism 
\[ u:V  \to \ell_{\infty}(V')
\]
defined by 
\[u(v)_{\alpha}= \alpha(v)
\] for each $\alpha \in V'-\{0\}$.
It is a strict monomorphism in $\ttBan_{k}$ because if $(c_{\alpha})_{\alpha \in V'-\{0\}}$ is the image of $v$ then
\[
\|(c_{\alpha})_{\alpha \in V'-\{0\}}\| = \sup_{\alpha \in V'-\{0\}} \frac{| c_{\alpha} |}{\|\alpha\|}=  \sup_{\alpha \in V'-\{0\}} \frac{| \alpha(v)|}{\| \alpha\|} = \|v\|.
\]
At the last step we have used the pseudo-reflexivity of $V$ (see Section 11 of \cite{Schneider} for this property and its validity in the spherically complete context).
Now, suppose that $k$ is a general non-Archimedean field and $V\in \ttBan_{k}.$ There is a strict monomorphism over $k$ of $V$ into the completed tensor product of $V$ over $k$ with a certain non-trivially valued field $K_r$ introduced in Proposition 2.1.2 of \cite{Ber1990}. If $k$ was already non-trivially valued, one could simply use $k$ itself instead of $K_r$. The field $K_r$ has the property that the completed tensor product with it preserves strict monomorphisms.  Notice also that a strict monomorphism in $\ttBan_{K_r}$ remains a strict monomorphism when considered as a morphism in $\ttBan_{k}$. Therefore by Lemma \ref{lem:BaseChangeInj} applied to the functor given given by the completed tensor product with $K_r$
\[\ttBan_{k} \to \ttBan_{K_r}
\]
and its right adjoint forgetful functor, it is enough to produce a strict monomorphism of the resulting element of $\ttBan_{K_r}$ into an injective of $\ttBan_{K_r}$. Therefore, we may assume that $k$ is non-trivially valued.
Let $K$ be a spherical completion of the fraction field of $\Sym_{\leq 1}(V)$. We are using the notation $\Sym_{\leq 1}(V)$ to refer to the symmetric algebra introduced in subsection \ref{Symmetric} computed in the category $\ttBan^{\leq 1}_{k}$. See 3.2.2 in \cite{Rob} for more information on the spherical completion of a field. The construction in \cite{Rob} refers to the case where the field is the algebraic closure of $\mathbb{Q}_{p}$ but the construction works for any non-Archimedean field. Because the norm on  $\Sym_{\leq 1}(V)$ is strictly multiplicative (as can be seen from the Gauss Lemma as in Proposition 1 of 5.1.2 in \cite{BGR}) the norm on $\Frac(\Sym_{\leq 1}(V))$ is as well.  Therefore $\Frac(\Sym_{\leq 1}(V))$ is a non-Archimedean valuation field. Properties of the spherical completion imply that the norm on $\Frac(\Sym_{\leq 1}(V))$ agrees with the norm on its image (a closed subspace of) $K$ and that $K$ is spherically complete. The morphism $V \to K$ is a strict monomorphism in $\ttBan_{k}$ because it decomposes into strict monomorphisms in $\ttBan_{k}$
\[V \to \Sym_{\leq 1}(V) \to \Frac(\Sym_{\leq 1}(V)) \to K.
\] 
Consider the composition $u$ given by
\begin{equation}\label{equation:consider} V \to  K \to \ell_{\infty}(K').
\end{equation}
Using Lemmas \ref{lem:BaseChangeInj} and \ref{lem:InjObj} we see that $\ell_{\infty}(K')$ is injective in $\ttBan_{k}.$ Since $K$ is both non-trivially valued and spherically complete, the morphism $K \to \ell_{\infty}(K')$ is a strict monomorphism in $\ttBan_{K}$. Therefore, (\ref{equation:consider}) is a strict monomorphism in $\ttBan_{k}$.

\ \hfill $\Box$

\begin{defn}\label{defn:explInj} Let $\sI: \ttBan_{k} \to  \ttBan_{k}$ be the functor constructed in the proof of Lemma \ref{lem:BanEnoughInj} together with the functorial construction in Definition \ref{defn:ellInfty}. This is also extended to a strictly exact sequence 
\[0 \to M \to \sI(M)_{1} \to \sI(M)_{2} \to  \sI(M)_{3} \to \cdots
\]
where $\sI(M)_{0} = M$, $\sI(M)_{1} = \sI(M)$ and $\sI(M)_{i+1} = \sI(\sI(M)_{i}/\sI(M)_{i-1})$ for $i \geq 1$.
\end{defn}
\begin{rem}Much of the gymnastics involved in the proof of Lemma \ref{lem:BanEnoughInj} can be avoided if we use the fact that for any non-Archimedean valuation field $k$ there is a spherically complete field $K$ containing $k$ so that the completed tensor product with $K$ over $k$ preserves strict exact sequences and any such sequence embeds isometrically into its completed tensor product with $K$.  See \cite{Gruson} and Lemma 3.1 of \cite{Poi}. 
\end{rem}
\subsection{The closed structure in the category of Banach spaces}
\begin{lem}Let $k$ be a valuation field and let $U,V,W$ be Banach spaces over $k$. The natural equivalence of functors 
\[\ttVect^{op}_{k} \times \ttVect^{op}_{k} \times \ttVect_{k} \to \ttVect_{k}
\]
given by
\[\Hom_{k}(U \otimes_{k} V, W) \cong \Hom_{k}(U, \Hom_{k}(V,W)) 
\]
induces a natural equivalence of functors
\[\ttBan^{op}_{k} \times \ttBan^{op}_{k} \times \ttBan_{k} \to \ttBan_{k}
\]
given by morphisms of norm $1$
\[\uHom_{k}(U \wotimes_{k} V, W) \cong \uHom_{k}(U, \uHom_{k}(V,W)). 
\]
Therefore, 
\[\ttBan_{k}(U \wotimes_{k} V, W) \cong \ttBan_{k}(U, \uHom_{k}(V,W)),
\]
and
\[\ttBan^{\leq 1}_{k}(U \wotimes_{k} V, W) \cong \ttBan^{\leq 1}_{k}(U, \uHom_{k}(V,W)),
\]
showing that $U \mapsto U \wotimes_{k} V$ is left adjoint to $W \mapsto \uHom_{k}(V,W)$ in $\ttBan_{k}$ and $\ttBan^{\leq 1}_{k}$.
\end{lem}
{\bf Proof.} The proof goes along the lines of the corresponding statement for semi-normed spaces which was Lemma \ref{lem:SNrmAdj}. The only difference is that because we are now using the tensor product of Banach spaces defined in Definition \ref{defn:ProjTensBan}, we use the canonical bijection between bounded maps from $U \otimes_{k} V$ to $W$ in the category of semi-normed spaces and morphisms in the category of Banach spaces between $U \wotimes_{k} V$ and $W$.
\ \hfill $\Box$
\begin{obs}
Notice that the completion functor defined in Definition \ref{defn:Completion} is a morphism of closed symmetric monoidal categories from the category of semi-normed spaces to the category of Banach spaces. This morphism is essentially surjective.
\end{obs}
In the category $\ttBan_{k}$ it is easily checked that any finite dimensional object is isomorphic to $k_{\infty}^{\oplus n}$ for some positive integer $n$, this indicates the vector space $k^{\oplus n}$ equipped with the norm $\|(v_1, \dots, v_n) \| = \max_{i=1}^{n} |v_i|.$ However, in the category $\ttBan^{\leq 1}_{k}$ this is not the case. 
\begin{defn}\label{defn:OneDim}
For each $r \in \mathbb{R}_{+}$, define the one dimensional Banach spaces $k_{r}$ over $k$ to be $k$ equipped with the norm $c \mapsto r |c|$.
\end{defn}
 We have $k_{r_1}\wotimes_{k} k_{r_2} \cong k_{r_1 r_2}$. For any Banach space $W$ over $k$ we have 
\begin{equation}\label{equation:HomBall}\ttBan^{\leq 1}_{k}(k_s,W) \cong \{w \in W | \|w\| \leq s\}
\end{equation} where the morphism on the left is determined by sending $1$ to $w$. Notice that if $r_1, r_2 \in \mathbb{R}_{>0}$ then the one dimensional Banach spaces $k_{r_i}$ over $k$ are isomorphic in $\ttBan^{\leq 1}_{k}$ if and only if $\frac{r_1}{r_2} \in |k^{\times}|$.

\subsection{Completion}
\begin{defn}\label{defn:Completion} Let $k$ be a complete valuation field. When we speak of the categories $\ttSNrm^{\leq 1}_{k}$ or $\ttBan^{\leq 1}_{k}$ assume that we are in the non-Archimedean setting.
There is a completion functor 
\[\straightC_{k}:\ttSNrm_{k} \to \ttBan_{k}
\]
defined in the standard way by taking equivalence classes of Cauchy sequences. We sometimes denote the completion of a morphism $\phi:V \to W$ by $\widehat{\phi}:\widehat{V} \to \widehat{W}$.
\end{defn}
\begin{rem}\label{rem:CompletionColim}
The completion functor $\straightC_k$ from Definition \ref{defn:Completion} is left adjoint to the forgetful functor $\ttBan_{k} \to \ttSNrm_{k}$. It commutes with colimits. Note that the morphism $V \to \straightC_{k}(V)$ is the colimit of the category whose objects are pairs $(f,W)$ where $W \in \ttBan_{k}$ and $f \in \ttSNrm_{k}(V,W)$ and whose morphisms are commuting triangles. By Proposition 5, Chapter 1 of \cite{BGR} the completion functor commutes with kernels of admissible morphisms as well. The underlying topological space of $\straightC_{k}(V)$ is the closure of the image of $V$ in $\straightC_{k}(V)$.
\end{rem}
The completion functor is not fully faithful but it respects the other structures we have introduced
\begin{lem}\label{lem:CompletionCompat}For any $M,N \in \ttSNrm_{k}$ there are natural isomorphisms 
\[\straightC_{k}{\uHom_{k}(M,N)} \cong \uHom_{k}(\straightC_{k}(M), \straightC_{k}(N))
\]
and
\[\straightC_{k}(M \otimes N) \cong \straightC_{k}(M) \wotimes_{k} \straightC_{k}(N).
\]
\end{lem}
\ \hfill $\Box$

There is a commutative diagram of adjoint pairs
\begin{equation}\label{equation:FreeAdj}
\xymatrix{
\ttVect_{k} 
\ar@/^/[rr]^{}  
& & \ttComm(\ttVect_{k}) 
\ar@/^/[ll]_{}
\\
\ttSNrm_{k} 
\ar@/^/[rr]^{}  \ar[d] \ar[u]
& & \ttComm(\ttSNrm_{k}) \ar[d] \ar[u]
\ar@/^/[ll]_{}
\\
\ttBan_{k} 
\ar@/^/[rr]^{}  
& & \ttComm(\ttBan_{k}) 
\ar@/^/[ll]_{}
}
\end{equation}
where the top arrow in each row takes a Banach space to the (completion in the Banach case) of the symmetric algebra over it and the bottom arrow is a forgetful functor. The arrows going up are also forgetful functors and the arrows going down are completion functors (defined in Definition \ref{defn:Completion}).
\subsection{Banach algebras and modules} Let $k$ be a complete, valued field. \begin{rem}
The objects of $\ttComm(\ttBan_{k})$ differ from the standard definition of Banach $k$-algebras in that for each algebra object $\mA$ there is a universal constant 
$C$ such that $\|ab\| \leq C \|a\|\|b\|$ for all $a,b \in \mA$
but we do not insist that $C=1$. 
\end{rem}In this subsection, we discuss the categories $\ttMod(\mA)$ (and $\ttMod^{\leq 1}(\mA)$ in the non-Archimedean case) for any  $\mA \in \ttComm(\ttBan_{k})$. \begin{defn}For $k$ a non-Archimedean complete valued field, the category $\ttMod^{\leq 1}(\mA)$ has the same objects as $\ttMod(\mA)$ but its morphisms are those which have norm less than or equal to one.
\end{defn}

Let us sumerize some of their properties.

\begin{lem}\label{lem:AdjBan}
The categories $\ttMod(\mA)$ and  $\ttMod^{\leq 1}(\mA)$ are closed symmetric monoidal category having finite limits and colimits. The closed structures will be denoted $\uHom_{\mA}.$ The monoidal structures will be denoted $\widehat{\otimes}_{\mA}$. Note that for any $\mE,\mF,\mG \in \ttMod(\mA)$ we have the adjunction isomorphisms
\[\uHom_{\mA}(\mE \widehat{\otimes}_{\mA} \mF, \mG) \cong \uHom_{\mA}(\mE, \uHom_{\mA}(\mF,\mG))
\]
which behave as in Lemma \ref{lem:SNrmAdj}.
In both of these categories, the functor of tensoring with $\mF$ is therefore left adjoint to the functor of mapping from $\mF$ and therefore the former is right exact and the later is left exact.  In all the cases where we have said that limits and colimits in $\ttMod(\mA)$ or $\Mod^{c}(\mA)$ exist, they can be computed simply by equipping the models we have given for the corresponding limits and colimits in $\ttBan_{k}$ or $\ttBan^{\leq 1}_{k}$ with the obvious $\mA$ module structures. The analogue of Lemma \ref{lem:BanProps} holds in  $\ttMod(\mA)$. 
\end{lem}
{\bf Proof.} Define 

\[\uHom_{\mA}(\mF,\mG) = \{\phi \in \uHom(\mF,\mG) \ \ | \ \ \phi(af) = a\phi(f) \ \ \forall a \in \mA , f \in \mF\}.
\]
The rest of these properties are easily checked.
\ \hfill $\Box$
\begin{lem}The functor $\ootimes_{\mA}$ defined in Definition \ref{defn:TensA} agrees with the standard notion $\wotimes_{\mA}$ in the literature. So the tensor product $E\ootimes_{\mA}F$ in $\ttMod(\mA)$ or $\ttMod^{\leq 1}(\mA)$ is isomorphic in the non-Archimedean case to be the completion of the Banach space $E \otimes_{\mA}F$ with respect to the semi-norm 
\[\|u\|=\inf\{\max \{\|e_{i}\|\|f_{i}\| \ \ i=1 , \dots ,  n \ \ \} \ \  | \ \ u=\sum_{i=1}^{n} e_{i}\otimes f_{i}\}
\]
along with the action induced by $(a,e\otimes f) \mapsto ae \otimes f = e \otimes af$. In the Archimedean case it is 
\[\|u\|=\inf\{\sum_{i=1}^{n} \|e_{i}\|\|f_{i}\|  \ \  | \ \ u=\sum_{i=1}^{n} e_{i}\otimes f_{i}\}.
\]
\end{lem}
{\bf Proof.} The claim follows immediately from (1) Remark \ref{rem:CompletionColim}, (2) the fact that 
\begin{equation}\label{eqn:TensDef2}E \otimes_{\mA} F = \colim [\xymatrix{E \otimes_{k} \mA \otimes_{k} F \ar@/^/[rr]^{} \ar@/_/[rr]_{} & & E \otimes_{k} F}]
\end{equation}
where the completion of the diagram $\xymatrix{E \otimes_{k} \mA \otimes_{k} F \ar@/^/[rr]^{} \ar@/_/[rr]_{} & & E \otimes_{k} F}$ is the diagram in Equation \ref{equation:TensDef} which defined $E \ootimes_{\mA} F$ as a colimit similarly to Equation \ref{eqn:TensDef2}
and (3) the fact that taking colimits in different orders leads to naturally isomorphic answers.
\ \hfill $\Box$

\begin{rem}\label{rem:properties}

In the non-Archimedean setting, the products and coproducts in the category $\ttMod^{\leq 1}(\mA)$ have interesting properties. Consider for instance the coproduct of some elements $V_{i}$ in the category $\ttMod^{\leq 1}(\mA)$. Conceretly, it is the Banach space  \[\{
(v_i)_{i \in I} \in \prod_{i \in I}V_{i} \ \ | \ \lim_{i \in I} \|v_{i}\| =0 \}  
\] equipped with the norm 
\[\|(v_i)_{i \in I} \| = \sup_{i \in I} \|v_{i}\|.
\] and the obvious action of $\mA$. For any other object $W$, we have natural isomorphisms (see Proposition 8, page 76 of \cite{BGR})
\begin{equation}\label{eqn:coprodtens}(\coprod_{i\in I}V_{i})\wotimes_{\mA} W \cong \coprod_{i\in I}(V_{i} \wotimes_{\mA} W)
\end{equation}
and for any families $\{V_{i}\}_{i \in I}$ and $\{W_{i}\}_{i \in I}$ we have natural isomorphisms
\[\coprod_{i\in I}(V_{i} \times W_{i}) \cong (\coprod_{i\in I}V_{i}) \times (\coprod_{i\in I} W_{i}) 
\]
because finite products and coproducts are isomorphic. The inclusions $V_{i} \to \coprod_{i\in I}V_{i}$ are strict monomorphisms and if the $V_{i}$ are $\wotimes_{\mA}$-acyclic then so is $\coprod_{i\in I}V_{i}$.
\end{rem}

\begin{thm}Let $k$ be a valuation field. Choose $\mA \in \ttComm(\ttBan_{k})$. The quasi-abelian category $\Mod(\mA)$ has enough $\wotimes_{\mA}$-acyclices, projectives and injectives. \end{thm}
{\bf Proof.} In order to see the statements about projectives there are two ways to go. The first is simply to consider that the proof that $\ttBan_{k}$ has enough projectives and every object has a contracting strict epimorphism from a projective object. All the constructions involved in that proof actually can be carried out directly in the categories $\Mod(\mA)$  more or less replacing $k$ by $\mA$. For the other proof one can just tensor all the constructions with $\mA$. Therefore, it is an immediate combination of Lemmas \ref{lem:HelpWithEnough}, \ref{lem:BanEnoughProj} and Lemma \ref{lem:TensStrictEpi}. 
We now consider the statements about injectives. They follow immediately from Lemmas \ref{lem:HelpWithEnough}, \ref{lem:BanEnoughInj} and \ref{lem:HomStrictMon}.
\ \hfill $\Box$
\begin{defn}\label{defn:preferred}
We define a preferred resolution in the closed symmetric monoidal quasi-abelian category $\ttMod(\mA)$ of any object $\mM$. Let
\[\sP(\mM) = \coprod_{m\in \mM-\{0\}}\mA_{\|m\|}
\]
where the coproduct is computed in the category $\ttMod^{\leq 1}(\mA).$ Note that $\sP(\mM)$ is both projective and $\wotimes_{\mA}$-acyclic in the closed symmetric monoidal quasi-abelian category $\ttMod(\mA)$. We have a strict epimorphism 
\begin{equation}\label{equation:strictepi}\phi:\sP(\mM) \to \mM
\end{equation}
in $\ttMod(\mA)$ defined by 
\[(a_{m})_{m \in \mM-\{0\}} \mapsto \sum_{m \in \mM-\{0\}} a_{m}m.
\]
We define a strictly exact sequence 
\[\cdots \to \sP(\mM)_{3} \to \sP(\mM)_{2} \to \sP(\mM)_{1} \to \mM \to 0
\]
with $\sP(\mM)_{0}=\mM$, $\sP(\mM)_{1}=\sP(\mM)$ and $\sP(\mM)_{i+1}=\sP(\ker[\sP(\mM)_{i} \to \sP(\mM)_{i-1}])$ for $i>1$.
\end{defn}

\begin{lem}For any $V \in \ttMod(\mA)$ the morphism 
\[\sP(\mM)\wotimes_{\mA}V\to \mM\wotimes_{\mA}V
\]
is a strict epimorphism.
\end{lem}

\ \hfill $\Box$

We will now introduce another family of objects which one could use for the computation of derived tensor products. Let $\mA \in \ttComm(\ttBan_{k})$ and $M \in \ttMod(\mA).$
Since the category $\ttMod(\mA)$ has enough projectives (any object has a strict epimorphism from a projective object), the functor $\sF(N) = N \wotimes_{\mA} M$ is left derivable. As projectives which are coproducts of $\mA$ are $\sF$-acyclic, and any object has a strict epimorphism from such a coproduct, we get from Lemma \ref{lem:acyclic} that the class of $\sF$-acyclics is $\sF$-projective.

\section{Category theory background}
Let $\ttC$ be a category with fiber products.  A full subcategory 
$\ttD \subset \ttC$ is called dense when the restriction of the Yoneda embedding
\[\ttC \to \ttPr(\ttD)= \Hom(\ttD^{op},\ttSet) 
\]
is fully faithful. This implies (see ExpI, Prop, 7.2 of \cite{SGA4}) that any $c \in \ttC$ is the colimit of the canonical diagram $\ttD/c \to \ttC$. 
Recall that a site with underlying category $\ttC$ consists for every $U \in \ttC$ of a collection $S_{U}$ of covering families $\{U_{i}\to U\}_{i\in I}$ in $\ttC$ including isomorphisms, closed under compositions and  pullbacks with respect to arbitrary morphisms in $\ttC$.  These are called covering families of $U$ and define a {\it Grothendieck pretopology} on $\ttC.$ We denote this site by $\ttC^{S}$. Let $\ttC^{S}$ and $\ttD^{T}$ be sites and assume that they are subcanonical, this means that representable presheaves are sheaves.
\begin{defn}\label{defn:ContCocont} A functor $\su:\ttC \to \ttD$ is called continuous (with respect to the set of covers $S$ and $T$) if the image of every covering family in $S$ is a covering family in $T$ and for any $W \to V$ in $\ttC$ the morphism 
\[\su(W \times_{V} V_{i}) \to \su(W) \times_{\su(V)} \su(V_i)
\]
is an isomorphism. A functor $\su:\ttC \to \ttD$ is called cocontinuous (with respect to the set of covers $S$ and $T$) if for every $c \in \ttC$ and any covering $\{d_{j}\to \su(c)\}_{j \in J}$ of $\ttD$ there exists a covering $\{c_{i} \to c\}_{i\in I}$ such that the family of maps $\{\su(c_i) \to \su(c)\}_{i\in I}$ refines the covering  $\{d_{j}\to \su(c)\}_{j \in J}$.
\end{defn}

 The category of presheaves $\ttPr(\ttC)$ is simply the category of functors from $\ttC^{op}$ to the category of sets. The category of sheaves with respect to a topology $S$ on $\ttC$ is be denoted  $\ttSh(\ttC^{S})$.

Say $\su:\ttC \to \ttD$ is any functor.  The pullback functor on presheaves will be written $\su^{-1}:\ttPr(\ttD) \to \ttPr(\ttC)$.  It has left and right adjoints $\su_{!}$ and $\su_{*}$ respectively.  It is shown in SGA4 \cite{SGA4} I.5.6 that the functor $\su_{!}$ is fully faithful if and only if $\su$ is fully faithful.

The category of sheaves with respect to the topology $S$ on $\ttC$ will be denoted  $\ttSh(\ttC^{S})$.
We use $\su_{s}$ to denote the composition $\ttSh(\ttC^{S}) \to \ttPr(\ttC) \stackrel{\su_{!}}\to \ttPr(\ttD) \to \ttSh(\ttD^{T})$ where the first map is the fully faithful inclusion of sheaves into presheaves and the final map is sheafification. It is known (see \cite{stacks5}) that if $\su$ is continuous that $\su^{-1}$ preserves sheaves and by Lemma 7.14.3 of \cite{stacks2} that $\su_{s}$ is a left adjoint to $\su^{-1}: \ttSh(\ttD^{T})\to  \ttSh(\ttC^{S}).$ From \cite{stacks6} and the cocontinuity we get a natural equivalence between $\su^{-1} \su_{s}$ and the identity and therefore we have the following 

\begin{lem}\label{lem:extension}Suppose we are given a continuous, cocontinuous, fully faithful functor $u: \ttC\to \ttD$ between the categories underlying sub-canonical sites $\ttC^{S}$ and $\ttD^{T}$. Then the functor $\su_{s}: \ttSh(\ttC^{S}) \to  \ttSh(\ttD^{T})$ is fully faithful.
\end{lem}

 \begin{rem}Notice that in the situation of Lemma \ref{lem:extension}, $\su_{s}$ commutes with colimits. Also by Lemma 7.14.5 of \cite{stacks3}, $\su_{s}$ preserves representable sheaves in the sense that $\su_{s}(h_{c})=h_{\su(c)}$ for every $c\in \ttC$. We have \[(\su^{-1}\su_{s}(h_{c}))(c') = \su_{s}(h_{c}) (\su(c')) =h_{\su(c)}(\su(c'))= \Hom(\su(c'), \su(c)) 
\]
for every $c,c' \in \ttC$. Therefore since $\su$ is fully faithful, we have natural isomorphisms $\su^{-1}\su_{s}(h_{c}) \cong h_{c}$ for every $c\in \ttC$. 
 \end{rem}

\bibliographystyle{amsalpha}

\end{document}